\documentclass[11pt]{amsart}
\usepackage{amsmath,amsthm,amssymb,bm}
\usepackage{xcolor}
\usepackage{indentfirst}
\usepackage{graphicx}
\usepackage{tikz}
\usetikzlibrary{decorations.pathreplacing, decorations.pathmorphing, decorations.shapes}
\usetikzlibrary{arrows,calc}
\usepackage{mathrsfs,esint}
\usepackage{comment}
\usepackage{enumerate}
\usepackage{hyperref}
\hypersetup{
    colorlinks = true,
    linkcolor = red,
    anchorcolor = green!50!black,
    citecolor = green!50!black,
    filecolor = green!50!black,
    urlcolor = blue
}
\usepackage{cleveref}

\usepackage[top=1in, bottom=1in, left=1.25in, right=1.25in]{geometry}
\usepackage{setspace}

\usepackage{float}

\usepackage{subcaption}
\usepackage{pgfplots}
\pgfplotsset{compat=newest, compat/show suggested version=false}
\usepackage{tikz}
\usetikzlibrary{intersections}

\usepackage{diagbox}
\usepackage{algorithm, algorithmicx, algpseudocode}
\usepackage{tabularx}

\usepackage[titletoc,toc,title]{appendix}

\newtheorem{thm}{Theorem}[section]
\newtheorem{theorem}{Theorem}[section]

\newtheorem{remark}[thm]{\textit{Remark}}

\newcommand{\titleShort}{Fast Solver for FFPE}
\newcommand{\titleLong}{A Fast and Accurate Solver for the Fractional Fokker-Planck Equation with Dirac-Delta Initial Conditions}

\newcommand{\HIGHLIGHTA}[1]{\textcolor{red}{#1}}
\newcommand{\drift}{b}
\newcommand{\diffusionOrdinary}{D_{\textnormal{o}}}
\newcommand{\diffusionFractional}{D_{\textnormal{f}}}

\newcommand{\integralExpP}{K^{+}}
\newcommand{\integralExpM}{K^{-}}
\newcommand{\integralExpPM}{K^{\pm}}

\newcommand{\integralSinExpP}{T^{+}}
\newcommand{\integralSinExpM}{T^{-}}
\newcommand{\integralSinExpPM}{T^{\pm}}

\title[\titleShort]{\titleLong}

\thanks{}

\author{Qihao~Ye}
\address{Department of Mathematics, University of California, San Diego, CA 92093, United States} 
\email{q8ye@ucsd.edu}

\author{Xiaochuan~Tian}
\address{Department of Mathematics, University of California, San Diego, CA 92093, United States} 
\email{xctian@ucsd.edu}

\author{Dong~Wang}
\address{School of Science and Engineering, The Chinese University of Hong Kong, Shenzhen, Guangdong 518172, China \& Shenzhen International Center for Industrial and Applied Mathematics, Shenzhen Research Institute of Big Data, Guangdong 518172, China \& Shenzhen Loop Area Institute, Guangdong 518172, China}
\email{wangdong@cuhk.edu.cn}

\date{}

\numberwithin{equation}{section}

\begin{document}

\newcommand{\beq}{\begin{equation}}
\newcommand{\eeq}{\end{equation}}
\newcommand{\sgn}{\operatorname{sgn}}
\newcommand{\dist}{\text{dist}}
\newcommand{\argmin}{\text{argmin}}

%%%%%%%%%% MATHCAL %%%%%%%%%%
\def\cA{\mathcal{A}}
\def\cB{\mathcal{B}}
\def\cC{\mathcal{C}}
\def\cD{\mathcal{D}}
\def\cE{\mathcal{E}}
\def\cF{\mathcal{F}}
\def\cG{\mathcal{G}}
\def\cT{\mathcal{T}}
\def\cH{\mathcal{H}}
\def\cI{\mathcal{I}}
\def\cJ{\mathcal{J}}
\def\cK{\mathcal{K}}
\def\cL{\mathcal{L}}
\def\cM{\mathcal{M}}
\def\cN{\mathcal{N}}
\def\cO{\mathcal{O}}
\def\cP{\mathcal{P}}
\def\cR{\mathcal{R}}
\def\cS{\mathcal{S}}
\def\cT{\mathcal{T}}
\def\cU{\mathcal{U}}
\def\cV{\mathcal{V}}
\def\cW{\mathcal{W}}
\def\cX{\mathcal{X}}
\def\cY{\mathcal{Y}}
\def\cZ{\mathcal{Z}}

%%%%%%%%%% MATHBB %%%%%%%%%%
\def\A{\mathbb{A}}
\def\B{\mathbb{B}}
\def\C{\mathbb{C}}
\def\D{\mathbb{D}}
\def\E{\mathbb{E}}
\def\F{\mathbb{F}}
\def\G{\mathbb{G}}
\def\H{\mathbb{H}}
\def\I{\mathbb{I}}
\def\J{\mathbb{J}}
\def\K{\mathbb{K}}
\def\L{\mathbb{L}}
\def\M{\mathbb{M}}
\def\N{\mathbb{N}}
\def\O{\mathbb{O}}
\def\P{\mathbb{P}}
\def\Q{\mathbb{Q}}
\def\R{\mathbb{R}}
\def\S{\mathbb{S}}
\def\T{\mathbb{T}}
\def\U{\mathbb{U}}
\def\V{\mathbb{V}}
\def\W{\mathbb{W}}
\def\X{\mathbb{X}}
\def\Y{\mathbb{Y}}
\def\Z{\mathbb{Z}}

%%%%%%%%%% MATHSCR %%%%%%%%%%
\def\sA{\mathscr{A}}
\def\sB{\mathscr{B}}
\def\sC{\mathscr{C}}
\def\sD{\mathscr{D}}
\def\sE{\mathscr{E}}
\def\sF{\mathscr{F}}
\def\sG{\mathscr{G}}
\def\sH{\mathscr{H}}
\def\sI{\mathscr{I}}
\def\sJ{\mathscr{J}}
\def\sK{\mathscr{K}}
\def\sL{\mathscr{L}}
\def\sM{\mathscr{M}}
\def\sN{\mathscr{N}}
\def\sO{\mathscr{O}}
\def\sP{\mathscr{P}}
\def\sQ{\mathscr{Q}}
\def\sR{\mathscr{R}}
\def\sS{\mathscr{S}}
\def\sT{\mathscr{T}}
\def\sU{\mathscr{U}}
\def\sV{\mathscr{V}}
\def\sW{\mathscr{W}}
\def\sX{\mathscr{X}}
\def\sY{\mathscr{Y}}
\def\sZ{\mathscr{Z}}

%%%%%%%%%% MATHFRAK %%%%%%%%%%
\def\fA{\mathfrak{A}}
\def\fB{\mathfrak{B}}
\def\fC{\mathfrak{C}}
\def\fD{\mathfrak{D}}
\def\fE{\mathfrak{E}}
\def\fF{\mathfrak{F}}
\def\fG{\mathfrak{G}}
\def\fH{\mathfrak{H}}
\def\fI{\mathfrak{I}}
\def\fJ{\mathfrak{J}}
\def\fK{\mathfrak{K}}
\def\fL{\mathfrak{L}}
\def\fM{\mathfrak{M}}
\def\fN{\mathfrak{N}}
\def\fO{\mathfrak{O}}
\def\fP{\mathfrak{P}}
\def\fQ{\mathfrak{Q}}
\def\fR{\mathfrak{R}}
\def\fS{\mathfrak{S}}
\def\fT{\mathfrak{T}}
\def\fU{\mathfrak{U}}
\def\fV{\mathfrak{V}}
\def\fW{\mathfrak{W}}
\def\fX{\mathfrak{X}}
\def\fY{\mathfrak{Y}}
\def\fZ{\mathfrak{Z}}

%%%%%%%%%%%%% BM %%%%%%%%%%%%%%%%
\def\ab{\bm a}
\def\bb{\bm b}
\def\cb{\bm c}
\def\db{\bm d}
\def\eb{\bm e}
\def\fb{\bm f}
\def\gb{\bm g}
\def\hb{\bm h}
\def\ib{\bm i}
\def\jb{\bm j}
\def\kb{\bm k}
\def\lb{\bm l}
\def\mb{\bm m}
\def\nb{\bm n}
\def\ob{\bm o}
\def\pb{\bm p}
\def\qb{\bm q}
\def\rb{\bm r}
\def\sb{\bm s}
\def\tb{\bm t}
\def\ub{\bm u}
\def\vb{\bm v}
\def\wb{\bm w}
\def\xb{\bm x}
\def\yb{\bm y}
\def\zb{\bm z}

\def\xib{\bm\xi}

%%%%%%%%%% Greek letters %%%%%%%%%

\newcommand{\al}{\alpha}
\newcommand{\be}{\beta}
\newcommand{\ga}{\gamma}
\newcommand{\del}{\delta}
\newcommand{\ep}{\epsilon}
\newcommand{\vep}{\varepsilon}
\newcommand{\ze}{\zeta}
\newcommand{\vthe}{\vartheta}
\newcommand{\io}{\iota}
\newcommand{\ka}{\kappa}
\newcommand{\la}{\lambda}
\newcommand{\vrho}{\varrho}
\newcommand{\sig}{\sigma}
\newcommand{\up}{\upsilon}
\newcommand{\vphi}{\varphi}
\newcommand{\om}{\omega}

\newcommand{\Ga}{{\Gamma}}
\newcommand{\Del}{\Delta}
\newcommand{\The}{\Theta}
\newcommand{\La}{\Lambda}
\newcommand{\Sig}{\Sigma}
\newcommand{\Up}{\Upsilon}
\newcommand{\Om}{\Omega}

\begin{abstract}
The classical Fokker-Planck equation (FPE) is a key tool in physics for describing systems influenced by drag forces and Gaussian noise, with applications spanning multiple fields.
We consider the fractional Fokker-Planck equation (FFPE), which models the time evolution of probability densities for systems driven by L\'evy processes, relevant in scenarios where Gaussian assumptions fail. 
The paper presents an efficient and accurate numerical approach for the free-space FFPE with constant coefficients and Dirac-delta initial conditions.
This method utilizes the integral representation of the solutions and enables the efficient handling of very high-dimensional problems using fast algorithms.
Our work is the first to present a high-precision numerical solver for the free-space FFPE with Dirac-delta initial conditions.
In addition to Dirac-delta initial data, we demonstrate the effectiveness of our method for initial conditions given by sums of Gaussians.
This opens the door for future research on more complex scenarios, including those with variable coefficients and other types of initial conditions.

\end{abstract}



\maketitle
\begin{center}
{\small \color{purple}
This work has been published in the \emph{SIAM Journal on Scientific Computing}.\\
Please refer to \href{https://doi.org/10.1137/24M1682907}{the official publication} for citation purposes.\\
The source code is available at \href{https://github.com/ACMathX/FFPEDDIC}{GitHub}.}
\end{center}
\onehalfspacing

%%%%%%%%%%%%%%%%%%%%%%%%%%%Section 1 %%%%%%%%%%%%%%%%%%%%%%%%%%

\section{Introduction}\label{sec:introduction}
The classical Fokker-Planck equation (FPE) is a fundamental equation in physics that describes the dynamic behavior of systems under the influence of drag forces and Gaussian random noise.
It is widely applied in various fields, including statistical mechanics, plasma physics, social science, economics, neuroscience, information theory, and data science \cite{faisal2008noise,helbing2001traffic,hull1990pricing,nicholson1983introduction,pavliotis2014stochastic,risken1996fokker,sato2014approximation,stephan2017stochastic,thomas2006elements}.
However, there is growing interest in models that account for non-Gaussian random noise characterized by L\'evy processes.
Such models have found applications in fields like physics, biology, and economics, highlighting their significance in accurately modeling real-world systems where standard Brownian motion falls short \cite{benson2000fractional,gao2022data,kao2010random,metzler2000random,reynolds2007displaced,scalas2000fractional,tesfay2021dynamics}.
Motivated by this, we consider the fractional Fokker-Planck equation (FFPE), which describes the time evolution of probability density functions for L\'evy processes.

As an initial step in this challenging area, we consider the following free-space FFPE with constant coefficients:
\begin{equation*}
    \frac{\partial p}{\partial t} ( \boldsymbol{x}, t )
    = - \boldsymbol{\drift} \cdot \nabla p ( \boldsymbol{x}, t )
    + \diffusionOrdinary \Delta p ( \boldsymbol{x}, t )
    - \diffusionFractional \left ( - \Delta \right )^{\alpha} p ( \boldsymbol{x}, t ),
\end{equation*}
where $p(\boldsymbol{x}, t)$ denotes the probability density function at position $\boldsymbol{x}\in \R^d $ at time $t>0$, and $\boldsymbol{\drift} \in \mathbb{R}^{d}$, $\diffusionOrdinary \geq 0$ and $\diffusionFractional > 0$ represent the drift coefficient, the ordinary diffusion coefficient, and fractional diffusion coefficient, respectively.
Note that the ordinary diffusion coefficient is allowed to vanish in this work, which is a scenario that poses particular challenges for the numerical evaluation of solutions.
Here $(- \Delta)^{\alpha}$ ($\alpha \in (0, 1)$) is the fractional Laplacian operator, which is the infinitesimal generator of the pure L\'evy $2 \alpha$-stable process \cite{applebaum2009levy,di2012hitchhikerʼs}. 
In our research, we focus on developing an efficient and accurate numerical method for solving the aforementioned FFPE with Dirac-delta initial conditions given by $ p(\boldsymbol{x}, 0) = \delta_{\boldsymbol{x}_{0}} ( \boldsymbol{x} )$ for a given $\boldsymbol{x}_{0} \in \mathbb{R}^{d}$.
This provides the fundamental solution of the FFPE, which can be used to solve other types of initial conditions.
Due to the nonlocal term, the fundamental solution of the FFPE is not explicitly known.
Therefore, an efficient and accurate numerical approach is crucial and could be highly beneficial in numerous studies.

The numerical approaches to the FPE have been the subject of extensive discussion.
Among the various methods, finite difference methods, finite element methods, and Lagrangian particle methods have been thoroughly investigated \cite{carrillo2019blob,degond1990deterministic,hu2023positivity,hu2022energetic,jordan1998variational,junge2017fully,srinivasan2018positivity,sun2018discontinuous}.
These discussions often center on the unique advantages each method offers such as structure preservation and energy dissipation properties.
More recently, there has been a surge in discussions regarding solvers for high-dimensional FPE \cite{chen2023convex,chertkov2021solution,hu2024score,liu2022neural,tang2023solving,zhai2022deep}, paralleling the increasing interest in the discussion of methods for high-dimensional PDEs.
For the discretization of equations involving the fractional Laplacian, a variety of methods have been developed.
These include quadrature-based methods \cite{duan2023quadrature,duo2018novel,huang2014numerical,tian2015class}, finite element methods \cite{acosta2019finite,acosta2017fractional,bonito2019numerical,burrage2012efficient,liu2021diagonal,sun2022numerical}, spectral methods \cite{bueno2014fourier,duo2018fast,zhou2024novel}, rational approximation techniques \cite{duan2023pade,harizanov2020analysis}, among others.
For a comprehensive introduction to fractional equations and the related numerical methods, readers can refer to the review papers \cite{bonito2018numerical,d2020numerical}.
These works focus on bounded domains with regular initial conditions in low-dimensional settings. Neural network approaches and Monte Carlo methods have also been explored \cite{pang2019fpinns, zan2020stochastic,zeng2023adaptive}.
However, neural network approaches have not demonstrated high accuracy and do not address Dirac-delta initial conditions, while Monte Carlo methods require a large number of sample paths.
Our work is the first to present a high-precision numerical solver for the free-space FFPE with Dirac-delta initial conditions, capable of handling very high-dimensional problems.

In this work, we start by representing the fundamental solution $p(\boldsymbol{x}, t)$ as a $d$-dimensional integral through Fourier transform, which further reduces to the one-dimensional integral
\begin{equation*}
    p(\boldsymbol{x}, t)
    = \frac{1}{y^{(d - 2) / 2}} \int_{0}^{\infty} \left ( \frac{r}{2 \pi} \right )^{d / 2} J_{(d - 2) / 2}(y r) \exp \left ( - \left ( \diffusionOrdinary r^{2} + \diffusionFractional r^{2 \alpha} \right ) t \right ) \, d r,
\end{equation*}
where $y := |\boldsymbol{x} - \boldsymbol{x}_{0} - \boldsymbol{\drift} t| > 0$ and $J_{\nu}$ is the Bessel function of the first kind.
Our method then focuses on the fast numerical evaluation of the integral in this form.
The fast evaluation of such integral relies on the careful treatment of the near-origin integrals using expansion techniques or re-weighting techniques and far-field integrals using the windowing function techniques \cite{bruno2017windowed,monro2008super}.
This is combined with other techniques, such as the scaling methods specific to the solutions to this problem.
Additionally, the case $y = 0$ is treated separately.
Our approach enables the efficient handling of very high-dimensional problems using fast algorithms.
Another major advantage of our method is that it does not require time stepping, thus eliminating the errors commonly associated with time-stepping schemes.
For moderate values of $t$, the evaluation of $p(\boldsymbol{x}, t)$ can achieve machine precision.
As $t > 0$ gets smaller, the evaluation of $p(\boldsymbol{x}, t)$ becomes more challenging as the solution approaches a Dirac-delta measure.
The accuracy and feasibility of our method for a given $t$ depend on the ordinary and fractional diffusion coefficients $\diffusionOrdinary$ and $\diffusionFractional$, the fractional diffusion power $\alpha$, and the dimension $d$ of the problem.
Comprehensive studies on these dependencies are presented in the numerical experiments in \Cref{sec:numerical_results}.
For further reading on fast algorithms for evaluating integrals, readers may refer to \cite{gao2022kernel,greengard2023dual,greengard1987fast,ying2004kernel}. 
While the main focus of this paper is on the FFPE with constant coefficients and Dirac-delta initial conditions, we also present numerical results for initial data given by sums of Gaussians, which serve as a natural adaptation of our method.
Such tests are of practical interest, as general functions can often be approximated by Gaussian radial basis functions \cite{Buhmann_2003,lowe1988multivariable,shao2025solving}.
Future work may consider other types of equations, such as those with variable coefficients, and other types of initial conditions, such as characteristic functions for the nonlocal threshold dynamics \cite{caffarelli2010convergence,jiang2018efficient}.

This paper is organized as follows.
In \Cref{sec:FFPE_Dirac_IC}, we derive the integral representation of the solution to the FFPE using Fourier analysis, along with the scaling laws, which form the basis of our numerical methods.
\Cref{sec:numerical_approach} details our numerical approach.
We begin by introducing techniques for handling singular and far-field integrals, followed by the introduction of our main solver in \Cref{sec:FFPE_solver}.
Our numerical approach is concisely outlined in \Cref{rmk:Dirac_concise_explanation}.
\Cref{sec:numerical_results} presents the numerical results, demonstrating the efficiency and accuracy of our method.
Finally, in \Cref{sec:conclusion}, we provide our conclusions and suggest avenues for further research.

%%%%%%%%%%%%%%%%%%%%%%%%%%%Section 2 %%%%%%%%%%%%%%%%%%%%%%%%%%

\section{FFPE with Dirac-delta initial conditions}\label{sec:FFPE_Dirac_IC}
We consider the FFPE in arbitrary dimensions:
\begin{equation}\label{eq:Dirac_general_FFPE_hd}
    \left \{
    \begin{aligned}
        \frac{\partial p}{\partial t} ( \boldsymbol{x}, t )
        &= - \boldsymbol{\drift} \cdot \nabla p ( \boldsymbol{x}, t )
        + \diffusionOrdinary \Delta p ( \boldsymbol{x}, t )
        - \diffusionFractional \left ( - \Delta \right )^{\alpha} p ( \boldsymbol{x}, t )\\
        p \left ( \boldsymbol{x}, 0 \right )
        &= \delta_{\boldsymbol{x}_{0}} ( \boldsymbol{x} )
    \end{aligned}
    \right .,
\end{equation}
where $d$ is the dimension, $\boldsymbol{\drift} \in \mathbb{R}^{d}$, $\diffusionOrdinary \in \mathbb{R}_{\geq 0}$, $\diffusionFractional \in \mathbb{R}_{\geq 0}$, $\diffusionOrdinary + \diffusionFractional > 0$, $\alpha \in (0, 1)$, $\boldsymbol{x}_{0} \in \mathbb{R}^{d}$, $\boldsymbol{x} \in \mathbb{R}^{d}$ and $t \in \mathbb{R}_{+}$.

\begin{remark}
    We focus on cases where diffusion does not vanish, i.e., $ \diffusionOrdinary + \diffusionFractional > 0$. 
    In the case where only ordinary diffusion is present, i.e., $\diffusionFractional = 0$, \Cref{eq:Dirac_general_FFPE_hd} is reduced to the classical FPE with the solution given by:
    \begin{equation*}\label{eq:Dirac_special_no_fractional_diffusion}
        p(\boldsymbol{x}, t)
        = \frac{1}{\left ( \sqrt{4 \pi \diffusionOrdinary t} \right )^{d}} \exp \left ( - \frac{|\boldsymbol{x} - \boldsymbol{x}_{0} - \boldsymbol{\drift} t|^{2}}{4 \diffusionOrdinary t} \right )
        = \frac{1}{\left ( \sqrt{4 \pi \diffusionOrdinary t} \right )^{d}} \exp \left ( - \frac{y^{2}}{4 \diffusionOrdinary t} \right ).
    \end{equation*}
    Therefore, our primary emphasis is addressing the scenarios where $\diffusionFractional > 0$.
\end{remark}
Here and for the rest of this paper we denote $y := | \boldsymbol{x} - \boldsymbol{x}_{0} - \boldsymbol{\drift} t |$.
In the rest of this section, we present the mathematical tools that underpin the numerical approach that will be discussed in \Cref{sec:numerical_approach}.

\subsection{Integral representation of the solution}
Apply the Fourier transform on both sides of \Cref{eq:Dirac_general_FFPE_hd}, we have
\begin{equation*}
    \left \{
    \begin{aligned}
        \hat{p}_{t}(\boldsymbol{k}, t)
        &= - i \langle \boldsymbol{k}, \boldsymbol{\drift} \rangle \hat{p}(\boldsymbol{k}, t)
        - \diffusionOrdinary |\boldsymbol{k}|^{2} \hat{p}(\boldsymbol{k}, t)
        - \diffusionFractional |\boldsymbol{k}|^{2 \alpha} \hat{p}(\boldsymbol{k}, t)\\
        \hat{p}(\boldsymbol{k}, 0)
        &= \exp \left ( - i \langle \boldsymbol{k}, \boldsymbol{x}_{0} \rangle \right )
    \end{aligned}
    \right .,
\end{equation*}
which implies that
\begin{equation*}
    \begin{aligned}
        \hat{p}(\boldsymbol{k}, t)
        &= \exp \left ( - i \langle \boldsymbol{k}, \boldsymbol{x}_{0} + \boldsymbol{\drift} t \rangle \right ) \exp \left ( - \left ( \diffusionOrdinary |\boldsymbol{k}|^{2} + \diffusionFractional |\boldsymbol{k}|^{2 \alpha} \right ) t \right ).
    \end{aligned}
\end{equation*}
Hence
\begin{equation}\label{eq:Dirac_explicit_integration_general_FFPE_hd}
    \begin{aligned}
        (2 \pi)^{d} p(\boldsymbol{x}, t)
        &= \int_{\mathbb{R}^{d}} \hat{p}(\boldsymbol{k}, t) \exp \left ( i \langle \boldsymbol{k}, \boldsymbol{x} \rangle \right ) \, d k\\
        &= \int_{\mathbb{R}^{d}} \exp \left ( i \langle \boldsymbol{k}, \boldsymbol{x} - \boldsymbol{x}_{0} - \boldsymbol{\drift} t \rangle \right ) \exp \left ( - \left ( \diffusionOrdinary |\boldsymbol{k}|^{2} + \diffusionFractional |\boldsymbol{k}|^{2 \alpha} \right ) t \right ) \, d \boldsymbol{k}\\
        &= \int_{\mathbb{R}^{d}} \cos \left ( \langle \boldsymbol{k}, \boldsymbol{x} - \boldsymbol{x}_{0} - \boldsymbol{\drift} t \rangle \right ) \exp \left ( - \left ( \diffusionOrdinary |\boldsymbol{k}|^{2} + \diffusionFractional |\boldsymbol{k}|^{2 \alpha} \right ) t \right ) \, d \boldsymbol{k}.
    \end{aligned}
\end{equation}

Notice that for any vector $\boldsymbol{x} - \boldsymbol{x}_{0} - \boldsymbol{\drift} t$, there exists a rotation matrix $R \in \mathbb{R}^{d \times d}$ that transforms $\boldsymbol{x} - \boldsymbol{x}_{0} - \boldsymbol{\drift} t$ into $y \boldsymbol{e}_{1}$, where $y=|\boldsymbol{x} - \boldsymbol{x}_{0} - \boldsymbol{\drift} t|$ and $\boldsymbol{e}_{1}$ denotes the first standard unit vector, defined by having its first component equal to one and all other components equal to zero.
As a result, \Cref{eq:Dirac_explicit_integration_general_FFPE_hd} can be reformulated as follows:
\begin{equation}\label{eq:Dirac_explicit_integration_general_FFPE_hd_with_rotation}
    \begin{aligned}
        (2 \pi)^{d} p(\boldsymbol{x}, t)
        &= \int_{\mathbb{R}^{d}} \cos \left ( \langle \boldsymbol{k}, \boldsymbol{x} - \boldsymbol{x}_{0} - \boldsymbol{\drift} t \rangle \right ) \exp \left ( - \left ( \diffusionOrdinary |\boldsymbol{k}|^{2} + \diffusionFractional |\boldsymbol{k}|^{2 \alpha} \right ) t \right ) \, d \boldsymbol{k}\\
        &= \int_{\mathbb{R}^{d}} \cos \left ( \langle R \boldsymbol{k}, R ( \boldsymbol{x} - \boldsymbol{x}_{0} - \boldsymbol{\drift} t ) \rangle \right ) \exp \left ( - \left ( \diffusionOrdinary |\boldsymbol{k}|^{2} + \diffusionFractional |\boldsymbol{k}|^{2 \alpha} \right ) t \right ) \, d \boldsymbol{k}\\
        &= \int_{\mathbb{R}^{d}} \cos \left ( y \langle R \boldsymbol{k}, \boldsymbol{e}_{1} \rangle \right ) \exp \left ( - \left ( \diffusionOrdinary |R \boldsymbol{k}|^{2} + \diffusionFractional |R \boldsymbol{k}|^{2 \alpha} \right ) t \right ) \, d (R \boldsymbol{k})\\
        &= \int_{\mathbb{R}^{d}} \cos \left ( y \langle \boldsymbol{k}, \boldsymbol{e}_{1} \rangle \right ) \exp \left ( - \left ( \diffusionOrdinary |\boldsymbol{k}|^{2} + \diffusionFractional |\boldsymbol{k}|^{2 \alpha} \right ) t \right ) \, d \boldsymbol{k}.
    \end{aligned}
\end{equation}
From the above equation, for fixed diffusion coefficients $\diffusionOrdinary$ and $\diffusionFractional $, $p( \boldsymbol{x}, t)$ only depends on the displacement $y$ and $t$.
Therefore we denote $\tilde{p}(y, t) := p( \boldsymbol{x}, t)$.

Additionally, for $d \geq 2$, we employ the $d$-dimensional polar coordinate system, which consist of a radial coordinate $r$ and $d - 1$ angular coordinates $\theta_{1}, \ldots, \theta_{d - 1}$, where the angles $\theta_{1}, \ldots, \theta_{d - 2}$ range over $[0, \pi)$ radians and $\theta_{d - 1}$ ranges over $[0, 2 \pi)$ radians. It follows that
\begin{equation*}
    \boldsymbol{k}
    = \left (
    \begin{aligned}
        &r \cos(\theta_{1})\\
        &r \sin(\theta_{1}) \cos(\theta_{2})\\
        &r \sin(\theta_{1}) \sin(\theta_{2}) \cos(\theta_{3})\\
        &\quad \vdots\\
        &r \sin(\theta_{1}) \cdots \sin(\theta_{d - 2}) \cos(\theta_{d - 1})\\
        &r \sin(\theta_{1}) \cdots \sin(\theta_{d - 2}) \sin(\theta_{d - 1})
    \end{aligned}
    \right )
\end{equation*}
and
\begin{equation*}
    \begin{aligned}
        d \boldsymbol{k}
        = r^{d - 1} \sin^{d - 2}(\theta_{1}) \sin^{d - 3}(\theta_{2}) \cdots \sin(\theta_{d - 2}) \, d r \, d \theta_{1} \, d \theta_{2} \cdots \, d \theta_{d - 1}.
    \end{aligned}
\end{equation*}
Upon converting into $d$-dimensional polar coordinates, \Cref{eq:Dirac_explicit_integration_general_FFPE_hd_with_rotation} with $d\geq 2$ is transformed as follows:
\begin{equation}\label{eq:Dirac_explicit_integration_general_FFPE_hd_polar}
    \begin{aligned}
        &\quad \, (2 \pi)^{d} p(\boldsymbol{x}, t)
        = \int_{\mathbb{R}^{d}} \cos \left ( y \langle \boldsymbol{k}, \boldsymbol{e}_{1} \rangle \right ) \exp \left ( - \left ( \diffusionOrdinary |\boldsymbol{k}|^{2} + \diffusionFractional |\boldsymbol{k}|^{2 \alpha} \right ) t \right ) \, d \boldsymbol{k}\\
        &= S_{d - 2} \int_{0}^{\pi} \sin^{d - 2}(\theta) \int_{0}^{\infty} r^{d - 1} \cos \left ( \cos(\theta) y r \right ) \exp \left ( - \left ( \diffusionOrdinary r^{2} + \diffusionFractional r^{2 \alpha} \right ) t \right ) \, d r \, d \theta,
    \end{aligned}
\end{equation}
where $S_{d - 2} = 2 \pi^{\frac{d - 1}{2}} / \Gamma \left ( \frac{d - 1}{2} \right )$ represents the surface area of the unit $(d - 2)$-sphere embedded in $(d - 1)$-dimensional Euclidean space and $\Gamma$ denotes the gamma function
\begin{equation*}
    \Gamma(z)
    = \int_{0}^{\infty} t^{z - 1} \exp(- t) \, d t.
\end{equation*}

Using Poisson's representation for Bessel function of the first kind for $\nu > - 1 / 2$,
\begin{equation*}
     J_{\nu}(z)
     = \frac{(z / 2)^{\nu}}{\pi^{1 / 2} \Gamma(\nu + 1 / 2)} \int_{0}^{\pi} \sin^{2 \nu}(\theta) \cos(\cos(\theta) z) \, d \theta,
\end{equation*}
one can reformulate \Cref{eq:Dirac_explicit_integration_general_FFPE_hd_polar} as
\begin{equation}\label{eq:Dirac_general_FFPE_Bessel}
    \begin{aligned}
        &\quad \, \tilde{p}(y, t)
        = p(\boldsymbol{x}, t)\\
        &= \frac{1}{y^{(d - 2) / 2}} \int_{0}^{\infty} \left ( \frac{r}{2 \pi} \right )^{d / 2} J_{(d - 2) / 2}(y r) \exp \left ( - \left ( \diffusionOrdinary r^{2} + \diffusionFractional r^{2 \alpha} \right ) t \right ) \, d r
    \end{aligned}
\end{equation}
for any $d \geq 2$ and $y \neq 0$.
For $d=1$, we simply use \Cref{eq:Dirac_explicit_integration_general_FFPE_hd_with_rotation} and write
\begin{equation}\label{eq:Dirac_general_FFPE_1d}
    \tilde{p}(y, t)
    = p(x, t)
    = \frac{1}{\pi} \int_{0}^{\infty} \cos(y r) \exp \left ( - \left ( \diffusionOrdinary r^{2} + \diffusionFractional r^{2 \alpha} \right ) t \right ) \, d r.
\end{equation}
Notice that when $y = 0$, one can use \Cref{eq:Dirac_explicit_integration_general_FFPE_hd_with_rotation} to get
\begin{equation}\label{eq:Dirac_general_FFPE_zero_displacement}
    \begin{aligned}
        \tilde{p}(0, t) 
        = \frac{S_{d - 1}}{2 \pi} \int_{0}^{\infty} \left ( \frac{r}{2 \pi} \right )^{d - 1} \exp \left ( - \left ( \diffusionOrdinary r^{2} + \diffusionFractional r^{2 \alpha} \right ) t \right ) \, d r
    \end{aligned}
\end{equation}
for any $d \geq 1$.
Actually, \Cref{eq:Dirac_general_FFPE_Bessel} serves as the ultimate solution to the FFPE as it also holds for $d = 1$, though it is less computationally efficient compared to \Cref{eq:Dirac_general_FFPE_1d}.
Additionally, by taking the limit $y \to 0$, \Cref{eq:Dirac_general_FFPE_zero_displacement} can be obtained.

In \Cref{sec:numerical_approach}, we will employ representations provided in \Cref{eq:Dirac_general_FFPE_Bessel,eq:Dirac_general_FFPE_1d,eq:Dirac_general_FFPE_zero_displacement} to calculate the solution to \Cref{eq:Dirac_general_FFPE_hd}.
During this process, we will confront two primary challenges. 
The first concerns the behavior near the origin, where the integrand becomes singular due to the term $|r|^{2 \alpha}$ when $\alpha$ lies in $(0, 1)$.
In particular, the integrand $\exp( - \diffusionFractional |r|^{2 \alpha} t )$ is not differentiable at $r = 0$ when $\alpha \in (0, 1 / 2)$, and fails to be twice differentiable when $\alpha \in (1 / 2, 1)$. 
This lack of smoothness presents challenges for accurate numerical integration.
A graphical illustration of this near-origin behavior is provided in \Cref{fig:Dirac_singularity_in_integral}.
The second challenge lies in the far-field region, where the integrand may exhibit a slow rate of decay, especially when the ordinary diffusion vanishes.
Note that the formulas also indicate that solutions in higher dimensions are harder to approximate, as the integrand decays more slowly with increasing $d$ due to the presence of the term $r^{d / 2}$ in \Cref{eq:Dirac_general_FFPE_Bessel} and $r^{d - 1}$ in \Cref{eq:Dirac_general_FFPE_zero_displacement}.

\begin{remark}
    Observe that in \Cref{eq:Dirac_general_FFPE_Bessel,eq:Dirac_general_FFPE_zero_displacement}, $2 \pi$ is incorporated into the integrand to avoid scenarios where a large quantity is divided by another large quantity when $d$ is large.
    Additionally, for large $d$, the decay rate in the far-field integration is decelerated by the positive power of $r$ in the integrand, dividing by $2 \pi$ helps to slightly alleviate this effect.
\end{remark}

\subsection{The scaling law}\label{sec:scaling_law}
In this subsection, we present scaling laws essential for the numerical method to be discussed in \Cref{sec:numerical_approach}.
To establish the scaling law, we first express the dependence on $\diffusionOrdinary$ in the solution, denoted as $\tilde{p} (y, t; \diffusionOrdinary) = p(\boldsymbol{x}, t; \diffusionOrdinary)$.
Starting from \Cref{eq:Dirac_explicit_integration_general_FFPE_hd_with_rotation}, we derive the scaling law as follows:
\begin{equation}\label{eq:Dirac_general_FFPE_hd_scaling}
    \begin{aligned}
        &\quad \, \tilde{p}(y, t; \diffusionOrdinary)
        = p(\boldsymbol{x}, t; \diffusionOrdinary)\\
        &= \frac{1}{(2 \pi)^{d}} \int_{\mathbb{R}^{d}} \cos \left ( y \xi_{1} \right ) \exp \left ( - \left ( \diffusionOrdinary |\boldsymbol{\xi}|^{2} + \diffusionFractional |\boldsymbol{\xi}|^{2 \alpha} \right ) t \right ) \, d \boldsymbol{\xi}\\
        &= \resizebox{0.85\linewidth}{!}{$\displaystyle
        \frac{1}{(2 \pi)^{d}} \int_{\mathbb{R}^{d}} \cos \left ( y \xi_{1} \right ) \exp \left ( - \left ( \left ( \frac{t}{T} \right )^{1 - \frac{1}{\alpha}} \diffusionOrdinary \left ( \frac{t}{T} \right )^{\frac{2}{2 \alpha}} |\boldsymbol{\xi}|^{2} + \diffusionFractional \left ( \frac{t}{T} \right ) |\boldsymbol{\xi}|^{2 \alpha} \right ) T \right ) \, d \boldsymbol{\xi}
        $}\\
        &= \resizebox{0.85\linewidth}{!}{$\displaystyle
        \frac{(t/T)^{- \frac{d}{2 \alpha}}}{(2 \pi)^{d}} \int_{\mathbb{R}^{d}} \cos \left ( \left ( \frac{t}{T} \right )^{- \frac{1}{2 \alpha}} y k_{1} \right ) \exp \left ( - \left ( \left ( \frac{t}{T} \right )^{1 - \frac{1}{\alpha}} \diffusionOrdinary |\boldsymbol{k}|^{2} + \diffusionFractional |\boldsymbol{k}|^{2 \alpha} \right ) T \right ) \, d \boldsymbol{k}
        $}\\
        &= \left ( \frac{t}{T} \right )^{- \frac{d}{2 \alpha}} \tilde{p} \left ( \left ( \frac{t}{T} \right )^{- \frac{1}{2 \alpha}} y, T; \left ( \frac{t}{T} \right )^{1 - \frac{1}{\alpha}} \diffusionOrdinary \right ),
    \end{aligned}
\end{equation}
for any $T > 0$.

\begin{figure}[htp]
    \centering
    \includegraphics[width=\textwidth]{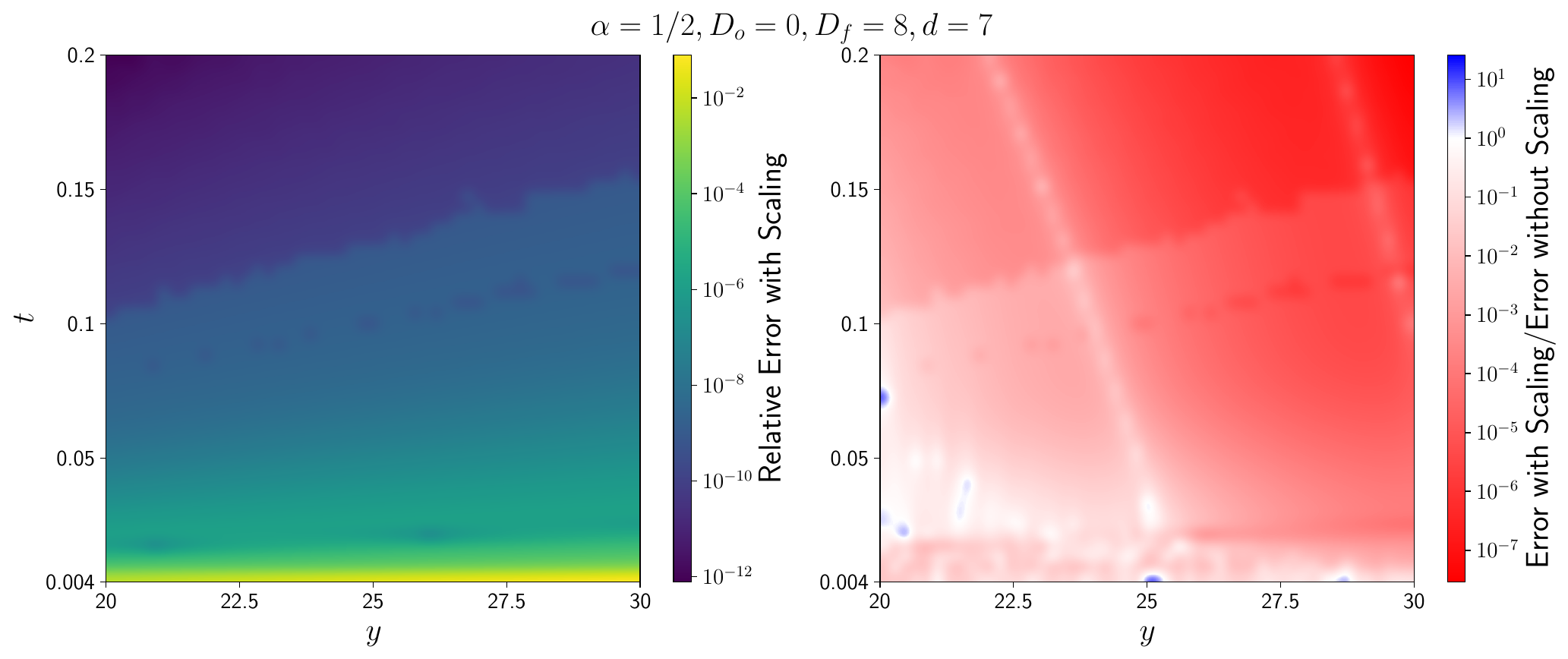}
    \caption{
        Relative error comparison with and without scaling for $\alpha = 1 / 2$, $\diffusionOrdinary = 0$, $\diffusionFractional = 8$, and $d = 7$.
        The left plot shows the relative error with scaling, while the right plot shows the ratio of the error with scaling to the error without scaling.
        In the right plot, areas where the ratio is smaller than 1 (indicating that scaling performs better) are shown in red, and areas where the ratio is larger than 1 (indicating that scaling performs worse) are shown in blue.
        The color bars indicate the magnitude of the errors and their ratios on a logarithmic scale.
        These plots illustrate the effectiveness of scaling in reducing the error across different values of $y$ and $t$ for large $y$.
    }\label{fig:Dirac_force_scaling_comparison}
\end{figure}

The scaling mentioned above allows us to adjust the coefficient $\diffusionOrdinary$ or the displacement $y$ to different scales to aid the numerical evaluation of integrals.
This is a crucial aspect of our approach.
In \Cref{sec:numerical_approach}, we will introduce a technique for treating far-field integration, which is particularly effective when $y$ is on the order of $O(1)$.
To demonstrate the importance of the scaling technique, we present the error comparisons of the algorithm with and without the scaling technique in \Cref{fig:Dirac_force_scaling_comparison}.
The results clearly show that for large values of $y$, the scaling technique offers significant advantages.
Details of the numerical algorithm are provided in \Cref{sec:numerical_approach}.

\subsection{Zero displacement}\label{sec:zero_displacement}
In the case of zero displacement, i.e., $y = 0$, the scaling identity \eqref{eq:Dirac_general_FFPE_hd_scaling} reduces to
\begin{equation*}
    \tilde{p} (0,t; \diffusionOrdinary)
    = \left ( \frac{t}{T} \right )^{- \frac{d}{2 \alpha}} \tilde{p} \left ( 0, T; \left ( \frac{t}{T} \right )^{1 - \frac{1}{\alpha}} \diffusionOrdinary \right ). 
\end{equation*}
Using \Cref{eq:Dirac_general_FFPE_zero_displacement}, we further get
\begin{equation*}
    \tilde{p}(0, t; \diffusionOrdinary)
    = \frac{S_{d - 1}}{(2 \pi)^{d}}
    \left ( \frac{t}{T} \right )^{- \frac{d}{2 \alpha}} \int_{0}^{\infty} r^{d - 1} \exp \left ( - \left ( \left ( \frac{t}{T} \right )^{1 - \frac{1}{\alpha}} \diffusionOrdinary r^{2} + \diffusionFractional r^{2 \alpha} \right ) T \right ) \, d r,
\end{equation*} 
for any $T > 0$.

If $\diffusionOrdinary > 0$, one can select appropriate $T > 0$ such that $\left ( \frac{t}{T} \right )^{1 - \frac{1}{\alpha}} \diffusionOrdinary T$ is sufficiently large to allow for the straightforward computation of the above integral using truncation.
To be more precise, one can estimate the error caused by truncation as follows:
\begin{equation*}
    \begin{aligned}
        &\quad \, \int_{R}^{\infty} r^{d - 1} \exp \left ( - \left ( \left ( \frac{t}{T} \right )^{1 - \frac{1}{\alpha}} \diffusionOrdinary r^{2} + \diffusionFractional r^{2 \alpha} \right ) T \right ) \, d r\\
        &\leq \int_{R}^{\infty} r^{d - 1} \min \left\{ \exp \left ( - \xi_{1} r^{2} \right ) , \exp \left ( - \xi_{2} r^{2 \alpha} \right )\right\} \, d r,
    \end{aligned}
\end{equation*}
where $R > 0$, $\xi_{1} = \left ( \frac{t}{T} \right )^{1 - \frac{1}{\alpha}} \diffusionOrdinary T$ and $\xi_{2} = \diffusionFractional T$.
Utilizing the upper incomplete Gamma function, we have
\begin{equation*}
    \int_{R}^{\infty} r^{q} \exp \left ( - \xi r^{p} \right ) \, d r
    = \frac{\gamma_{\textnormal{upper}} \left ( \frac{q + 1}{p}, \xi R^{p} \right )}{\xi^{\frac{q + 1}{p}} p},
\end{equation*}
for any $p$ and $q$.
For a reasonably large $R > 0$, one can then compute the minimum values of $\xi_{1}$ and $\xi_{2}$ such that the above quantity drops below a certain specified accuracy threshold.
In this way, one can determine the appropriate $T$ to use in practice.

In the case $\diffusionOrdinary = 0$, one can find the explicit formula of the solution using special functions. In particular,
\begin{equation*}
    \begin{aligned}
        \int_{0}^{\infty} r^{d - 1} \exp \left ( - r^{2 \alpha} \tau \right ) \, d r
        &= \frac{1}{2 \alpha} \tau^{- \frac{d}{2 \alpha}} \Gamma \left ( \frac{d}{2 \alpha} \right )
        = \frac{1}{d} \tau^{- \frac{d}{2 \alpha}} \Gamma \left ( \frac{d}{2 \alpha} + 1 \right ).
    \end{aligned}
\end{equation*}
Therefore,
\begin{equation}\label{eq:Dirac_FFPE_general_hd_special_case_no_displacement_no_ordinary_diffusion}
    \begin{aligned}
        \tilde{p}(0, t; 0)
        &= \frac{S_{d - 1}}{(2 \pi)^{d} d} \left ( \diffusionFractional t \right )^{- \frac{d}{2 \alpha}} \Gamma \left ( \frac{d}{2 \alpha} + 1 \right ).
    \end{aligned}
\end{equation}
We note that the special functions $\Gamma(\cdot)$ and $\gamma_{\textnormal{upper}}(\cdot, \cdot)$ used in the above expressions are standard and can be evaluated efficiently using rapidly converging series or continued fractions.
These functions are readily available in most numerical computing environments, allowing for direct use without the need for custom implementation.

%%%%%%%%%%%%%%%%%%%%%%%%%%%Section 3 %%%%%%%%%%%%%%%%%%%%%%%%%%

\section{Numerical approaches}\label{sec:numerical_approach}
In this section, we present the solvers for the FFPE using the integral representations.
For numerical computations, we use \Cref{eq:Dirac_general_FFPE_1d} for one-dimensional cases and \Cref{eq:Dirac_general_FFPE_Bessel} for dimensions higher than one.
The treatment for the special case of zero displacement follows the discussion in \Cref{sec:zero_displacement}.
When the scaling law is needed, we apply \Cref{eq:Dirac_general_FFPE_hd_scaling}.

The general solver is presented in \Cref{sec:FFPE_solver}, while necessary integration techniques for near-origin and far-field integrals are presented in \Cref{sec:singular_integration} and \Cref{sec:far_field_integration}), respectively.
For a succinct summary of the methodology, one is advised to consult \Cref{rmk:Dirac_concise_explanation}.

\begin{remark}\label{rmk:Dirac_concise_explanation}
    Our methodology is summarized as follows:
    \begin{equation}
        \tilde{p}(y, t)
        = \underbrace
        {
            \underbrace
            {
                \underbrace
                {
                    \tilde{g}_{\text{local1}}(y, t)
                }
                _{
                    \text{\Cref{alg:Dirac_integration_with_singularity}}
                }
                +
                \underbrace
                {
                    \tilde{g}_{\text{local2}}(y, t)
                }
                _{
                    \text{\Cref{alg:Dirac_integration_for_slow_decay}}
                }
            }
            _{
                \text{\Cref{alg:Dirac_basic_solver_for_general_FFPE_with_zero_displacement,alg:Dirac_basic_solver_for_general_FFPE}}
            }
            +
            \underbrace
            {
                \tilde{g}_{\text{far}}(y, t)
            }
            _{
                \text{error}
            }
        }
        _{
            \text{\Cref{alg:Dirac_solver_for_general_FFPE_in_hd}}
        }.
    \end{equation}
    It is important to note that \Cref{alg:Dirac_solver_for_general_FFPE_in_hd} is specifically designed to handle different cases using \Cref{alg:Dirac_basic_solver_for_general_FFPE_with_zero_displacement,alg:Dirac_basic_solver_for_general_FFPE}, potentially incorporating a scaling approach.
\end{remark}

\subsection{Techniques for the singular integration}\label{sec:singular_integration}
\begin{figure}[htp]
    \centering
    \includegraphics[scale=0.6]{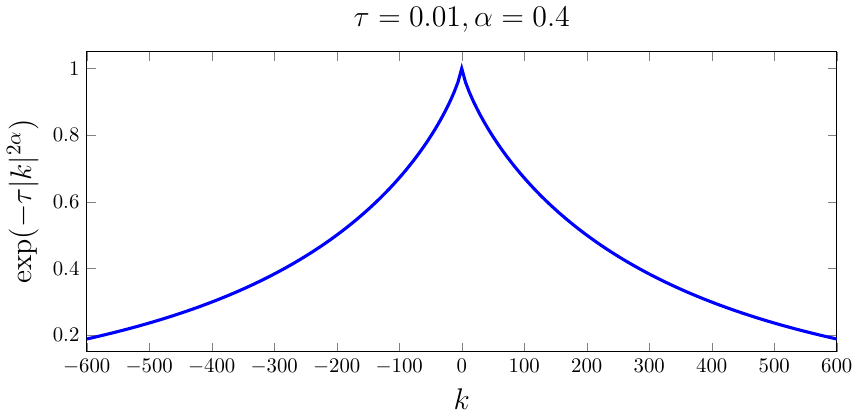}
    \caption{
        A graphical representation of the function $\exp(- \tau |k|^{2\alpha})$, where $\tau = 0.01$ and $\alpha = 0.4$.
        The plot highlights the near-origin singularity and the far-field slow decay.
    }\label{fig:Dirac_singularity_in_integral}
\end{figure}

It is important to observe that for $\alpha$ within the interval $(0, 1 / 2)$, the function $\exp \left( - |k|^{2 \alpha} \tau \right)$ lacks differentiability at $k = 0$.
Moreover, when $\alpha$ falls within the range $(1 / 2, 1)$, the function $\exp \left( - |k|^{2 \alpha} \tau \right)$ is not second-order differentiable at $k = 0$.
For a graphical representation of the singularity near the origin, refer to \Cref{fig:Dirac_singularity_in_integral}.

To mitigate the challenges posed by singularity of integrands near origin, we introduce two techniques for integrals of the form
\begin{equation}\label{eq:singular_integral}
    \int_{0}^{L} f(z) \exp \left ( - z^{2 \alpha} \tau \right ) \, d z,
\end{equation}
where $f$ denotes a smooth function and $L\in (0, 1]$.
For solving the FFPE, the function $f$ in \Cref{eq:singular_integral} is given explicitly in \Cref{alg:Dirac_basic_solver_for_general_FFPE_with_zero_displacement,alg:Dirac_basic_solver_for_general_FFPE}, and those expressions arise from the solution representations given by \Cref{eq:Dirac_general_FFPE_Bessel,eq:Dirac_general_FFPE_1d,eq:Dirac_general_FFPE_zero_displacement}.
Note that for solving the FFPE, we fix $L = 1$ in our algorithm later.
The first method involves approximating $\exp \left ( - z^{2 \alpha} \tau \right )$ by its Taylor polynomial and then applying the Gauss-Jacobi quadrature for each integral.
This method is effective for smaller values of $\tau$ when Taylor polynomials converge fast.
This method is named the expansion technique with further details in \Cref{sec:expansion_technique}.
The second method, suitable for cases when $\tau$ has larger values, relies on calculating the weights $\{ w_{j} \}$ such that the approximation $\int_{0}^{L} f(z) \exp \left ( - z^{2 \alpha} \tau \right ) \, d z \approx \sum_{j} w_{j} f(s_{j})$ is precise up to certain order of polynomials.
This method is referred to as the re-weighting technique and is detailed in \Cref{sec:reweighting_technique}.
Combining the two techniques, the algorithm for approximating integrals of the form \eqref{eq:singular_integral} is summarized in \Cref{alg:Dirac_integration_with_singularity}.
\begin{algorithm}[htbp]
    \caption{Integration with Singularity}\label{alg:Dirac_integration_with_singularity}
    \begin{algorithmic}[1]
        \Require A smooth function denoted by $f$, an upper limit of integration labeled as $L \in (0, 1]$, a fractional exponent signified by $\alpha \in (0, 1)$, and a pseudo-temporal coefficient represented by $\tau > 0$
        \Ensure Computation of the integral $\displaystyle \int_{0}^{L} f(z) \exp \left ( - z^{2 \alpha} \tau \right ) \, d z$

        \State $N, K_{1}, K_{2}, K_{3}, K_{4} \gets 16, 17, 9, 6, 4$
            \Comment{Adjustable}

        \If{$\tau > \textit{1}$ \textbf{or} precomputed}
            \State Apply the \Call{ReweightingTechnique}{} to ascertain the quadrature points $\{ s_{j} \}_{j = 1}^{N}$ and the weights $\{ w_{j} \}_{j = 1}^{N}$
            \Statex \Comment{Can be precomputed for repetitive evaluations under identical $L$, $\alpha$, and $\tau$}
            \State \Return $\sum_{j = 1}^{N} w_{j} f(s_{j})$
        \Else
            \State \Return Outcome via the \Call{ExpansionTechnique}{} with $K$ selected as follows
            \Statex
            \begin{center}
                \begin{tikzpicture}[scale=1]
                    \draw[->] (0, 0) -- (10, 0) node[right] {$\tau$};
\coordinate (tick1) at (3, 0);
\coordinate (tick2) at (5, 0);
\coordinate (tick3) at (7, 0);
\coordinate (tick4) at (9, 0);

\foreach \tick in {tick1, tick2, tick3, tick4}
{
    \draw (\tick) -- ++(0, +0.2);
}

\node[below] at (tick1) {\textit{1e-3}};
\node[below] at (tick2) {\textit{1e-2}};
\node[below] at (tick3) {\textit{1e-1}};
\node[below] at (tick4) {\textit{1}};

\node[above] at (2, 0) {$K_{4}$};
\node[above] at (4, 0) {$K_{3}$};
\node[above] at (6, 0) {$K_{2}$};
\node[above] at (8, 0) {$K_{1}$};

                \end{tikzpicture}
            \end{center}
        \EndIf
    \end{algorithmic}
\end{algorithm}

\subsubsection{The expansion technique}\label{sec:expansion_technique}
The Gauss-Jacobi quadrature rule can be directly applied to compute integrals of the form
\begin{equation*}
    \int_{0}^{L} f(z) z^{\beta} \, d z.
\end{equation*}
By expanding the exponential function, we derive
\begin{equation*}
    \begin{aligned}
        \int_{0}^{L} f(z) \exp \left ( - z^{2 \alpha} \tau \right ) \, d z
        \approx \sum_{\kappa = 0}^{K} \frac{(- \tau)^{\kappa}}{\kappa!} \int_{0}^{L} f(z) z^{2 \alpha \kappa} \, d z.
    \end{aligned}
\end{equation*}
Notice that the accuracy of the approximation increases with larger values of $K$.
For a given $\tau > 0$, the recommended value of $K$ can be found in \Cref{tab:Dirac_expansion_relation}, and is chosen as the smallest $K$ such that $\tau^{K + 1} / (K + 1)!$ falls below machine precision.
This is effective particularly when the integral $\int_{0}^{L} |f(z)| \, d z$ remains relatively small.
We suggest employing the expansion technique when $\tau \leq 1$.
\begin{table}[htbp]
    \centering
    \resizebox{0.9\linewidth}{!}{
    \begin{tabular}{|c|c|c|c||c|c|c|c|}
        \hline
        &&&&&&&\\[-10pt]
        $\tau$ & $K$ & $\tau^{K} / K!$ & $\tau^{K + 1} / (K + 1)!$ & $\tau$ & $K$ & $\tau^{K} / K!$ & $\tau^{K + 1} / (K + 1)!$\\
        \hline
        \textit{1e-3} & $4$ & \textit{4.1667e-14} & \textit{8.3333e-18}
        & \textit{1} & $17$ & \textit{2.8115e-15} & \textit{1.5619e-16}\\
        \textit{1e-2} & $6$ & \textit{1.3889e-15} & \textit{1.9841e-18}
        & \textit{3} & $27$ & \textit{7.0031e-16} & \textit{7.5033e-17}\\
        \textit{1e-1} & $9$ & \textit{2.7557e-15} & \textit{2.7557e-17}
        & \textit{10} & $51$ & \textit{6.4470e-16} & \textit{1.2398e-16}\\
        \hline
    \end{tabular}
    }
    \caption{For a given $\tau > 0$, the recommended value of $K$ can be found in this table. The value of $K$ is chosen as the smallest $K$ such that $\tau^{K + 1} / (K + 1)!$ falls below machine precision. The table lists values of $\tau$ and their corresponding $K$ values, along with the computed values of $\tau^{K} / K!$ and $\tau^{K + 1} / (K + 1)!$.}\label{tab:Dirac_expansion_relation}
\end{table}
For values of $\tau$ greater than 1, we suggest adopting the alternative method described in the following subsection.

\subsubsection{The re-weighting technique}\label{sec:reweighting_technique}
We assume the following quadrature approximation:
\begin{equation}\label{eq:reweighting_technique}
    \int_{0}^{L} f(z) \exp \left ( - z^{2 \alpha} \tau \right ) \, d z
    \approx \sum_{j = 1}^{N} w_{j} f(s_{j}).
\end{equation}
The ideal Gaussian quadrature rule seeks $\{ s_{j} \}$ and $\{ w_{j} \}$ such that \Cref{eq:reweighting_technique} yields an exact result for $f$ being a polynomial of degree $2 N - 1$ or less.
However, orthogonal polynomials with respect to the weight function $\exp \left ( - z^{2 \alpha} \tau \right )$ are hard to obtain.
We therefore choose to fix the quadrature points $\{ s_{j} \}$ and solve a linear system for the weights $\{ w_{j} \}$. This approach requires that equation \Cref{eq:reweighting_technique} be exact for polynomials up to degree $N - 1$.
For example, one may fix $\{ s_{j} \} \subset [0, L]$ as the roots of the Legendre polynomials and find $\{ w_{j} \}$ such that
\begin{equation}\label{eq:polynomialreproducation}
    \int_{0}^{L} f(z) \exp \left ( - z^{2 \alpha} \tau \right ) \, d z
    = \sum_{j = 1}^{N} w_{j} f(s_{j})
    \qquad \forall f \in \mathcal{P}_{N - 1}(\mathbb{R}),
\end{equation}
where $\mathcal{P}_{N - 1}(\mathbb{R})$ denotes the polynomial space of order $N - 1$.
In practice, we recommend using Legendre polynomials to form a basis for $\mathcal{P}_{N - 1}(\mathbb{R})$ for better conditioning.
For a given polynomial $f$, the left-hand side of \Cref{eq:polynomialreproducation} can be approximated accurately by high-order quadrature and then $\{ w_{j} \}$ can be solved from a linear system.
We emphasize that the computation of weights $\{ w_{j} \}$ is performed only once for fixed parameters $L$, $\alpha$, and $\tau$, and the resulting values can be reused for all subsequent evaluations under the same settings.
Since this is the typical case in most applications, the cost of computing these weights is negligible in practice.
This reuse strategy is also reflected in \Cref{alg:Dirac_integration_with_singularity}.

\subsection{Techniques for the far-field integration}\label{sec:far_field_integration}
\begin{figure}[htp]
    \centering
    \includegraphics[scale=0.6]{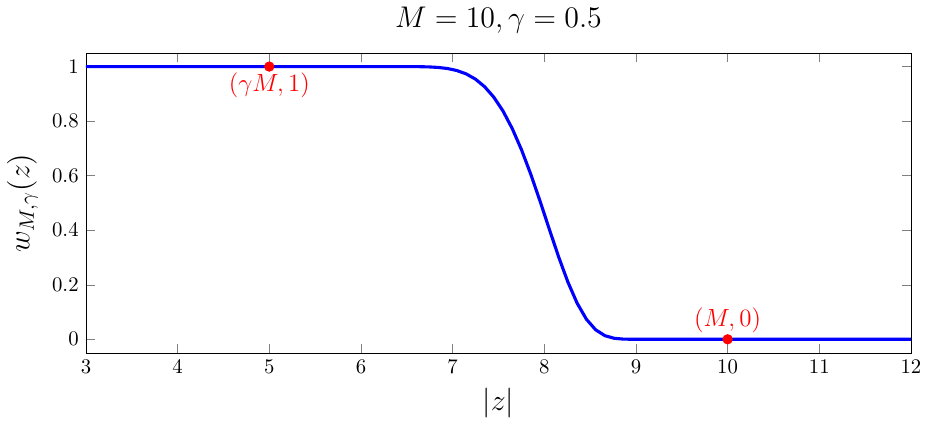}
    \caption{An example of the windowing function $w_{M,\gamma}(z)$ with parameters $M = 10$ and $\gamma = 0.5$. The key points $(\gamma M, 1)$ and $(M, 0)$ are marked in red. This illustrates how the function smoothly transitions from $1$ to $0$ around the interval $[\gamma M, M]$, demonstrating the behavior of the windowing function in truncating values outside the specified range.}\label{fig:Dirac_windowing_function_example}
\end{figure}

The function $\exp \left( - |k|^{2 \alpha} \tau \right )$ has a slow decay as $k\to \infty$ (see \Cref{fig:Dirac_singularity_in_integral}), especially when $\alpha$ is small with the range $(0, 1)$.
This slow decay presents difficulties for numerical integration in the far field, as direct numerical integration would require a significantly large truncation domain.
To address this challenge, we partition the integral into local and far-field parts through the application of a windowing function \cite{bruno2017windowed}, which serves as a partition of unity:
\begin{equation}\label{eq:partition_of_unity}
    \begin{aligned}
        &\quad \, \tilde{g}(y, t)
        := \int_{0}^{\infty} f(z; y) \exp \left ( - z^{2 \alpha} \tau \right ) \, d z\\
        &= \int_{0}^{\infty} f(z; y) \exp \left ( - z^{2 \alpha} \tau \right ) w_{M, \gamma}(z) \, d z\\
        &+ \int_{0}^{\infty} f(z; y) \exp \left ( - z^{2 \alpha} \tau \right ) \left ( 1 - w_{M, \gamma}(z) \right ) \, d z\\
        &=: \tilde{g}_{\text{local}}(y, t) + \tilde{g}_{\text{far}}(y, t),
    \end{aligned}
\end{equation}
where $f$ denotes a smooth function, the oscillation frequency of which is dictated by $y$.
The windowing function $w_{M, \gamma}$ is defined as
\begin{equation*}
    w_{M, \gamma}(z)
    = \left \{
    \begin{aligned}
        &1, && s \leq 0\\
        &\exp \left ( - 2 \frac{\exp(- 1 / s^{2})}{(1 - s)^{2}} \right ), && 0 < s < 1\\
        &0, && s \geq 1
    \end{aligned}
    \right .,
\end{equation*}
where the parameter $s$ is given by
\begin{equation*}
    s(z)
    = \frac{|z| - \gamma M}{M - \gamma M},
\end{equation*}
with $M$ and $\gamma$ being two constants that satisfy $M > 0$ and $0 < \gamma < 1$.
The windowing function $w_{M, \gamma}(z)$ is designed to be $1$ for $|z| \leq \gamma M$ and $0$ for $|z| \geq M$, with a smooth transition between these values.
Additionally, all derivatives gradually approach zero at the boundary points $|z| = M$ and $|z| = \gamma M$.
An illustrative example of the windowing function is provided in \Cref{fig:Dirac_windowing_function_example}.
The windowing function is designed such that $\tilde{g}_{\text{far}}$ is small and can be neglected in practice.
In our experiments, setting $\gamma = 0.5$ is proved effective, ensuring that the transition of the windowing function is neither too steep nor too gradual. Consequently, we fix $\gamma = 0.5$ in the subsequent work and denote $w_M := w_{M, 0.5}$ for simplicity.

To address the singularity at the origin, the local part in \Cref{eq:partition_of_unity} is further divided into two parts:
\begin{equation}
    \begin{aligned}
        &\quad \ \tilde{g}_{\text{local}}(y, t)\\
        &= \int_{0}^{\infty} f(z; y) \exp \left ( - z^{2 \alpha} \tau \right ) w_{M, \gamma}(z) \, d z
        = \int_{0}^{M} f(z; y) \exp \left ( - z^{2 \alpha} \tau \right ) w_{M, \gamma}(z) \, d z\\
        &=\int_{0}^{L} f(z; y) \exp \left ( - z^{2 \alpha} \tau \right ) \, d z
        + \int_{L}^{M} f(z; y) \exp \left ( - z^{2 \alpha} \tau \right ) w_{M, \gamma}(z) \, d z\\
        &=: \tilde{g}_{\text{local1}}(y, t) + \tilde{g}_{\text{local2}}(y, t),
    \end{aligned}
\end{equation}
where we implicitly used $L \leq \gamma M$ in the above derivation, a condition that holds in all our experiments.
For the evaluation of $\tilde{g}_{\text{local1}}$, we employ the techniques described in \Cref{sec:singular_integration}, and for evaluating $\tilde{g}_{\text{local2}}$, Gauss–Legendre quadrature is used, where the number of quadrature points is linearly proportional to $M$. 

The subsequent simplified example illustrates the effectiveness of using the windowing function by contrasting two approaches: one involving direct truncation of the integral, and the other incorporating truncation after the application of the windowing function. Consider the following integrals by choosing $L = 1$, $f(z; y) = \cos(z)$, $\tau = 0.01$, and $\alpha = 0.4$: 
\begin{equation*}
    \begin{aligned}
        I
        &= \int_{1}^{\infty} \cos( z ) \exp( - 0.01 z^{0.8} ) \, d z,\\
        I_{1}(M)
        &= \int_{1}^{M} \cos( z ) \exp( - 0.01 z^{0.8} ) \, d z,\\
        I_{2}(M)
        &= \int_{1}^{M} \cos( z ) \exp( - 0.01 z^{0.8} ) w_{M}(z) \, d z.\\
    \end{aligned}
\end{equation*}
Define the errors of the two approaches by $E_{1}(M) = |I - I_{1}(M)|$ and $E_{2}(M) = |I - I_{2}(M)|$.
The results, shown in the referenced \Cref{tab:Dirac_windowing_function}, underscore the efficacy of incorporating the windowing function.
A more detailed evaluation, including comparisons among several different window functions, is provided in \Cref{appendix:window_function_comparison}.
\begin{table}[htp]
    \centering
    \resizebox{0.9\linewidth}{!}{
    \begin{tabularx}{\linewidth}{c|XXXXX}
        $M$ & $80$ & $160$ & $320$ & $640$ & $1280$\\
        \noalign{\hrule height 1pt}
        $E_{1}(M)$ & \textit{7.12097e-01} & \textit{1.24460e-01} & \textit{1.56845e-01} & \textit{1.33684e-01} & \textit{4.59293e-02}\\
        $E_{2}(M)$ & \textit{3.80618e-03} & \textit{1.14729e-04} & \textit{3.26691e-08} & \textit{5.06262e-14} & \textit{4.44089e-16}\\
    \end{tabularx}
    }
    \caption{
        An example showcasing the effectiveness of utilizing the windowing function, demonstrating that $E_{2}(M)$ decays spectrally in terms of $M$ while $E_{1}(M)$ decays very slowly comparably.
    }\label{tab:Dirac_windowing_function}
\end{table}

We now estimate the error of the far-field integration, i.e., how small $\tilde{g}_{\text{far}}$ is. For $y \neq 0$, the solution is given by the integrals in \Cref{eq:Dirac_general_FFPE_Bessel} and \Cref{eq:Dirac_general_FFPE_1d}.
Notice that the asymptotic behavior of the Bessel function $J_{\nu}(x)$ for large $x$ is given by $\sqrt{\frac{2}{\pi x}} \cos \left ( x - \frac{(2 \nu + 1) \pi}{4} \right )$. Therefore, analyzing the case when $f(z; y) = z^{l} e^{i y z}$ for $l \geq 0$ is sufficient and representative for our purposes.
To this end, we estimate the convergence of $\int_{\gamma M}^{\infty} u(z) (1 - w_{M, \gamma}(z)) e^{i k z} \, d z$ to zero as $M \to \infty$, where $u(z) = z^{l} \exp \left ( - \tau z^{q} \right )$ for $l \geq 0$ and $q \in (0, 2]$, with $k > 0$ and $\gamma \in (0, 1)$ fixed.
The estimate is given below while the proof is given in \Cref{sec:proofoftheorem}.
\begin{theorem}\label{thm:farfield}
    Let $u(z) = z^{l} \exp \left ( - \tau z^{q} \right )$ for $l \geq 0$ and $q \in (0, 2]$. Fix $k > 0$ and $\gamma \in (0, 1)$. Then for any $n \in \mathbb{N}$,
    \begin{equation*} 
        \begin{split}
        &\quad \, \left | \int_{\gamma M}^{\infty} u(z) (1 - w_{M, \gamma}(z)) e^{i k z} \, d z \right |\\
        &\leq
        \left \{
        \begin{aligned}
            &\frac{\mathcal{C}_{1}(n, l) M^{l}}{k^{n} (\gamma M)^{n - 1} \gamma} \min \left \{ \frac{1}{1 - \tau q M^{q}}, n + 1 \right \} \exp \left ( - \gamma^{q} \tau M^{q} \right )
            &&\text{for } \tau q M^{q} \leq 1,\\
            &\mathcal{C}_{2}(n, l) M^{l + 1} \left ( \frac{\tau q M^{q - 1}}{k \gamma} \right )^{n} \exp \left ( - \gamma^{q} \tau M^{q} \right )
            &&\text{for } \tau q M^{q} > 1.
        \end{aligned}
        \right .
        \end{split}
    \end{equation*}
    $\mathcal{C}_{1}(n, l)$ and $\mathcal{C}_{2}(n, l) $ are defined as $\mathcal{C}_{1}(n, l) = \sum_{j = 0}^{n} \binom{n}{j} \mathcal{A}(n - j, l) \mathcal{B}(j)$ and $\mathcal{C}_{2}(n, l) = (n + 1) \sum_{j = 0}^{n} \binom{n}{j} \mathcal{A}(n - j, l) \mathcal{B}(j) \left ( \tau q M^{q} \right )^{- j}$, where 
    \begin{align*}
        \mathcal{A}(n, l)
        := \max_{k \in \{ 0, 1, \cdots, n \}} \binom{n}{k} (l + n - 1)_{n - k}
        \quad \text{and} \quad
        \mathcal{B}(n)
        := \max_{x \in [1, 1 / \gamma]} \left | \frac{d^{n}}{d x^{n}} w_{\frac{1}{\gamma}, \gamma}(x) \right |.
    \end{align*}
    Here $(l)_{m}$ is the rising Pochhammer symbol defined as $(l)_{0} := 1$ and $(l)_m := l (l + 1) \cdots (l + m - 1)$ for $m \geq 1$.
\end{theorem}
The theorem has several implications.
First, In the regime where $\tau q M^{q} \leq 1$, the convergence rate is super-algebraic, meaning that the error decays at a rate faster than any polynomial rate of convergence.
In the regime where $\tau q M^{q} > 1$, the windowing function technique has a beneficial effect for $q \in (0, 1)$. This is due to the fact that $M^{q - 1} \to 0$ as $M \to \infty$, resulting in accelerated convergence when using the windowing function technique in this regime.
Furthermore, as the frequency $k$ increases, the convergence becomes faster.
The dependence of the convergence rate on the dimensionality of the problem is captured the term $M^{l}$.
For a $d$-dimensional problem, by \Cref{eq:Dirac_general_FFPE_Bessel}, we have $l = (d - 1) / 2$.
Consequently, as the dimension increases, the convergence rate becomes slower.
In practice, it is hard to determine the precise $M$ needed to achieves a given accuracy.
Therefore, we use a strategy of progressively doubling the value of $M$ until the difference between successive iterations becomes sufficiently small, or until a predefined maximum limit for $M$ is reached.
This approach is used in \Cref{alg:Dirac_integration_for_slow_decay}.

\begin{algorithm}[htbp]
    \caption{Integration for Slow Decay}\label{alg:Dirac_integration_for_slow_decay}
    \begin{algorithmic}[1]
        \Require A smooth oscillatory function denoted by $f$, a lower limit of integration labeled as $L \in (0, 1]$, a fractional exponent signified by $\alpha \in (0, 1)$, and a pseudo-temporal coefficient represented by $\tau > 0$
        \Ensure Computation of the integral $\displaystyle \int_{L}^{\infty} f(z; y) \exp \left ( - z^{2 \alpha} \tau \right ) \, d z$
        
        \State $M, M_{\text{max}}, \varepsilon \gets 80, 5120, \textit{1e-14}$
            \Comment{Adjustable}
        \State $I_{\text{previous}}, I_{\text{current}} \gets 0, \infty$
        \While{$|I_{\text{current}} - I_{\text{previous}}| > \varepsilon$ \textbf{and} $M \leq M_{\text{max}}$}
            \State $I_{\text{previous}} \gets I_{\text{current}}$
            \State $I_{\text{current}} \gets$ Apply quadrature to $\displaystyle \int_{L}^{M} f(z; y) \exp(- z^{2 \alpha} \tau) w_{M}(z) \, d z$
            \Statex \Comment{Utilize precomputed quadrature points and weights for consistent $L$, scaling the number of quadrature points with $M$}
            \State $M \gets 2 M$
        \EndWhile
        \If{$|I_{\text{current}} - I_{\text{previous}}| > \varepsilon$}
            \State \textbf{raise} a flag (without stopping)
        \EndIf
        \State \Return $I_{\text{current}}$
    \end{algorithmic}
\end{algorithm}

\subsection{The FFPE solver}\label{sec:FFPE_solver}
Integrating \Cref{alg:Dirac_integration_with_singularity} and \Cref{alg:Dirac_integration_for_slow_decay}, the foundational algorithms for calculating the solution to \Cref{eq:Dirac_general_FFPE_hd} are presented as \Cref{alg:Dirac_basic_solver_for_general_FFPE_with_zero_displacement} (zero displacement) and \Cref{alg:Dirac_basic_solver_for_general_FFPE} (nonzero displacement).
Incorporating the scaling law, we merge \Cref{alg:Dirac_basic_solver_for_general_FFPE_with_zero_displacement,alg:Dirac_basic_solver_for_general_FFPE} into the FFPE solver shown as \Cref{alg:Dirac_solver_for_general_FFPE_in_hd}.
The typical method of conducting the experiment involves invoking $\Call{\Cref{alg:Dirac_solver_for_general_FFPE_in_hd}}{y, t, \diffusionOrdinary, \diffusionFractional, \alpha}$, with potential initialization to store components that are repeatedly used under identical $\diffusionOrdinary, \diffusionFractional, \alpha$.

\begin{algorithm}[htbp]
    \caption{Fundamental Algorithm for the General $d$-Dimensional FFPE with Zero Displacement}\label{alg:Dirac_basic_solver_for_general_FFPE_with_zero_displacement}
    \begin{algorithmic}[1]
        \Require Temporal parameter $t > 0$, coefficients $\diffusionOrdinary \geq 0$, $\diffusionFractional > 0$, and a fractional exponent $\alpha \in (0, 1)$
        \Ensure Computation of $\tilde{p}(0, t)$ as the solution to \Cref{eq:Dirac_general_FFPE_hd} defined at \Cref{eq:Dirac_general_FFPE_zero_displacement}

        \State $f \gets \text{Function } r \mapsto \left ( \frac{r}{2 \pi} \right )^{d - 1} \exp \left ( - \diffusionOrdinary r^{2} t \right )$
        \State $L, \tau \gets 1, \diffusionFractional t$; $g \gets $ \Call{\Cref{alg:Dirac_integration_for_slow_decay}}{$f$, $L$, $\alpha$, $\tau$} $+$ \Call{\Cref{alg:Dirac_integration_with_singularity}}{$f$, $L$, $\alpha$, $\tau$}
        \State \Return $\frac{S_{d - 1}}{2 \pi} g$
    \end{algorithmic}
\end{algorithm}

\begin{algorithm}[htbp]
    \caption{Fundamental Algorithm for the General $d$-Dimensional FFPE with Nonzero Displacement}\label{alg:Dirac_basic_solver_for_general_FFPE}
    \begin{algorithmic}[1]
        \Require Magnitude of the displacement $y$, temporal parameter $t > 0$, coefficients $\diffusionOrdinary \geq 0$, $\diffusionFractional > 0$, and a fractional exponent $\alpha \in (0, 1)$
        \Ensure Computation of $\tilde{p}(y, t)$ as the solution to \Cref{eq:Dirac_general_FFPE_hd}

        \If{$d$ is $1$}
            \State $f \gets \text{Function } r \mapsto \cos \left ( y r \right ) \exp \left ( - \diffusionOrdinary r^{2} t \right )$
                \Comment{\Cref{eq:Dirac_general_FFPE_1d}}
        \Else
            \State $f \gets \text{Function } r \mapsto \left ( \frac{r}{2 \pi} \right )^{d / 2} J_{(d - 2) / 2} \left ( y r \right ) \exp \left ( - \diffusionOrdinary r^{2} t \right )$
                \Comment{\Cref{eq:Dirac_general_FFPE_Bessel}}
        \EndIf
        \State $L, \tau \gets 1, \diffusionFractional t$; $g \gets $ \Call{\Cref{alg:Dirac_integration_for_slow_decay}}{$f$, $L$, $\alpha$, $\tau$} $+$ \Call{\Cref{alg:Dirac_integration_with_singularity}}{$f$, $L$, $\alpha$, $\tau$}
        \State \Return $g / \pi$ \textbf{if} $d$ is $1$ \textbf{else} $g / y^{(d - 2) / 2}$
    \end{algorithmic}
\end{algorithm}

\begin{algorithm}
    \caption{Solver for the General $d$-Dimensional FFPE}\label{alg:Dirac_solver_for_general_FFPE_in_hd}
    \begin{algorithmic}[1]
        \Require Magnitude of the displacement $y$, temporal parameter $t > 0$, coefficients $\diffusionOrdinary \geq 0$, $\diffusionFractional > 0$, and a fractional exponent $\alpha \in (0, 1)$
        \Ensure Computation of $\tilde{p}(y, t)$ as the solution to \Cref{eq:Dirac_general_FFPE_hd}

        \State $y_{\text{limit}} \gets 10$
            \Comment{Adjustable threshold for applying force scaling}
        \If{$y$ is equal to $0$}
            \State $t^{\text{scaled}} \gets $ compute the quantity based on \Cref{sec:zero_displacement}
            \If{$t^{\text{scaled}} > t$}
                \State $\diffusionOrdinary^{\text{scaled}} \gets \left ( \frac{t}{t^{\text{scaled}}} \right )^{1 - \frac{1}{\alpha}} \diffusionOrdinary$
                \State \Return $\left ( \frac{t}{t^{\text{scaled}}} \right )^{- \frac{d}{2 \alpha}}$ $\times$ \Call{\Cref{alg:Dirac_basic_solver_for_general_FFPE_with_zero_displacement}}{$t^{\text{scaled}}$, $\diffusionOrdinary^{\text{scaled}}$, $\diffusionFractional$, $\alpha$}
            \Else
                \State \Return \Call{\Cref{alg:Dirac_basic_solver_for_general_FFPE_with_zero_displacement}}{$t$, $\diffusionOrdinary$, $\diffusionFractional$, $\alpha$}
            \EndIf
        \EndIf
        \If{$y \leq y_{\text{limit}}$}
            \State $p \gets $ \Call{\Cref{alg:Dirac_basic_solver_for_general_FFPE}}{$y$, $t$, $\diffusionOrdinary$, $\diffusionFractional$, $\alpha$}
        \EndIf
            \Comment{Applying force scaling when $y$ is large}
        \If{$y > y_{\text{limit}}$ \textbf{or} $|I_{\text{current}} - I_{\text{previous}}| > \varepsilon$ remains true in \Call{\Cref{alg:Dirac_basic_solver_for_general_FFPE}}{}}
            \State $y^{\text{scaled}} \gets \pi / 2$;
            $t^{\text{scaled}} \gets \left ( \frac{y^{\text{scaled}}}{y} \right )^{2 \alpha} t$;
            $\diffusionOrdinary^{\text{scaled}} \gets \left ( \frac{y^{\text{scaled}}}{y} \right )^{2 - 2 \alpha} \diffusionOrdinary$
            \State $p \gets$ $\left ( \frac{y^{\text{scaled}}}{y} \right )^{d}$ $\times$ \Call{\Cref{alg:Dirac_basic_solver_for_general_FFPE}}{$y^{\text{scaled}}$, $t^{\text{scaled}}$, $\diffusionOrdinary^{\text{scaled}}$, $\diffusionFractional$, $\alpha$}
        \EndIf
        \State \Return $p$
    \end{algorithmic}
\end{algorithm}

%%%%%%%%%%%%%%%%%%%%%%%%%%%Section 4 %%%%%%%%%%%%%%%%%%%%%%%%%%

\section{Numerical results}\label{sec:numerical_results}
In this section, we conduct numerical experiments to solve the FFPE using the solver described in \Cref{alg:Dirac_solver_for_general_FFPE_in_hd}.
We first study the qualitative behavior of solutions across different dimensions and varying diffusion coefficients.
Next, we validate the numerical accuracy of our method, using the explicit formulas of solution derived in \Cref{appendix:special_cases_half_alpha} (for $\alpha = 1 / 2$) and \Cref{appendix:special_cases_pure_FFPE_rational_alpha} (for $\diffusionOrdinary = 0$ and $\alpha = 1 / 3$) approximated by high-precision implementations.
Finally, we analyze the dependence of the running time on different factors of our algorithm when applied to different cases.

All experiments are conducted using MATLAB R2023a on a desktop equipped with an 11th Generation Intel\textsuperscript{\textregistered} Core\textsuperscript{\texttrademark} i7-11700F CPU and DDR4 2$\times$32GB 3600MHz memory.
Code is available at \texttt{https://github.com/ACMathX/FFPEDDIC}.

\subsection{Solution display}
In this subsection, we study the qualitative behavior of solutions.
Notice that the figures plotted use log-scales on the vertical axis.
\begin{figure}[htp]
    \centering
    \includegraphics[width=\textwidth]{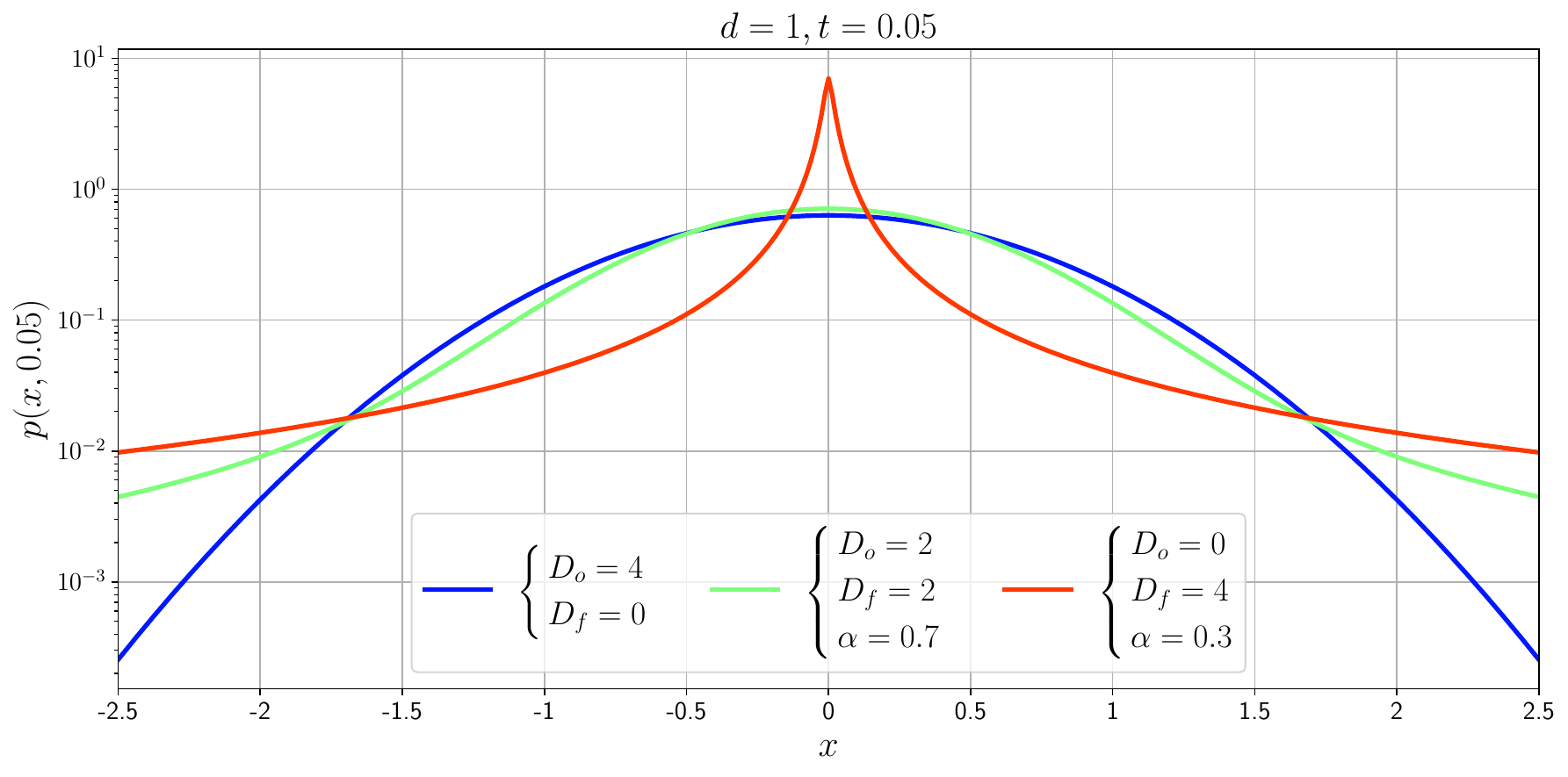}
    \caption{
        Solution representations in one dimension with $x_{0} = 0$ and $\drift = 0$ for $d = 1$ and $t = 0.05$.
        The plot shows $p(x, t)$ for different parameter settings: $\diffusionOrdinary = 4, \diffusionFractional = 0$ (pure ordinary diffusion, blue), $\diffusionOrdinary = 2, \diffusionFractional = 2, \alpha = 0.7$ (mixed diffusion, green), and $\diffusionOrdinary = 0, \diffusionFractional = 4, \alpha = 0.3$ (pure fractional diffusion, red).
        The graph illustrates how the density function $p(x, t)$ varies with $x$ under different diffusion conditions.
    }\label{fig:Dirac_solution_representation}
\end{figure}

\begin{figure}[htp]
    \centering
    \includegraphics[width=\textwidth]{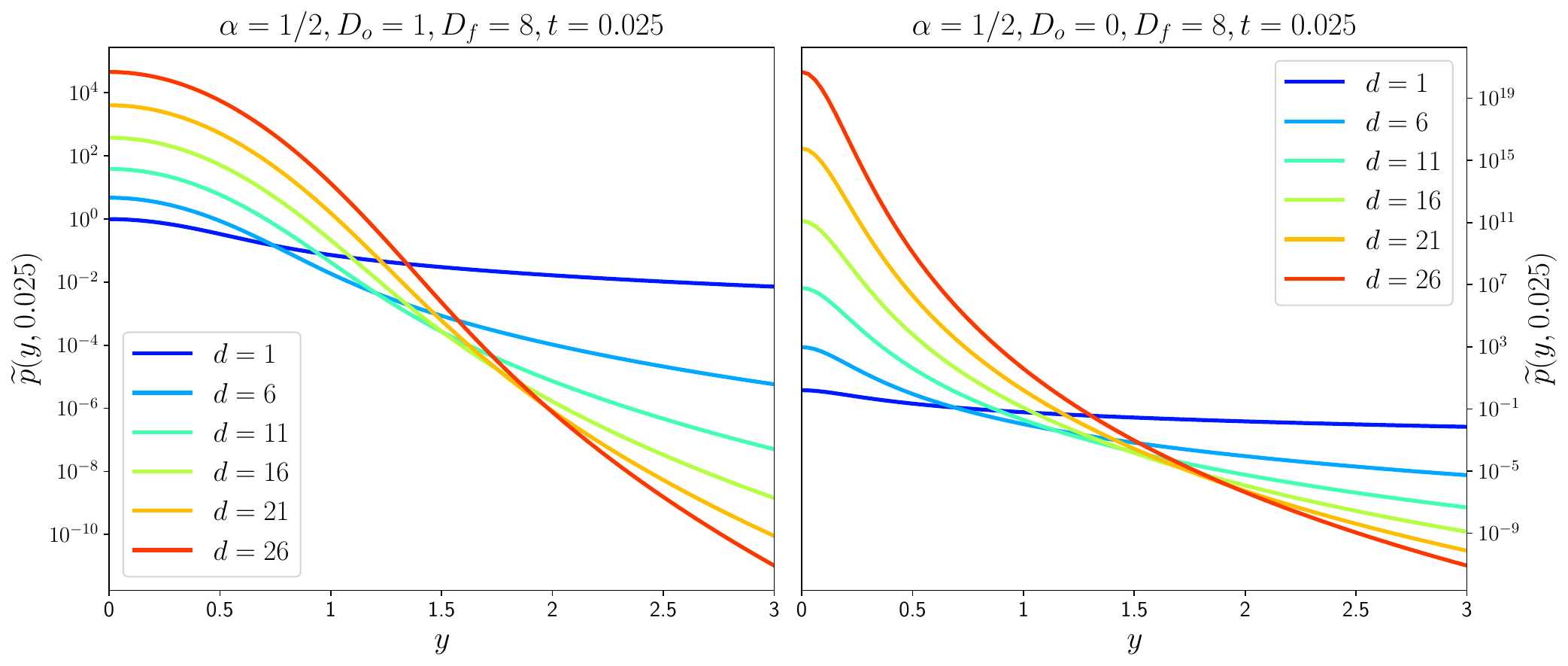}
    \caption{
        Special case examples for various values of $d$ and $y$ with $t = 0.025$.
        The left plot illustrates the behavior of $\tilde{p}(y, 0.025)$ with parameter settings $\alpha = 1 / 2$, $\diffusionOrdinary = 1$, and $\diffusionFractional = 8$ (mixed diffusion).
        The right plot shows $\tilde{p}(y, 0.025)$ with $\alpha = 1 / 2$, $\diffusionOrdinary = 0$, and $\diffusionFractional = 8$ (pure fractional diffusion).
        Different colors represent different values of $d$. The plots depict how $\tilde{p}(y, t)$ varies with $y$ for different dimensions $d$ under the given diffusion conditions.
    }\label{fig:Dirac_special_case_various_y}
\end{figure}

\begin{figure}[htp]
    \centering
    \includegraphics[width=\textwidth]{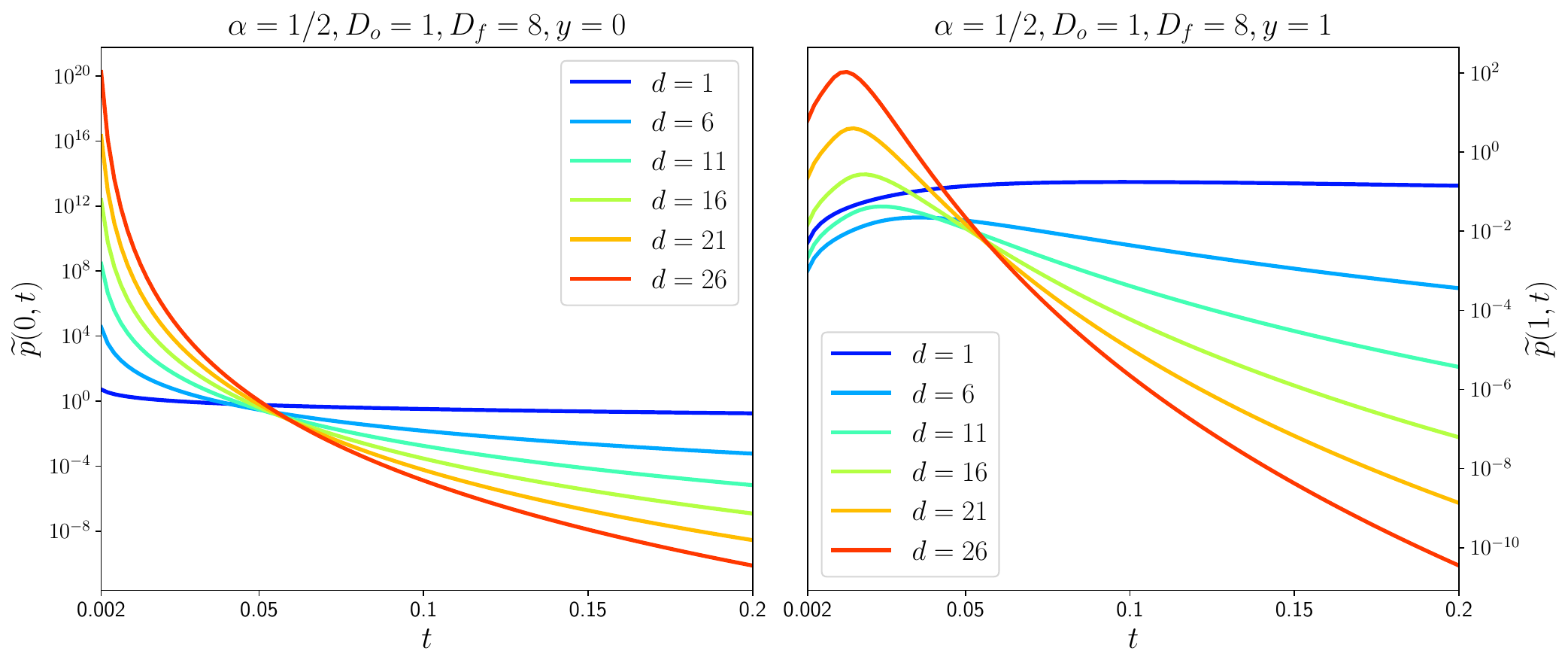}
    \caption{
        Special case examples for various values of $d$ and $t$ with $y = 0$ or $y = 1$.
        The left plot shows the behavior of $\tilde{p}(0, t)$ as $t$ varies, while the right plot shows the behavior of $\tilde{p}(1, t)$.
        The parameter settings are $\alpha = 1 / 2$, $\diffusionOrdinary = 1$, and $\diffusionFractional = 8$.
        Different colors represent different values of $d$.
        The plots illustrate the change in $\tilde{p}(y, t)$ over $t$ for these different dimensions, highlighting the differences in behavior for $y = 0$ and $y = 1$.
    }\label{fig:Dirac_special_case_various_t}
\end{figure}

\Cref{fig:Dirac_solution_representation} illustrates several representative curves for one-dimensional solutions, highlighting behaviors of the solutions with different diffusion coefficients.
Specifically, several distinctive characteristics of the FFPE are observed.
First, the density function for pure fractional diffusion exhibits heavier tails compared to the Gaussian distribution, indicating a higher probability of extreme values.
Additionally, the center of the density function for pure fractional diffusion is more concentrated and narrower than that of the Gaussian distribution.
These phenomena become more pronounced as the fractional diffusion coefficient increases or as the value of $\alpha \in (0, 1)$ decreases.
\Cref{fig:Dirac_special_case_various_y,fig:Dirac_special_case_various_t} show solutions for various $d$ with fixed $t$ or $y$ respectively under different parameters.
We observe that for smaller values of $t$ or larger dimensions $d$, the solution is closer to a Dirac-delta measure, making it more challenging to approximate accurately.

\subsection{Error analysis}
\Cref{table:mixed_diffusion_half_alpha,table:pure_fractional_diffusion_half_alpha,table:pure_fractional_diffusion_one_third_alpha} present the relationships between $d$, $y$, $t$, and the relative error for different parameters.
In these tables, the highlighted values in red indicate the cases where the numerical algorithm is not effective.
Before hitting the machine precision, the relative error generally grows with increasing $d$ and decreases with increasing $t$.
Pure fractional cases, particularly when $\alpha \in (0, 1)$ is smaller, are more challenging to approximate.
Consequently, the error tables exhibit larger relative errors for these cases.

\begin{table}[htp]
    \centering
    \resizebox{0.92\linewidth}{!}{%
    \begin{tabular}{|c|c|c|c|c|c|c|c|c|}
        \hline
        \diagbox{$t$}{$d$} & 1 & 5 & 9 & 13 & 17 & 21 & 25 & 29\\
        \hline
        0.004 & \textit{1.59e-13} & \textit{2.40e-11} & \textit{4.08e-10} & \textit{1.45e-08} & \textit{9.80e-08} & \textit{6.73e-07} & \textit{6.32e-06} & \textit{1.11e-04}\\
        \hline
        0.020 & \textit{8.13e-15} & \textit{3.93e-13} & \textit{1.75e-12} & \textit{2.55e-12} & \textit{8.39e-13} & \textit{3.34e-12} & \textit{3.32e-12} & \textit{1.35e-11}\\
        \hline
        0.040 & \textit{2.27e-15} & \textit{2.41e-14} & \textit{8.46e-14} & \textit{1.46e-13} & \textit{5.22e-14} & \textit{8.37e-14} & \textit{1.92e-14} & \textit{3.06e-14}\\
        \hline
        0.060 & \textit{1.52e-15} & \textit{3.83e-15} & \textit{5.13e-15} & \textit{8.80e-15} & \textit{6.39e-15} & \textit{1.13e-14} & \textit{6.84e-15} & \textit{8.71e-15}\\
        \hline
        0.080 & \textit{9.79e-16} & \textit{1.37e-15} & \textit{9.02e-16} & \textit{2.13e-15} & \textit{2.27e-15} & \textit{1.84e-15} & \textit{3.74e-15} & \textit{6.30e-15}\\
        \hline
        0.100 & \textit{4.37e-16} & \textit{9.90e-16} & \textit{1.44e-15} & \textit{2.10e-15} & \textit{1.99e-15} & \textit{2.23e-15} & \textit{3.97e-15} & \textit{4.72e-15}\\
        \hline
        0.120 & \textit{4.11e-16} & \textit{1.18e-15} & \textit{5.31e-16} & \textit{1.31e-15} & \textit{1.22e-15} & \textit{1.48e-15} & \textit{2.46e-15} & \textit{4.83e-15}\\
        \hline
        0.140 & \textit{3.54e-16} & \textit{6.72e-16} & \textit{1.13e-15} & \textit{1.51e-15} & \textit{2.62e-15} & \textit{2.60e-15} & \textit{4.95e-15} & \textit{2.22e-15}\\
        \hline
        0.160 & \textit{1.13e-15} & \textit{8.14e-16} & \textit{7.56e-16} & \textit{1.29e-15} & \textit{1.87e-15} & \textit{1.63e-15} & \textit{3.76e-15} & \textit{2.41e-15}\\
        \hline
        0.180 & \textit{5.28e-16} & \textit{7.86e-16} & \textit{7.21e-16} & \textit{1.40e-15} & \textit{9.22e-16} & \textit{1.26e-15} & \textit{2.58e-15} & \textit{3.56e-15}\\
        \hline
        0.200 & \textit{3.79e-16} & \textit{4.78e-16} & \textit{9.96e-16} & \textit{7.62e-16} & \textit{1.79e-15} & \textit{2.15e-15} & \textit{4.05e-15} & \textit{1.85e-15}\\
        \hline
    \end{tabular}
    }
    \caption{
        Maximum relative errors for $y \in [0, 2]$ for various $t$ and $d$.
        Parameters are set as $\diffusionOrdinary = 1, \diffusionFractional = 8, \alpha = 1 / 2$.
    }
    \label{table:mixed_diffusion_half_alpha}
\end{table}

\begin{table}[htp]
    \centering
    \resizebox{0.92\linewidth}{!}{%
    \begin{tabular}{|c|c|c|c|c|c|c|c|c|}
        \hline
        \diagbox{$t$}{$d$} & 1 & 5 & 9 & 13 & 17 & 21 & 25 & 29\\
        \hline 
        0.004 & \textit{7.09e-13} & \textit{2.21e-09} & \textit{2.23e-06} & \textit{3.23e-03} & \HIGHLIGHTA{\textit{1.65e+00}} & \HIGHLIGHTA{\textit{5.41e+02}} & \HIGHLIGHTA{\textit{2.27e+06}} & \HIGHLIGHTA{\textit{3.13e+09}}\\
        \hline
        0.020 & \textit{2.03e-14} & \textit{2.44e-12} & \textit{4.58e-11} & \textit{1.43e-08} & \textit{6.60e-07} & \textit{2.34e-05} & \textit{7.72e-04} & \textit{2.12e-02}\\
        \hline
        0.040 & \textit{4.60e-15} & \textit{3.18e-13} & \textit{1.70e-12} & \textit{1.11e-11} & \textit{1.38e-10} & \textit{6.61e-10} & \textit{1.51e-08} & \textit{3.59e-08}\\
        \hline
        0.060 & \textit{2.43e-15} & \textit{4.74e-14} & \textit{3.05e-13} & \textit{6.14e-13} & \textit{1.96e-12} & \textit{6.94e-12} & \textit{2.10e-10} & \textit{1.25e-09}\\
        \hline
        0.080 & \textit{9.95e-16} & \textit{1.39e-14} & \textit{6.86e-14} & \textit{1.42e-13} & \textit{7.80e-14} & \textit{3.43e-13} & \textit{1.85e-12} & \textit{7.07e-12}\\
        \hline
        0.100 & \textit{8.36e-16} & \textit{4.82e-15} & \textit{2.17e-14} & \textit{5.05e-14} & \textit{1.91e-14} & \textit{7.37e-14} & \textit{8.38e-14} & \textit{1.94e-13}\\
        \hline
        0.120 & \textit{5.11e-16} & \textit{1.98e-15} & \textit{7.02e-15} & \textit{1.91e-14} & \textit{1.02e-14} & \textit{3.04e-14} & \textit{1.65e-14} & \textit{2.20e-14}\\
        \hline
        0.140 & \textit{4.52e-16} & \textit{1.11e-15} & \textit{2.39e-15} & \textit{6.53e-15} & \textit{5.70e-15} & \textit{1.59e-14} & \textit{9.07e-15} & \textit{1.72e-14}\\
        \hline
        0.160 & \textit{1.06e-15} & \textit{1.37e-15} & \textit{7.99e-16} & \textit{2.34e-15} & \textit{3.44e-15} & \textit{8.19e-15} & \textit{7.07e-15} & \textit{1.10e-14}\\
        \hline
        0.180 & \textit{4.51e-16} & \textit{9.11e-16} & \textit{8.11e-16} & \textit{1.68e-15} & \textit{1.30e-15} & \textit{4.81e-15} & \textit{6.56e-15} & \textit{9.69e-15}\\
        \hline
        0.200 & \textit{4.67e-16} & \textit{7.73e-16} & \textit{9.59e-16} & \textit{1.72e-15} & \textit{1.88e-15} & \textit{2.82e-15} & \textit{5.49e-15} & \textit{7.12e-15}\\
        \hline
    \end{tabular}
    }
    \caption{
        Maximum relative errors for $y \in [0, 2]$ for various $t$ and $d$.
        Large quantities are highlighted in red.
        Parameters are set as $\diffusionOrdinary = 0, \diffusionFractional = 8, \alpha = 1 / 2$.
    }
    \label{table:pure_fractional_diffusion_half_alpha}
\end{table}

\begin{table}[htp]
    \centering
    \resizebox{0.92\linewidth}{!}{%
    \begin{tabular}{|c|c|c|c|c|c|c|c|c|}
        \hline
        \diagbox{$t$}{$d$} & 1 & 5 & 9 & 13 & 17 & 21 & 25 & 29\\
        \hline
        0.004 & \textit{3.85e-12} & \textit{6.40e-07} & \textit{9.04e-01} & \HIGHLIGHTA{\textit{1.80e+05}} & \HIGHLIGHTA{\textit{2.56e+10}} & \HIGHLIGHTA{\textit{2.88e+15}} & \HIGHLIGHTA{\textit{2.33e+20}} & \HIGHLIGHTA{\textit{1.35e+25}}\\
        \hline
        0.020 & \textit{1.72e-13} & \textit{2.01e-10} & \textit{5.13e-07} & \textit{1.46e-02} & \HIGHLIGHTA{\textit{1.31e+02}} & \HIGHLIGHTA{\textit{7.66e+05}} & \HIGHLIGHTA{\textit{1.45e+09}} & \HIGHLIGHTA{\textit{1.29e+14}}\\
        \hline
        0.040 & \textit{3.37e-14} & \textit{1.06e-11} & \textit{3.14e-09} & \textit{7.47e-06} & \textit{4.07e-03} & \HIGHLIGHTA{\textit{2.70e+00}} & \HIGHLIGHTA{\textit{9.86e+03}} & \HIGHLIGHTA{\textit{1.75e+07}}\\
        \hline
        0.060 & \textit{1.02e-14} & \textit{1.15e-12} & \textit{6.31e-10} & \textit{3.96e-08} & \textit{2.45e-05} & \textit{1.15e-02} & \HIGHLIGHTA{\textit{2.97e+00}} & \HIGHLIGHTA{\textit{1.86e+02}}\\
        \hline
        0.080 & \textit{2.01e-15} & \textit{3.90e-13} & \textit{1.71e-11} & \textit{7.22e-09} & \textit{3.00e-07} & \textit{5.62e-05} & \textit{1.97e-02} & \HIGHLIGHTA{\textit{3.97e+00}}\\
        \hline
        0.100 & \textit{1.83e-15} & \textit{1.74e-13} & \textit{2.24e-12} & \textit{7.16e-10} & \textit{4.34e-08} & \textit{1.42e-06} & \textit{9.65e-05} & \textit{2.89e-02}\\
        \hline
        0.120 & \textit{1.22e-15} & \textit{8.93e-14} & \textit{8.56e-13} & \textit{1.93e-11} & \textit{4.42e-09} & \textit{1.76e-07} & \textit{4.88e-06} & \textit{1.08e-04}\\
        \hline
        0.140 & \textit{9.93e-16} & \textit{4.85e-14} & \textit{3.91e-13} & \textit{2.79e-12} & \textit{3.84e-10} & \textit{1.73e-08} & \textit{5.10e-07} & \textit{1.22e-05}\\
        \hline
        0.160 & \textit{1.32e-15} & \textit{2.86e-14} & \textit{1.86e-13} & \textit{1.21e-12} & \textit{2.22e-11} & \textit{1.86e-09} & \textit{4.63e-08} & \textit{1.09e-06}\\
        \hline
        0.180 & \textit{5.71e-16} & \textit{1.76e-14} & \textit{1.09e-13} & \textit{4.87e-13} & \textit{4.66e-12} & \textit{1.84e-10} & \textit{3.44e-09} & \textit{7.90e-08}\\
        \hline
        0.200 & \textit{5.26e-16} & \textit{7.08e-15} & \textit{6.62e-14} & \textit{2.28e-13} & \textit{1.48e-12} & \textit{1.20e-11} & \textit{4.09e-10} & \textit{5.07e-09}\\
        \hline
    \end{tabular}
    }
    \caption{
        Maximum relative errors for $y \in [0, 2]$ for various $t$ and $d$.
        Large quantities are highlighted in red.
        Parameters are set as $\diffusionOrdinary = 0, \diffusionFractional = 8, \alpha = 1 / 3$.
    }
    \label{table:pure_fractional_diffusion_one_third_alpha}
\end{table}

\subsection{Running time study}
\begin{figure}[htp]
    \centering
    \includegraphics[width=\textwidth]{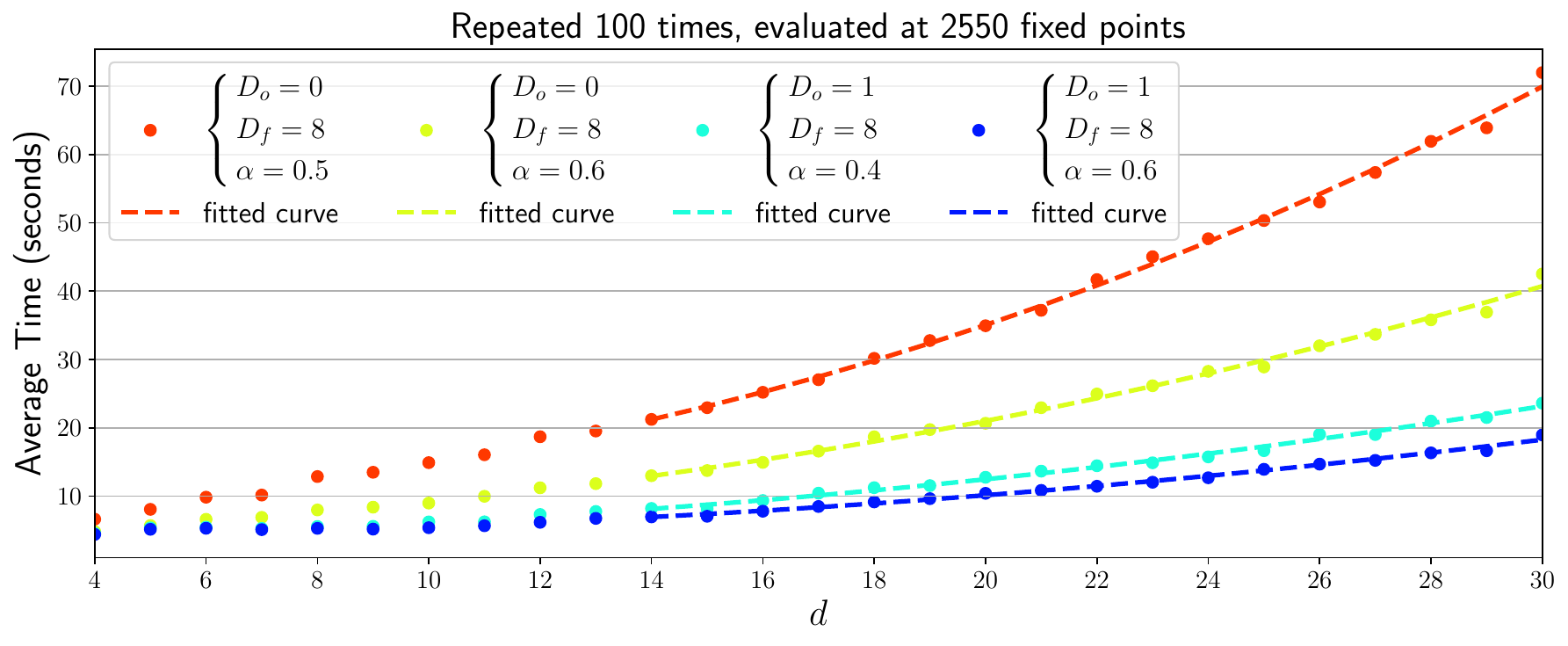}
    \caption{
        Average running time with various coefficients in high-dimensional cases.
        The running time is measured for computing all 2550 evaluations, repeated 100 times at fixed points that are evenly distributed over the space $\{ (y, t) | y \in [0, 2], t \in (0, 0.2] \}$.
        The different colors represent different coefficient sets as indicated in the legend.
        The fitted curves are quadratic, providing a representation of the growth trend in average time as the dimensionality $d$ increases.
    }\label{fig:Dirac_time_complexity}
\end{figure}

\begin{figure}[htp]
    \centering
    \includegraphics[width=\textwidth]{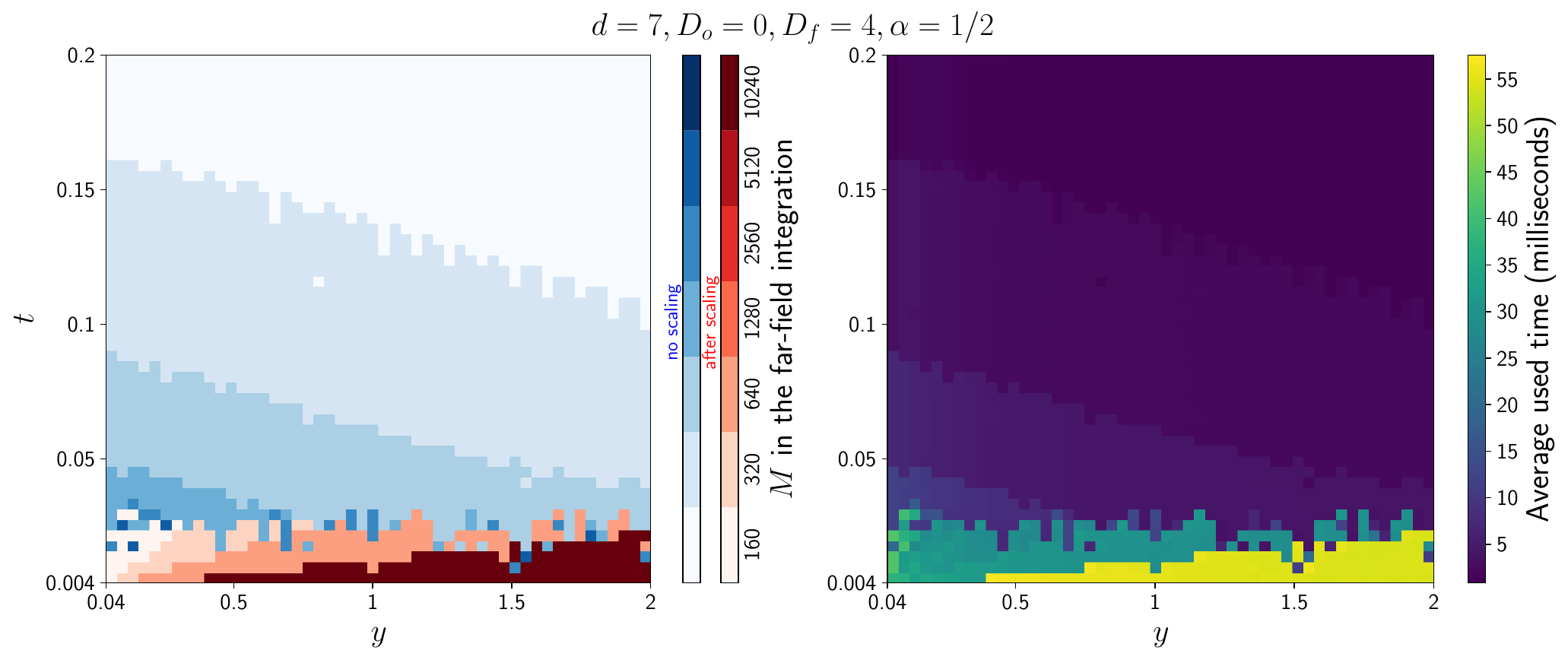}
    \caption{
        Average used time for the 2500 evaluations, which are evenly distributed over the space (excluding the 50 evaluations at $y = 0$), comparing the chosen $M$ for the far-field integration.
        The used time is measured 100 times, and the displayed values are the averages of these 100 tests.
        The left plot shows the $M$ landscape with two color bars: the blue bar represents $M$ values with direct convergence, and the red bar represents $M$ values where scaling was applied due to unsatisfied convergence criteria, followed by a re-evaluation to obtain the new $M$ after scaling.
        The smallest value of $M$ is $160$ because $M$ is doubled at least once in \Cref{alg:Dirac_integration_for_slow_decay}.
        The largest value of $M$ is $10240$, determined by the stopping criterion $M > M_{\text{max}}$ in \Cref{alg:Dirac_integration_for_slow_decay}.
        Note that $y$ is not large enough to engage the force scaling in this context.
        The right plot depicts the average used time in milliseconds for the evaluations.
    }\label{fig:Dirac_M_t_landscape}
\end{figure}

In this subsection, we study the computational cost of our algorithm across various dimensions and diffusion parameters by examining the average running time.
As \Cref{thm:farfield} indicates, the dimension $d$ influences the convergence rate of the far-field integration at $M \to \infty$.
With larger $d$, the convergence is expected to be slower.
\Cref{fig:Dirac_time_complexity} encapsulates the growth trend in average running time as the dimension $d$ increases.
While the increase in time with respect to $d$ appears to be quadratic, readers should also consult \Cref{table:mixed_diffusion_half_alpha,table:pure_fractional_diffusion_half_alpha,table:pure_fractional_diffusion_one_third_alpha} to compare the accuracy of the approximation for different cases.
As \Cref{fig:Dirac_M_t_landscape} indicates, the primary cost of our computation lies in the far-field integration approximation, where the selection of $M$ dominates the time consumption.

\subsection{Extension to Sum-of-Gaussians Initial Conditions}\label{sec:extension_to_sum_of_Gaussians_IC}
In this subsection, we extend our method to handle initial conditions given by sums of Gaussians.
More generally, the method can be applied to cases where the Fourier transform of the initial condition is easily computable, such as in the case of Gaussians or other radial basis functions with known Fourier representations \cite{Buhmann_2003}.
Since general functions can often be approximated by radial basis functions, our approach provides a practical means for computing solutions in such cases.

To illustrate this point, we first consider the initial condition given by a Gaussian centered at $\boldsymbol{x}_{0} \in \mathbb{R}^{d}$ with variance $\sigma^{2}$:
\begin{equation*}
    \left \{
        \begin{aligned}
            \frac{\partial p}{\partial t}(\boldsymbol{x}, t)
            &= - \boldsymbol{b} \cdot \nabla p(\boldsymbol{x}, t)
            + D_{\textnormal{o}} \Delta p(\boldsymbol{x}, t)
            - D_{\textnormal{f}} (- \Delta)^{\alpha} p(\boldsymbol{x}, t)\\
            p(\boldsymbol{x}, 0)
            &= \frac{1}{\left ( \sqrt{2 \pi \sigma^{2}} \right )^{d}} \exp \left ( - \frac{|\boldsymbol{x} - \boldsymbol{x}_{0}|^{2}}{2 \sigma^{2}} \right )
        \end{aligned}
    \right .
    .
\end{equation*}
Apply the Fourier transform on both sides, we have
\begin{equation*}
    \left \{
        \begin{aligned}
            \hat{p}_{t}(\boldsymbol{k}, t)
            &= - i \langle \boldsymbol{k}, \boldsymbol{b} \rangle \hat{p}(\boldsymbol{k}, t)
            - D_{\textnormal{o}} |\boldsymbol{k}|^{2} \hat{p}(\boldsymbol{k}, t)
            - D_{\textnormal{f}} |\boldsymbol{k}|^{2 \alpha} \hat{p}(\boldsymbol{k}, t)\\
            \hat{p}(\boldsymbol{k}, 0)
            &= \exp \left ( - i \langle \boldsymbol{k}, \boldsymbol{x}_{0} \rangle \right )
            \exp \left ( - \frac{\sigma^{2} |\boldsymbol{k}|^{2}}{2} \right )
        \end{aligned}
    \right .
    ,
\end{equation*}
which implies that
\begin{equation*}
    \begin{aligned}
        \hat{p}(\boldsymbol{k}, t)
        &= \exp \left ( - i \langle \boldsymbol{k}, \boldsymbol{x}_{0} + \boldsymbol{b} t \rangle \right )
        \exp \left ( - \frac{\sigma^{2} |\boldsymbol{k}|^{2}}{2} \right )
        \exp \left ( - \left ( D_{\textnormal{o}} |\boldsymbol{k}|^{2} + D_{\textnormal{f}} |\boldsymbol{k}|^{2 \alpha} \right ) t \right )\\
        &= \exp \left ( - i \langle \boldsymbol{k}, \boldsymbol{x}_{0} + \boldsymbol{b} t \rangle \right )
        \exp \left ( - \left ( \left [ D_{\textnormal{o}} + \frac{\sigma^{2}}{2 t} \right ] |\boldsymbol{k}|^{2} + D_{\textnormal{f}} |\boldsymbol{k}|^{2 \alpha} \right ) t \right ).
    \end{aligned}
\end{equation*}
In this case, the Gaussian initial condition leads to applying our original algorithm with a modified diffusion coefficient $D_{o} \mapsto D_{o} + \frac{\sigma^2}{2 t}$.

More generally, if the initial condition is a weighted sum of Gaussians,
\begin{equation*}
    p(\boldsymbol{x}, 0)
    = \sum_{j = 1}^{N} \frac{w_{j}}{\left ( \sqrt{2 \pi \sigma_{j}^{2}} \right )^{d}} \exp \left ( - \frac{|\boldsymbol{x} - \boldsymbol{x}_{j}|^{2}}{2 \sigma_{j}^{2}} \right ),
\end{equation*}
where $w_{j} > 0$, $\boldsymbol{x}_{j} \in \mathbb{R}^{d}$, $\sigma_{j} > 0$ for all $j$ and $\displaystyle \sum_{j = 1}^{N} w_{j} = 1$, denote $p_{j} \left ( \boldsymbol{x}, t; D_{o} + \frac{\sigma_{j}^{2}}{2 t} \right )$ as the solution to the equation with initial condition
\begin{equation*}
    p(\boldsymbol{x}, 0)
    = \frac{1}{\left ( \sqrt{2 \pi \sigma_{j}^{2}} \right )^{d}} \exp \left ( - \frac{|\boldsymbol{x} - \boldsymbol{x}_{j}|^{2}}{2 \sigma_{j}^{2}} \right ),
\end{equation*}
then the solution admits the decomposition
\begin{equation*}
    p(\boldsymbol{x}, t)
    = \sum_{j = 1}^{N} w_{j} p_{j} \left ( \boldsymbol{x}, t; D_{o} + \frac{\sigma_{j}^{2}}{2 t} \right ).
\end{equation*}
This linearity and decomposition allow us to efficiently compute solutions for more general initial data using the solver developed for Dirac-delta initial conditions.

\begin{figure}[htp]
    \centering
    \hfill
    \includegraphics[width=0.24\linewidth]{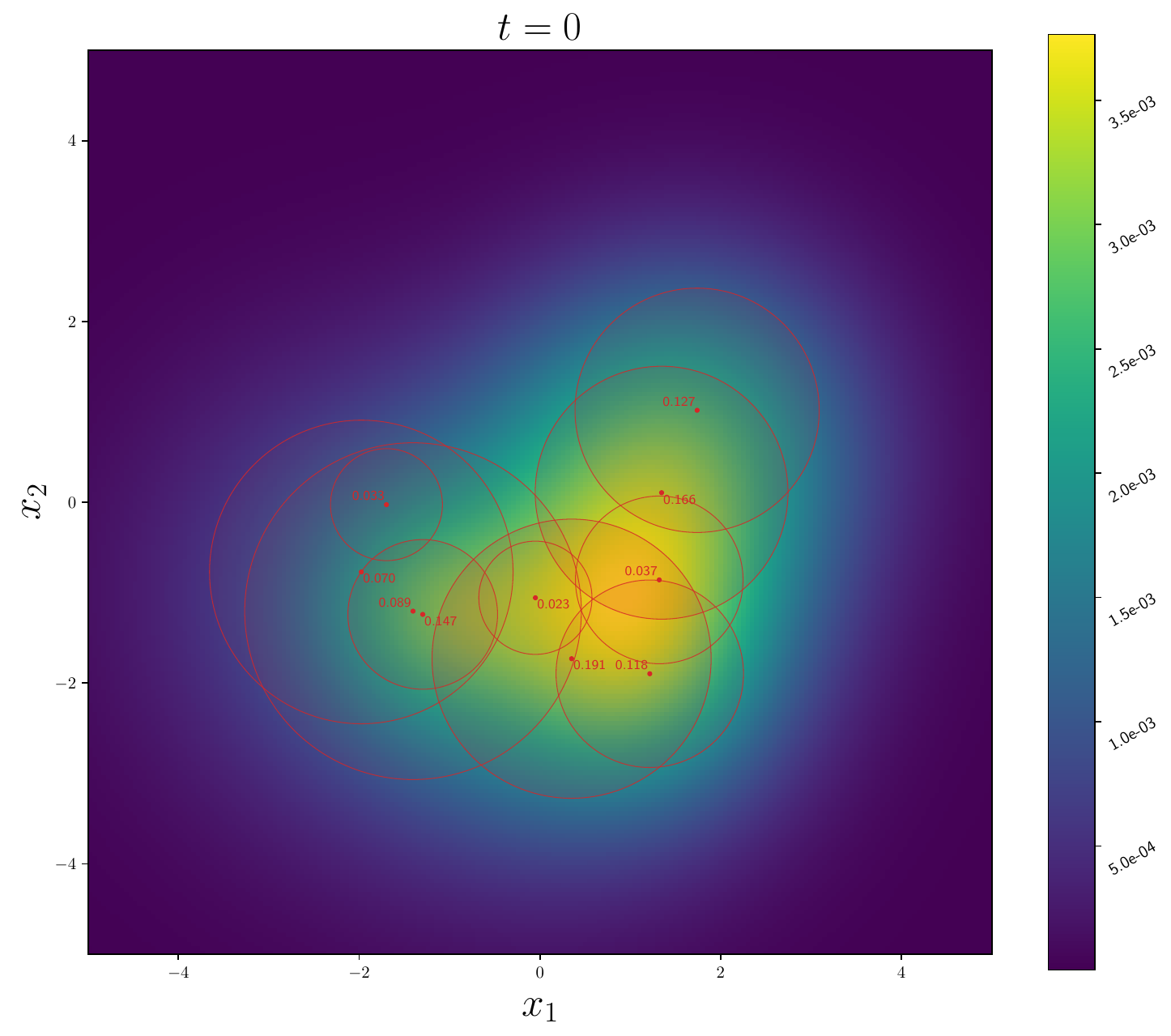}
    \hfill
    \includegraphics[width=0.24\linewidth]{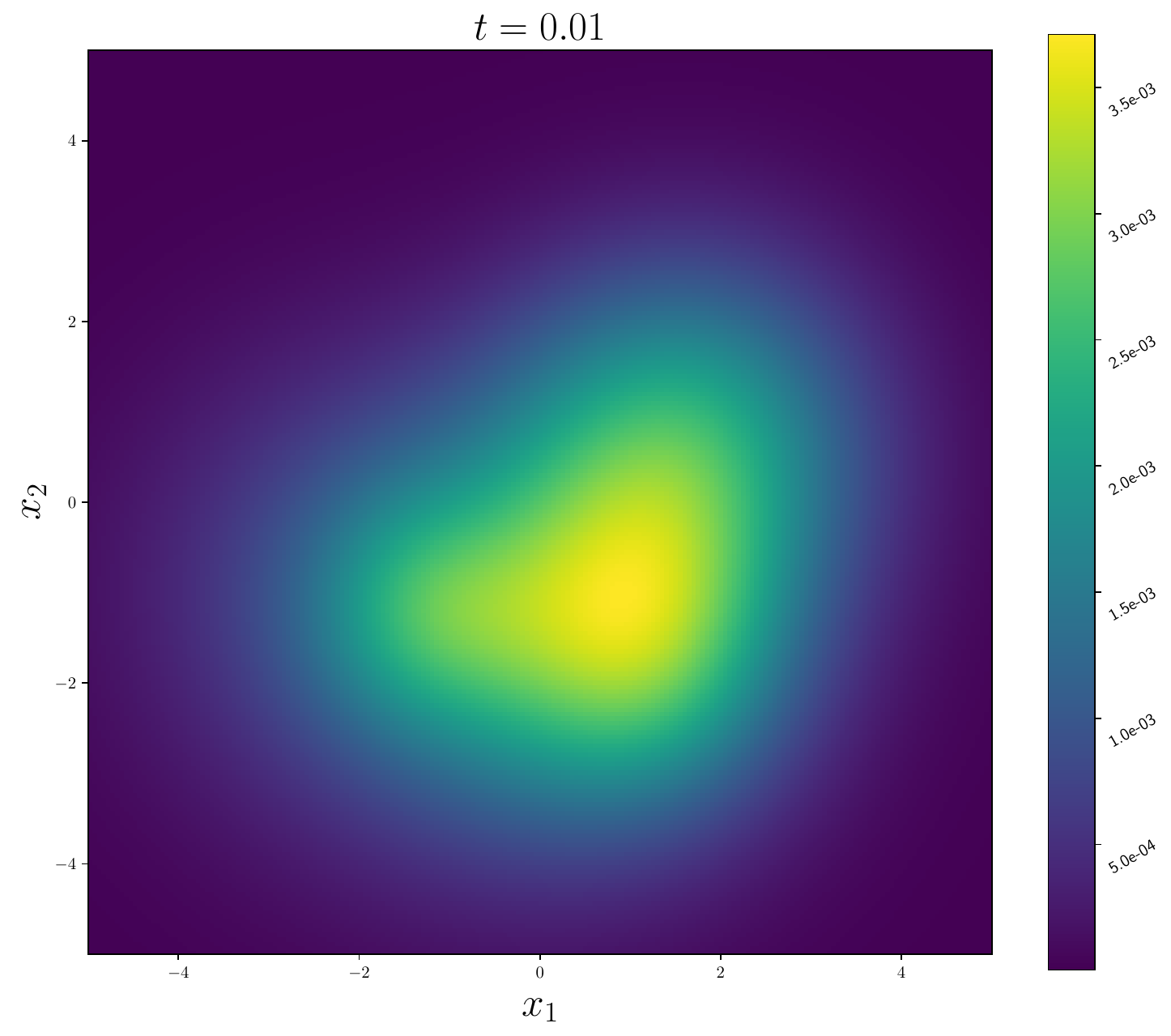}
    \hfill
    \includegraphics[width=0.24\linewidth]{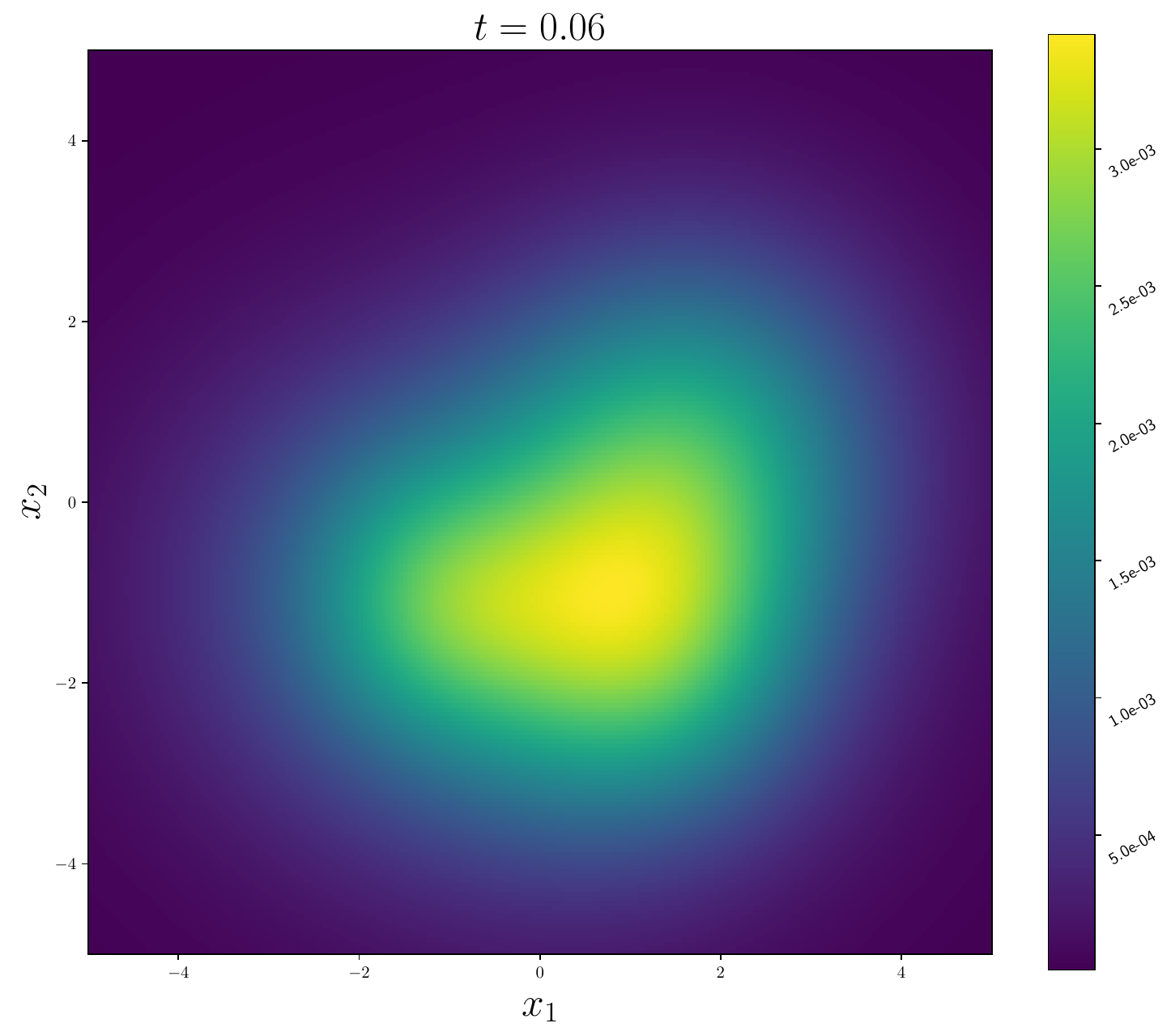}
    \hfill
    \includegraphics[width=0.24\linewidth]{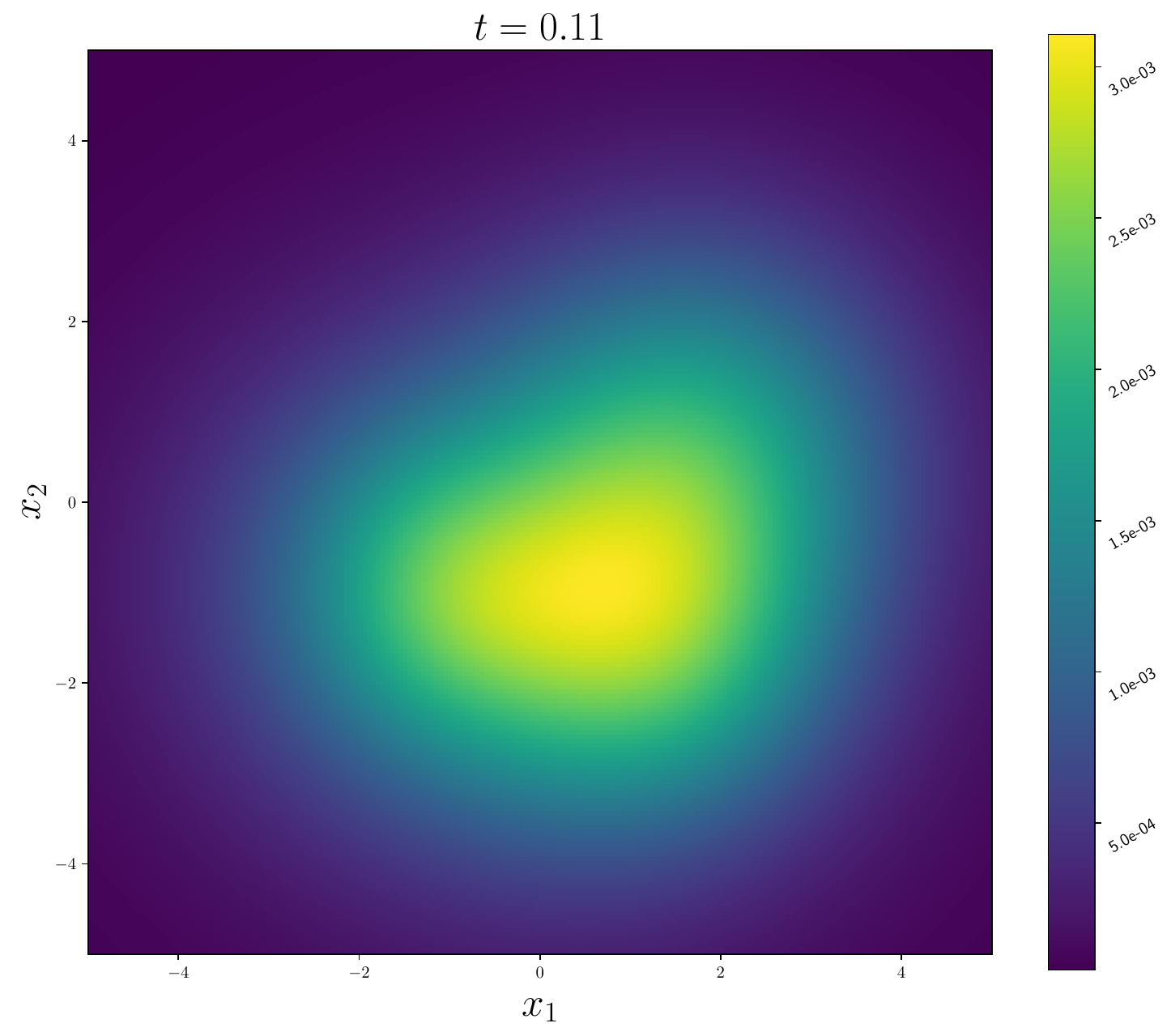}
    \hfill\\
    \hfill
    \includegraphics[width=0.24\linewidth]{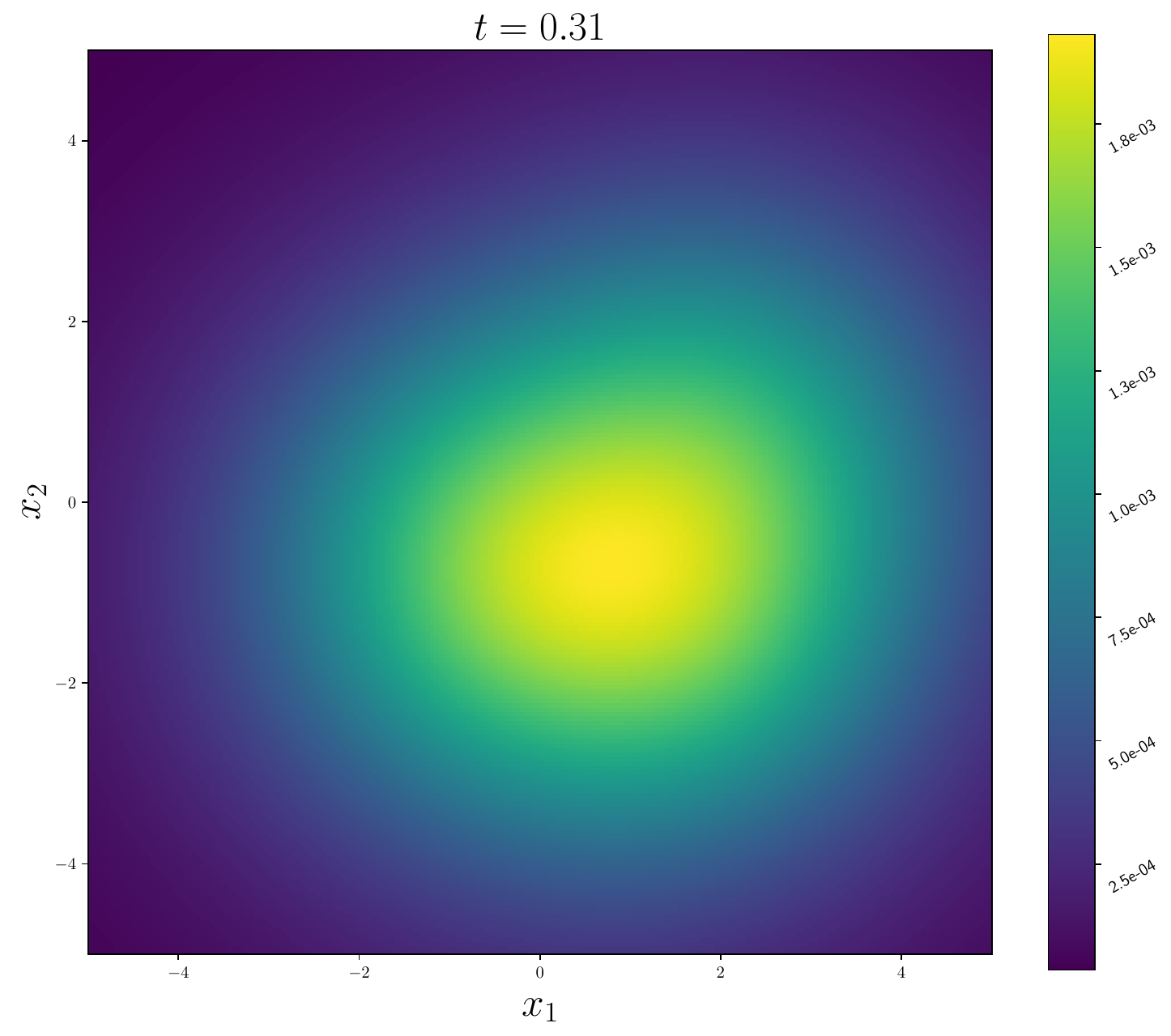}
    \hfill
    \includegraphics[width=0.24\linewidth]{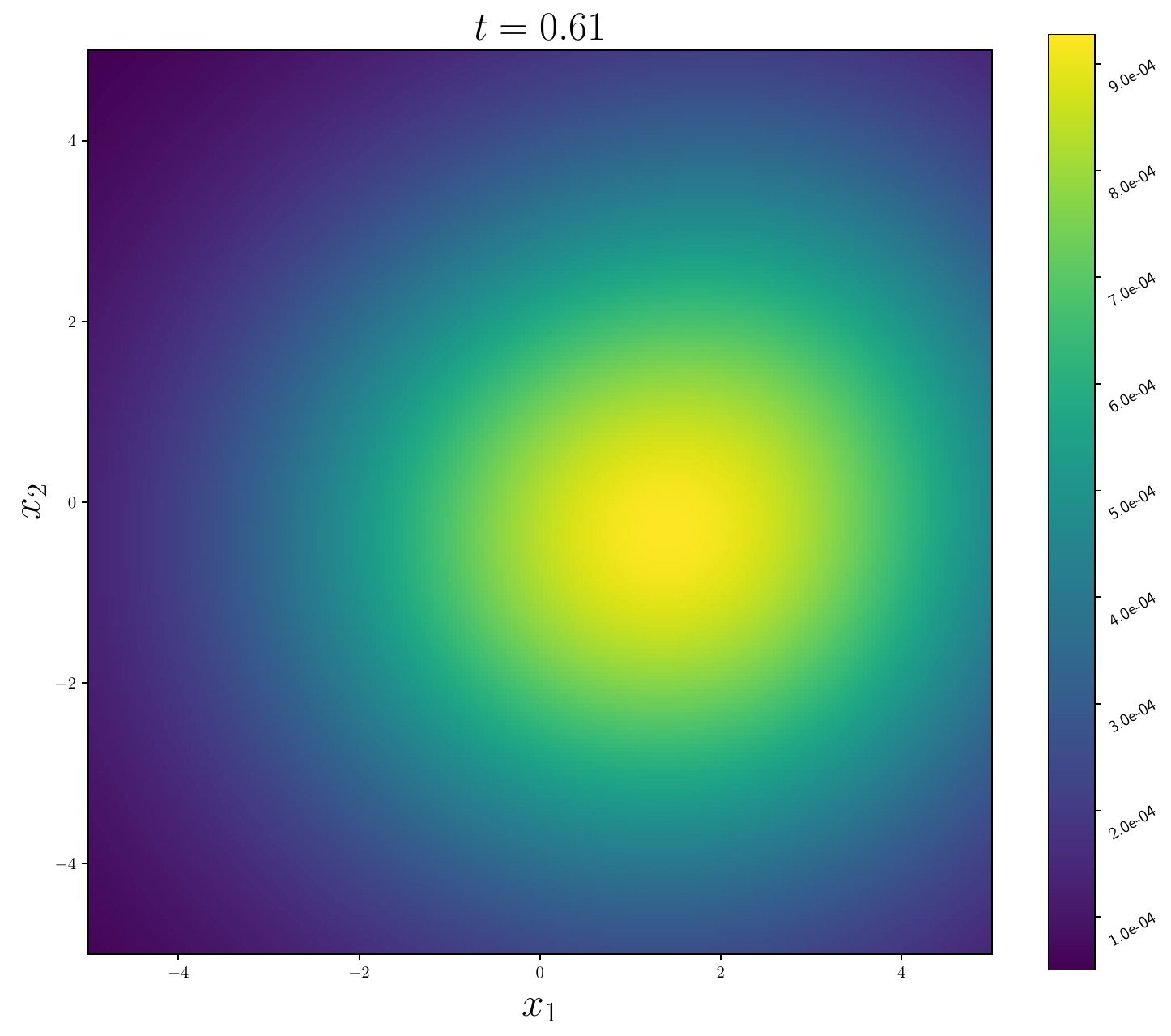}
    \hfill
    \includegraphics[width=0.24\linewidth]{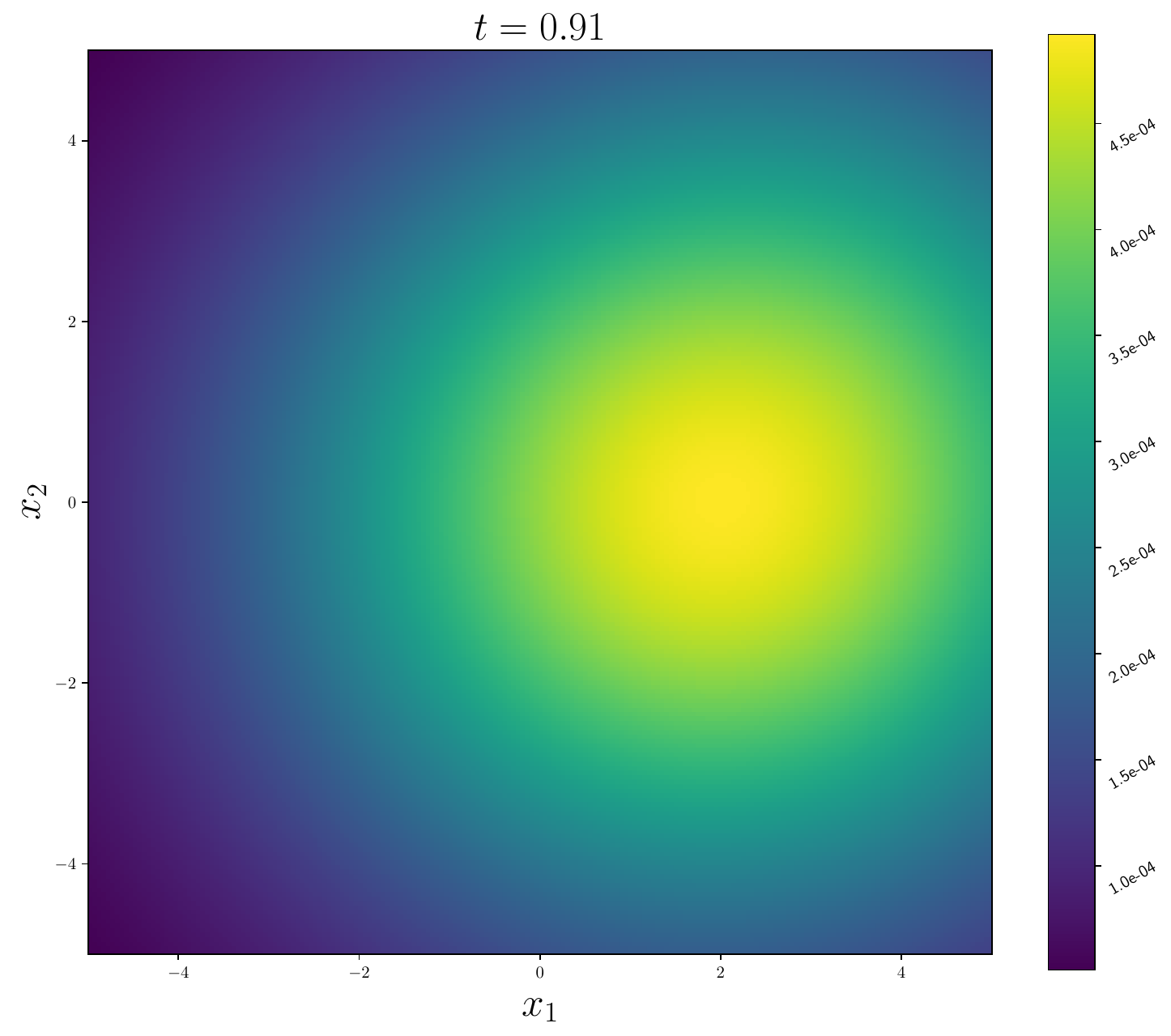}
    \hfill
    \includegraphics[width=0.24\linewidth]{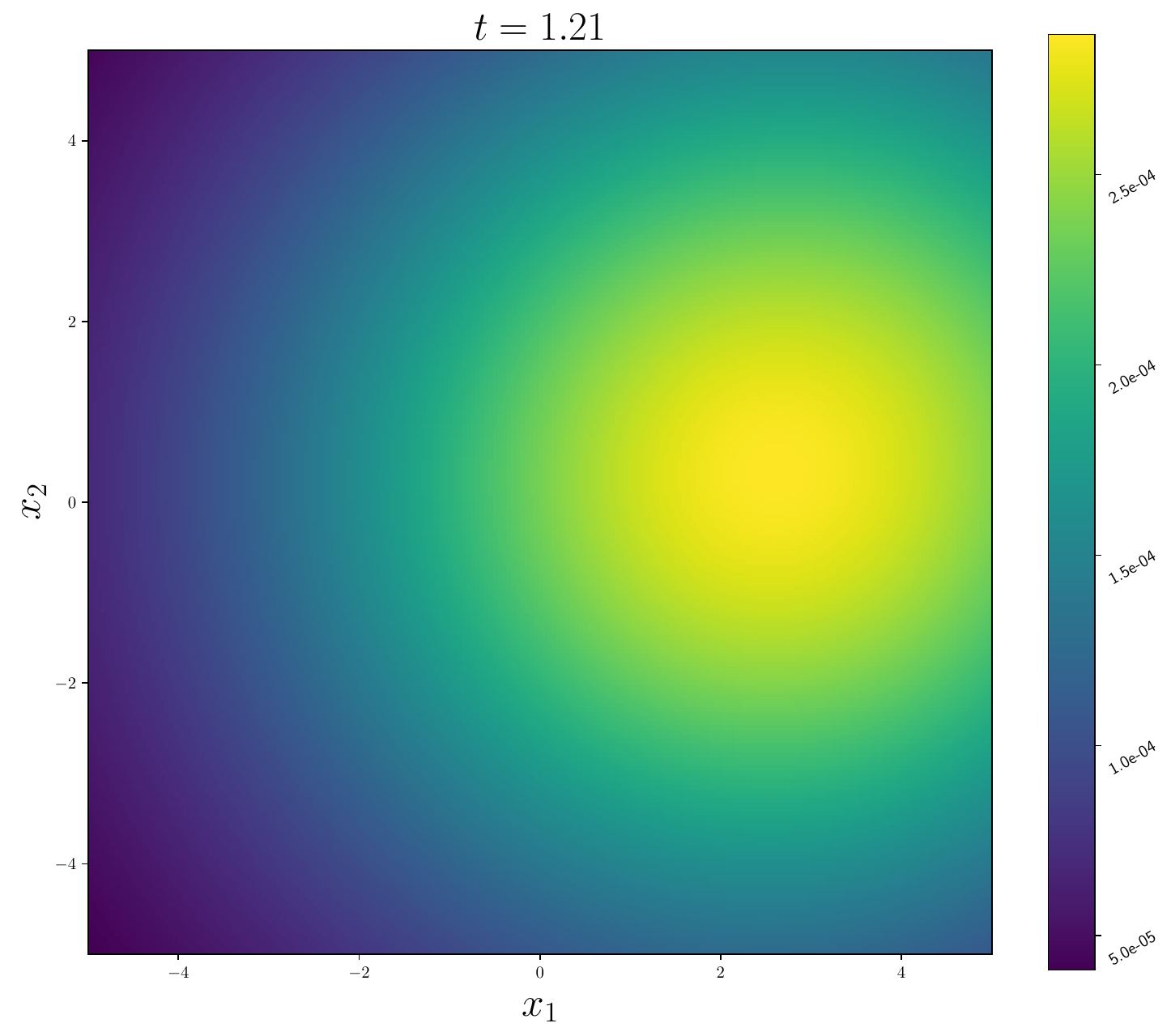}
    \hfill\\
    \caption{
        Numerical solution to a 3D FFPE with a sum of 10 Gaussians as the initial condition, projected onto the plane $x_{3} = 0$.
        Snapshots are shown at times $t = 0$, $0.01$, $0.06$, $0.11$, $0.31$, $0.61$, $0.91$, and $1.21$ (from left to right, top to bottom).
        Parameters: $\drift = (2, 1, -1)^{T}$, $\diffusionOrdinary = 1$, $\diffusionFractional = 4$, $\alpha = 1 / 2$.
        Red dots indicate the Gaussian centers, radii of red circles represent the standard deviation $\sigma$, and the numeric label near each center denotes its corresponding weight.
    }
    \label{fig:sum_of_Gaussians_test12}
\end{figure}
\Cref{fig:sum_of_Gaussians_test12} shows the numerical solution of the 3D FFPE with an initial condition given by the sum of 10 Gaussians.
For visualization, the centers of the initial Gaussians are placed on the plane $x_{3} = 0$, and the numerical solution is also projected onto this plane.
Although we only present a 3D example here, our FFPE solver is designed to efficiently handle high-dimensional problems, as illustrated in the previous subsections.
In high dimensions, our solver's computational cost scales primarily with the number of Gaussians rather than the spatial dimension, making it suitable for higher-dimensional applications with sums of Gaussians initial conditions.

%%%%%%%%%%%%%%%%%%%%%%%%%%%Section 5 %%%%%%%%%%%%%%%%%%%%%%%%%%

\section{Conclusion}\label{sec:conclusion}
We have developed an efficient and accurate numerical method for solving the free-space FFPE with constant coefficients and Dirac-delta initial conditions \Cref{eq:Dirac_general_FFPE_hd}. The method is concisely outlined in \Cref{rmk:Dirac_concise_explanation}.
The numerical method we proposed relies on \Cref{eq:Dirac_general_FFPE_Bessel}, which provides the single integral representation of the solution to the FFPE using Bessel functions of the first kind.
Fast evaluation of such integrals involves carefully handling near-origin integrals using expansion or re-weighting techniques and far-field integrals using windowing function techniques.
The algorithms also incorporate special techniques, such as the scaling law, to handle some extreme cases.
We have conducted numerical experiments with various parameter choices, including classical and fractional diffusion coefficients, fractional power, and dimensional settings. Additionally, we have performed empirical studies to assess the scalability of the method with respect to dimension.
The results shown in \Cref{sec:numerical_results} demonstrate the accuracy and efficiency of the proposed algorithms.

We note that the hyperparameters selected in our algorithms ($N, K_{1}, K_{2}, K_{3}, K_{4}$ in \Cref{alg:Dirac_integration_with_singularity}, $M, M_{\text{max}}, \epsilon$ in \Cref{alg:Dirac_integration_for_slow_decay}, $L$ in \Cref{alg:Dirac_basic_solver_for_general_FFPE_with_zero_displacement,alg:Dirac_basic_solver_for_general_FFPE} and $y_{\text{limit}}, \pi / 2$ in \Cref{alg:Dirac_solver_for_general_FFPE_in_hd}) are based on a balance of accuracy, efficiency, and generality.
These hyperparameters can be tuned to enhance accuracy or efficiency for specific cases, which is beyond the scope of this work. Investigating these details could be a future research.

Our solution provides the fundamental solution of the FFPE, which can potentially be used to solve other types of initial conditions (e.g., the sum-of-Gaussians case discussed in \Cref{sec:extension_to_sum_of_Gaussians_IC}) and approximate solutions to FFPE with variable parameters.
These applications are left for future work.
One important follow-up is to consider initial conditions that are characteristic functions of a set, which is useful for nonlocal threshold dynamics \cite{caffarelli2010convergence,jiang2018efficient}.
Another potential research direction is considering a variable $\alpha$.
Additionally, our numerical method accurately approximates integrals with oscillating and slowly decaying integrands, suggesting the possibility of extending our methodology to other similar integral forms.
Finally, the fractional Fokker-Planck solver can be employed to sample or recover probability distributions, which is highly useful in fields such as data science where data is driven by L\'evy processes.
We have explored this application in our recent work \cite{ye2025robust}.

%%%%%%%%%%%%%%%%%%%%%%%%Acknowledgement%%%%%%%%%%%%%%%%%%%%%%%%

\section*{Acknowledgement}
Q.~Ye and X.~Tian were supported in part by NSF DMS-2111608, NSF DMS-2240180 and the Alfred P. Sloan Fellowship.
D.~Wang was partially supported by National Natural Science Foundation of China (Grant No. 12101524), Guangdong Basic and Applied Basic Research Foundation (Grant No. 2023A1515012199), Shenzhen Science and Technology Innovation Program (Grant No. JCYJ20220530143803007, RCYX20221008092843046), Guangdong Provincial Key Laboratory of Mathematical Foundations for Artificial Intelligence (2023B1212010001), and Hetao Shenzhen-Hong Kong Science and Technology Innovation Cooperation Zone Project (No.HZQSWS-KCCYB-2024016).

We would like to thank the anonymous reviewers for their careful reading and constructive suggestions, which have helped improve the clarity and quality of this manuscript.

%%%%%%%%%%%%%%%%%%%%%%%%Code Availability%%%%%%%%%%%%%%%%%%%%%%%%

\section*{Code Availability}
Code is available at \texttt{https://github.com/ACMathX/FFPEDDIC}.

%%%%%%%%%%%%%%%%%%%%%%%%%%%Reference %%%%%%%%%%%%%%%%%%%%%%%%%%

\bibliographystyle{abbrv}
\bibliography{bib/ref}

%%%%%%%%%%%%%%%%%%%%%%%%%%% Appendix %%%%%%%%%%%%%%%%%%%%%%%%%%

\clearpage
\appendix

%%%%%%%%%%%%%%%%%%%%%%%%%%%Appendix A%%%%%%%%%%%%%%%%%%%%%%%%%%

\section{Proof of \texorpdfstring{\Cref{thm:farfield}}{Theorem 3.1}}\label{sec:proofoftheorem}
We prove \Cref{thm:farfield} in this section.
Notice that
\begin{equation*}
    \begin{split}
        &\int_{\gamma M}^{\infty} u(z) (1 - w_{M, \gamma}(z)) e^{i k z} \, d z\\
        = & \gamma M \int_{1}^{\infty} u(\gamma M z) (1 - w_{M, \gamma}(\gamma M z) ) e^{i k \gamma M z} \, d z\\
        = & \gamma M \int_{1}^{\infty} u(\gamma M z) (1 - w_{\frac{1}{\gamma}, \gamma}(z) ) e^{i k \gamma M z} \, d z
    \end{split}
\end{equation*}
Since $u(\gamma M z) ( 1 - w_{\frac{1}{\gamma}, \gamma}(z) )$ and all its derivatives vanish at $1$ and $\infty$, we further have
\begin{equation*}
    \begin{split}
        & \quad \, \gamma M \int_{1}^{\infty} u(\gamma M z) ( 1 - w_{\frac{1}{\gamma}, \gamma}(z) ) e^{i k \gamma M z} \, d z\\
        &= - \frac{1}{i k} \int_{1}^{\infty} \frac{d}{d z} \left ( u(\gamma M z) ( 1 - w_{\frac{1}{\gamma}, \gamma}(z) ) \right ) e^{i k \gamma M z} \, d z\\
        &= \frac{(- 1)^{n}}{(i k)^{n} (\gamma M)^{n - 1}} \int_{1}^{\infty} \frac{d^{n}}{d z^{n}} \left ( u(\gamma M z) ( 1 - w_{\frac{1}{\gamma}, \gamma}(z) ) \right ) e^{i k \gamma M z} \, d z\\
        &= \frac{(- 1)^{n}}{(i k)^{n} (\gamma M)^{n - 1}} \int_{1}^{1 / \gamma} \frac{d^{n}}{d z^{n}} \left ( u(\gamma M z) ( 1 - w_{\frac{1}{\gamma}, \gamma}(z) ) \right ) e^{i k \gamma M z} \, d z
    \end{split}
\end{equation*}
for any $n \in \mathbb{N}$.
Denote
\begin{equation*}
    \mathcal{B}(n)
    = \max_{x \in [1, 1 / \gamma]} \left | \frac{d^{n}}{d x^{n}}( 1 - w_{\frac{1}{\gamma}, \gamma}(x) ) \right |
    = \max_{x \in [1, 1 / \gamma]} \left | w^{(n)}_{\frac{1}{\gamma}, \gamma}(x) \right |.
\end{equation*}
Notice that $\frac{d^{n}}{d z^{n}} \left ( u(\gamma M z) \right ) = (\gamma M)^{n} u^{(n)}(\gamma M z)$, and for each $n\in \mathbb{N}$,
\begin{equation*}
    u^{(n)}(x)
    = \sum_{k = 0}^{n} C_{n}^{k}(l, q) (- \tau q)^{k} x^{l - n + k q} \exp( - \tau x^{q} ),
\end{equation*}
where $C_{n}^{k}(l, q)$ satisfies the estimate
\begin{equation*}
    \left | C_{n}^{k}(l, q) \right |
    \leq \binom{n}{k} (l + n - 1)_{n - k},
\end{equation*}
for $q \in (0, 2]$.
Therefore, we define
\begin{align*}
    \mathcal{A}(n, l)
    := \max_{k \in \{ 0, 1, \cdots, n \}} \binom{n}{k} (l + n - 1)_{n - k},
\end{align*}
and obtain
\begin{equation*}
    \begin{split}
        \max_{z \in [1, 1 / \gamma]} \left | \frac{d^{n}}{d z^{n}} \left ( u(\gamma M z) \right ) \right |
        &\leq (\gamma M)^{n} \mathcal{A}(n, l) \sum_{k = 0}^{n} (\tau q)^{k} M^{l - n + k q} \exp( - \gamma^{q} \tau M^{q} )\\
        &\leq \gamma^{n} M^{l} \mathcal{A}(n, l) \sum_{k = 0}^{n} \left ( \tau q M^{q} \right )^{k} \exp( - \gamma^{q} \tau M^{q} ).
    \end{split}
\end{equation*}
Notice that $\sum_{k = 0}^{n} \left ( \tau q M^{q} \right )^{k} \leq \min \{ \frac{1}{1 - \tau q M^{q}}, n+1 \}$ for $\tau q M^{q} \leq 1$ and $\sum_{k = 0}^{n} \left ( \tau q M^{q} \right )^{k} \leq (n + 1) (\tau q M^{q})^{n}$ for $\tau q M^{q} > 1$.
We therefore arrive at
\begin{equation*}
    \begin{aligned}
        &\quad \, \max_{z \in [1, 1 / \gamma]} \left | \frac{d^{n}}{d z^{n}} \left ( u( \gamma M z ) \right ) \right |\\
        &\leq
        \left \{
        \begin{aligned}
            &\gamma^{n} M^{l} \mathcal{A}(n, l) \min \left \{ \frac{1}{1 - \tau q M^{q}}, n + 1 \right \} \exp( - \gamma^{q} \tau M^{q} )
            &&\text{for } \tau q M^{q} \leq 1,\\
            &(n + 1) \gamma^{n} M^{l} \mathcal{A}(n, l) (\tau q M^{q})^{n} \exp( - \gamma^{q} \tau M^{q} )
            &&\text{for } \tau q M^{q} > 1.
        \end{aligned}
        \right .
    \end{aligned}
\end{equation*}
Combining the above estimates, we have
\begin{equation*}
    \begin{split}
        &\quad \, \left | \int_{\gamma M}^{\infty} u(z) (1 - w_{M, \gamma}(z)) e^{i k z} \, d z \right |\\
        &\leq
        \left \{
        \begin{aligned}
            &\frac{\mathcal{C}_{1}(n, l) M^{l}}{k^{n} (\gamma M)^{n - 1} \gamma} \min \left \{ \frac{1}{1 - \tau q M^{q}}, n + 1 \right \} \exp( - \gamma^{q} \tau M^{q})
            &&\text{for } \tau q M^{q} \leq 1,\\
            &
            \begin{aligned}
                &\frac{\mathcal{C}_{2}(n, l) M^{l} (\tau q M^{q})^{n}}{k^{n} (\gamma M)^{n - 1} \gamma} \exp( - \gamma^{q} \tau M^{q} )\\
                &\qquad = \mathcal{C}_{2}(n, l) M^{l + 1} \left ( \frac{\tau q M^{q - 1}}{k \gamma} \right )^{n} \exp( - \gamma^{q} \tau M^{q} )
            \end{aligned}
            &&\text{for } \tau q M^{q} > 1,
        \end{aligned}
        \right .
    \end{split}
\end{equation*}
where $\mathcal{C}_{1}(n, l)$ and $\mathcal{C}_{2}(n, l) $ are defined as $\mathcal{C}_{1}(n, l) = \sum_{j = 0}^{n} \binom{n}{j} \mathcal{A}(n - j, l) \mathcal{B}(j)$ and $\mathcal{C}_{2}(n, l) = (n + 1) \sum_{j = 0}^{n} \binom{n}{j} \mathcal{A}(n - j, l) \mathcal{B}(j) \left ( \tau q M^{q} \right )^{- j}$.
\qed

%%%%%%%%%%%%%%%%%%%%%%%%%%%Appendix B%%%%%%%%%%%%%%%%%%%%%%%%%%

\section{Explicit formulas of solutions for \texorpdfstring{$\alpha = 1 / 2$}{α = 1 / 2}}\label{appendix:special_cases_half_alpha}
For $\alpha = 1 / 2$, explicit formulas of solutions may be obtained using special functions. Indeed, it is well-known that with $\alpha = 1 / 2$ and in the absence of ordinary diffusion ($\diffusionOrdinary = 0$), the solution $p(\boldsymbol{x}, t) = \tilde{p}(y, t)$ to \Cref{eq:Dirac_general_FFPE_hd} is given by the multivariate Cauchy distribution \cite{samorodnitsky1994stable} 
\begin{equation}\label{eq:special_cases_half_alpha_Cauchy_distribution}
    \tilde{p}(y, t; 0)
    = \frac{\Gamma \left ( \frac{d + 1}{2} \right )}{\pi^{\frac{d + 1}{2}}} \frac{\diffusionFractional t}{\displaystyle \left [ \left ( \diffusionFractional t \right )^{2} + y^{2} \right ]^{\frac{d + 1}{2}}}.
\end{equation}
A simple derivation of the above formula can be found at \Cref{appendix:special_cases_pure_FFPE_rational_alpha_background}.

With ordinary diffusion present ($\diffusionOrdinary > 0$), we can derive the explicit formula for $d = 1$ as follows:
\begin{equation}\label{eq:Dirac_special_case_of_present_ordinary_diffusion_1d}
    \begin{aligned}
        \tilde{p}(y, t; \diffusionOrdinary)
        = \frac{1}{2 \sqrt{4 \pi \diffusionOrdinary t}} \left [ \textnormal{erfcx} \left ( \frac{\diffusionFractional t - y i}{2 \sqrt{\diffusionOrdinary t}} \right ) + \textnormal{erfcx} \left ( \frac{\diffusionFractional t + y i}{2 \sqrt{\diffusionOrdinary t}} \right ) \right ].
    \end{aligned}
\end{equation}
Here $i$ is the imaginary unit and \textnormal{erfcx} is the scaled complementary error function defined as
\begin{equation}\label{eq:erfcx}
    \textnormal{erfcx}(z)
    = \exp \left ( z^{2} \right ) \textnormal{erfc}(z)
    = \frac{2}{\sqrt{\pi}} \exp \left ( z^{2} \right ) \int_{z}^{\infty} \exp \left ( - t^{2} \right ) \, d t,
\end{equation}
for any complex number $z$. Notice that when $z$ is a complex number, the integral in the above definition is understood by expanding the integrand $\exp \left ( - t^{2} \right )$ into its Maclaurin series and integrating term by term. For $d \geq 2$, we can express the solution \Cref{eq:Dirac_explicit_integration_general_FFPE_hd_polar} by
\begin{equation}
    \tilde{p}(y, t; \diffusionOrdinary)
    = \frac{S_{d - 2}} {2 (2 \pi)^{d}} \int_{0}^{\pi} \sin^{d - 2}(\theta) \left( \integralExpP_{d - 1} (\cos(\theta) y) + \integralExpM_{d - 1}(\cos(\theta) y) \right) \, d \theta,
\end{equation}
where $\integralExpPM_{d}(z)$ ($d \geq 0$) can be found through the recurrence relation
\begin{equation}
\label{eq:special_cases_half_alpha_recurrence_relation_inner_integral}
    \integralExpPM_{d + 1}(z)
    = - \frac{\diffusionFractional t \pm z i}{2 \diffusionOrdinary t} \integralExpPM_{d}(z)
    + \frac{d}{2 \diffusionOrdinary t} \integralExpPM_{d - 1}(z),
    \qquad d \geq 1
\end{equation}
with $\integralExpPM_{0}$ and $\integralExpPM_{1}$ given by 
\begin{equation}
\label{eq:special_cases_half_alpha_recurrence_relation_inner_integral_initial_quantities}
    \left \{
    \begin{aligned}
        \integralExpPM_{0}(z)
        &= \frac{\sqrt{\pi}}{2 \sqrt{\diffusionOrdinary t}} \textnormal{erfcx} \left ( \frac{\diffusionFractional t \pm z i}{2 \sqrt{\diffusionOrdinary t}} \right ),\\
        \integralExpPM_{1}(z)
        &= \frac{1}{2 \diffusionOrdinary t} - \frac{\sqrt{\pi} \left ( \diffusionFractional t \pm z i \right )}{4 \left ( \diffusionOrdinary t \right )^{3 / 2}} \textnormal{erfcx} \left ( \frac{\diffusionFractional t \pm z i}{2 \sqrt{\diffusionOrdinary t}} \right ).
    \end{aligned}
    \right .
\end{equation}

Getting an explicit formulation for $\tilde{p}(y, t; \diffusionOrdinary)$ when $d \geq 2$ is difficult.
It turns out that when $d$ is an odd number, $\tilde{p}(y, t; \diffusionOrdinary)$ can be expressed by
\begin{equation}
    \begin{aligned}
        \tilde{p}(y, t; \diffusionOrdinary)
        = \frac{S_{d - 2}} {(2 \pi)^{d}} \integralSinExpP_{d - 2, d - 1},
    \end{aligned}
\end{equation}
where $\integralSinExpP_{p, q}$ can be found through the recurrence relation
\begin{equation}
\label{eq:special_cases_half_alpha_recurrence_relation_outer_integral}
    \begin{aligned}
        \integralSinExpP_{p, q}
        &= - \frac{(p - 1)(p - 3)}{y^{2}} \integralSinExpP_{p - 4, q - 2}
        + \frac{(p - 1)(p - 2)}{y^{2}} \integralSinExpP_{p - 2, q - 2},
        \qquad p \geq 4, q \geq 2
    \end{aligned}
\end{equation}
with $\integralSinExpP_{1, q}$ and $\integralSinExpP_{3, q}$ given by
\begin{equation}
\label{eq:special_cases_half_alpha_recurrence_relation_outer_integral_initial_quantities}
    \left \{
    \begin{aligned}
        \integralSinExpP_{1, q}
        &= \frac{1}{y i} \left ( \integralExpM_{q - 1}(y) - \integralExpP_{q - 1}(y) \right ),
        \quad q \geq 1,\\
        \integralSinExpP_{3, q}
        &= - \frac{2}{y^{2}} \left ( \integralExpM_{q - 2}(y) + \integralExpP_{q - 2}(y) \right )
        + \frac{2}{y^{3} i} \left ( \integralExpM_{q - 3}(y) - \integralExpP_{q - 3}(y) \right ),
        \quad q \geq 3.
    \end{aligned}
    \right .
\end{equation}

Further details and the derivation of these formulations are provided in \Cref{appendix:special_cases_half_alpha_derivation_part1,appendix:special_cases_half_alpha_derivation_part2,appendix:special_cases_half_alpha_derivation_part3,appendix:special_cases_half_alpha_derivation_part4}. 
In \Cref{sec:numerical_results}, these explicit formulas for the special case $\alpha = 1 / 2$ are used as a reference for verifying the accuracy of our numerical algorithm.

In \Cref{appendix:special_cases_pure_FFPE_rational_alpha}, we discuss about the formulation when $\diffusionOrdinary = 0$ and $\alpha \in (0, 1)$ being a rational number based on \Cref{eq:Dirac_general_FFPE_Bessel}.

\subsection{Deriving the recurrence relation (inner integral)}\label{appendix:special_cases_half_alpha_derivation_part1}
In this subsection, we derive
\Cref{eq:Dirac_special_case_of_present_ordinary_diffusion_1d,eq:special_cases_half_alpha_recurrence_relation_inner_integral}.
Denote
\begin{equation*}
    \integralExpPM_{d}(z)
    = \int_{0}^{\infty} r^{d} \exp( \mp r z i ) \exp \left ( - \left ( \diffusionOrdinary r^{2} + \diffusionFractional r \right ) t \right ) \, d r.
\end{equation*}
Then when $d \geq 1$ and $z \neq 0$, through integration by parts, we can derive
\begin{equation*}
    \begin{aligned}
        &\quad \, \integralExpPM_{d}(z)\\
        &= \left . \mp \frac{1}{z i} r^{d} \exp( \mp r z i ) \exp \left ( - \left ( \diffusionOrdinary r^{2} + \diffusionFractional r \right ) t \right ) \right |_{r = 0}^{\infty}\\
        &\pm \frac{1}{z i} \int_{0}^{\infty} \left ( d r^{d - 1} - 2 \diffusionOrdinary t r^{d + 1} - \diffusionFractional t r^{d} \right ) \exp( \mp r z i ) \exp \left ( - \left ( \diffusionOrdinary r^{2} + \diffusionFractional r \right ) t \right ) \, d r\\
        &= \pm \frac{d}{z i} \integralExpPM_{d - 1}(z) \mp \frac{2 \diffusionOrdinary t}{z i} \integralExpPM_{d + 1}(z) \mp \frac{\diffusionFractional t}{z i} \integralExpPM_{d}(z),
    \end{aligned}
\end{equation*}
hence \Cref{eq:special_cases_half_alpha_recurrence_relation_inner_integral} is shown.
With initial quantities $\integralExpPM_{0}(z)$ and $\integralExpPM_{1}(z)$ derived in \Cref{appendix:special_cases_half_alpha_derivation_part2}, we are able to obtain any $\integralExpPM_{d}(z)$ for $d \geq 0$.

Notice that in the one-dimensional case, we have directly:
\begin{equation*}
    \tilde{p}(y, t; \diffusionOrdinary)
    = \frac{1}{2 \pi} \left [ \integralExpP_{0}(y) + \integralExpM_{0}(y) \right ],
\end{equation*}
which leads to \Cref{eq:Dirac_special_case_of_present_ordinary_diffusion_1d}.

\subsection{Deriving the initial quantities (inner integral)}\label{appendix:special_cases_half_alpha_derivation_part2}
In this subsection, we derive \Cref{eq:special_cases_half_alpha_recurrence_relation_inner_integral_initial_quantities}.
As shown in \Cref{eq:erfcx}, the complementary error function is define as
\begin{equation*}
    \textnormal{erfc}(z)
    = \frac{2}{\sqrt{\pi}} \int_{z}^{\infty} \exp \left ( - t^{2} \right ) \, d t
    = 1 - \textnormal{erf}(z),
\end{equation*}
where the error function is defined as
\begin{equation*}
    \textnormal{erf}(z)
    = \frac{2}{\sqrt{\pi}} \int_{0}^{z} \exp \left ( - t^{2} \right ) \, d t.
\end{equation*}
A straightforward computation gives
\begin{equation*}
    \begin{aligned}
        &\quad \, \frac{\partial \textnormal{erf} \left ( \frac{2 \diffusionOrdinary t r + \diffusionFractional t \pm z i}{2 \sqrt{\diffusionOrdinary t}} \right )}{\partial r}
        = \frac{2}{\sqrt{\pi}} \frac{2 \diffusionOrdinary t}{2 \sqrt{\diffusionOrdinary t}} \exp \left ( - \left ( \frac{2 \diffusionOrdinary t r + \diffusionFractional t \pm z i}{2 \sqrt{\diffusionOrdinary t}} \right )^{2} \right )\\
        &= \frac{2 \sqrt{\diffusionOrdinary t}}{\sqrt{\pi}} \exp \left ( \mp r z i \right ) \exp \left ( - ( \diffusionOrdinary r^{2} + \diffusionFractional r ) t \right ) \exp \left ( - \left ( \frac{\diffusionFractional t \pm z i}{2 \sqrt{\diffusionOrdinary t}} \right )^{2} \right ),
    \end{aligned}
\end{equation*}
hence,
\begin{equation*}
    \begin{aligned}
        \integralExpPM_{0}(z)
        &= \int_{0}^{\infty} \exp( \mp r z i ) \exp \left ( - \left ( \diffusionOrdinary r^{2} + \diffusionFractional r \right ) t \right ) \, d r\\
        &= \frac{\sqrt{\pi}}{2 \sqrt{\diffusionOrdinary t}} \exp \left ( \left ( \frac{\diffusionFractional t \pm z i}{2 \sqrt{\diffusionOrdinary t}} \right )^{2} \right ) \left . \textnormal{erf} \left ( \frac{2 \diffusionOrdinary t r + \diffusionFractional t \pm z i}{2 \sqrt{\diffusionOrdinary t}} \right ) \right |_{r = 0}^{\infty}\\
        &= \frac{\sqrt{\pi}}{2 \sqrt{\diffusionOrdinary t}} \exp \left ( \left ( \frac{\diffusionFractional t \pm z i}{2 \sqrt{\diffusionOrdinary t}} \right )^{2} \right ) \textnormal{erfc} \left ( \frac{\diffusionFractional t \pm z i}{2 \sqrt{\diffusionOrdinary t}} \right )\\
        &= \frac{\sqrt{\pi}}{2 \sqrt{\diffusionOrdinary t}} \textnormal{erfcx} \left ( \frac{\diffusionFractional t \pm z i}{2 \sqrt{\diffusionOrdinary t}} \right ).
    \end{aligned}
\end{equation*}
Notice that
\begin{equation*}
    \begin{aligned}
        &\quad \, \frac{\partial \exp \left ( \mp r z i \right ) \exp \left ( - \left ( \diffusionOrdinary r^{2} + \diffusionFractional r \right ) t \right )}{\partial r}\\
        &= - \left ( 2 \diffusionOrdinary t r + \diffusionFractional t \pm z i \right ) \exp \left ( \mp r z i \right ) \exp \left ( - \left ( \diffusionOrdinary r^{2} + \diffusionFractional r \right ) t \right ),
    \end{aligned}
\end{equation*}
then we have
\begin{equation*}
    \begin{aligned}
        &\quad \, \integralExpPM_{1}(z)
        = \int_{0}^{\infty} r \exp \left ( \mp r z i \right ) \exp \left ( - \left ( \diffusionOrdinary r^{2} + \diffusionFractional r \right ) t \right ) \, d r\\
        &= \resizebox{0.85\linewidth}{!}{$\displaystyle
        - \frac{1}{2 \diffusionOrdinary t} \left . \left ( \frac{\sqrt{\pi}}{2 \sqrt{\diffusionOrdinary t}} \exp \left ( \left ( \frac{\diffusionFractional t \pm z i}{2 \sqrt{\diffusionOrdinary t}} \right )^{2} \right ) (\diffusionFractional t \pm z i) \textnormal{erf} \left ( \frac{2 \diffusionOrdinary t r + \diffusionFractional t \pm z i}{2 \sqrt{\diffusionOrdinary t}} \right ) + \exp \left ( \mp r z i \right ) \exp \left ( - \left ( \diffusionOrdinary r^{2} + \diffusionFractional r \right ) t \right ) \right ) \right |_{r = 0}^{\infty}
        $}\\
        &= - \frac{1}{2 \diffusionOrdinary t} \left ( \frac{\sqrt{\pi}}{2 \sqrt{\diffusionOrdinary t}} (\diffusionFractional t \pm z i) \textnormal{erfcx} \left ( \frac{\diffusionFractional t \pm z i}{2 \sqrt{\diffusionOrdinary t}} \right ) - 1 \right )\\
        &= \frac{1}{2 \diffusionOrdinary t} - \frac{\sqrt{\pi} (\diffusionFractional t \pm z i)}{4 \left ( \diffusionOrdinary t \right )^{3 / 2}} \textnormal{erfcx} \left ( \frac{\diffusionFractional t \pm z i}{2 \sqrt{\diffusionOrdinary t}} \right ).
    \end{aligned}
\end{equation*}

\subsection{Deriving the recurrence relation (outer integral)}\label{appendix:special_cases_half_alpha_derivation_part3}
In this subsection, we derive \Cref{eq:special_cases_half_alpha_recurrence_relation_outer_integral}.
Define
\begin{equation*}
    \integralSinExpPM_{p, q}
    = \int_{0}^{\pi} \sin^{p}(\theta) \integralExpPM_{q}(\cos(\theta) y) \, d \theta,
\end{equation*}
using the relation
\begin{equation*}
    \begin{aligned}
        &\quad \, \frac{\partial \integralExpPM_{q}(\cos(\theta) y)}{\partial \theta}\\
        &= \pm y i \sin(\theta) \int_{0}^{\infty} r^{q + 1} \exp \left ( \mp r \cos(\theta) y i \right ) \exp \left ( - \left ( \diffusionOrdinary r^{2} + \diffusionFractional r \right ) t \right ) \, d r\\
        &= \pm y i \sin(\theta) \integralExpPM_{q + 1}(\cos(\theta) y),
    \end{aligned}
\end{equation*}
we have the following recurrence relation achieved via integration by parts for $p \geq 4$ and $q \geq 2$:
\begin{equation*}
    \begin{aligned}
        &\quad \, \integralSinExpPM_{p, q}\\
        &= \left . \pm \frac{\sin^{p - 1}(\theta)}{y i} \integralExpPM_{q - 1}(\cos(\theta) y) \right |_{\theta = 0}^{\pi}\\
        &\quad \, \mp \frac{1}{y i} \int_{0}^{\pi} (p - 1) \sin^{p - 2}(\theta) \cos(\theta) \integralExpPM_{q - 1}(\cos(\theta) y) \, d \theta\\
        &= \mp \frac{p - 1}{y i} \int_{0}^{\pi} \sin^{p - 2}(\theta) \cos(\theta) \integralExpPM_{q - 1}(\cos(\theta) y) \, d \theta\\
        &= \mp \frac{p - 1}{y i}
        \Bigg [
            \left . \pm \frac{\sin^{p - 3}(\theta) \cos(\theta)}{y i} \integralExpPM_{q - 2}(\cos(\theta) y) \right |_{\theta = 0}^{\pi}\\
            &\qquad \quad \mp \frac{1}{y i} \int_{0}^{\pi} \left ( (p - 3) \sin^{p - 4}(\theta) \cos^{2}(\theta) - \sin^{p - 2}(\theta) \right ) \integralExpPM_{q - 2}(\cos(\theta) y) \, d \theta
        \Bigg ]\\
        &= - \frac{(p - 1)(p - 3)}{y^{2}} \integralSinExpPM_{p - 4, q - 2} + \frac{(p - 1) (p - 2)}{y^{2}} \integralSinExpPM_{p - 2, q - 2}.
    \end{aligned}
\end{equation*}
With initial quantities
$\integralSinExpPM_{1, q}$ and $\integralSinExpPM_{3, q}$
derived in \Cref{appendix:special_cases_half_alpha_derivation_part4}, we are able to obtain any $\integralSinExpPM_{p, q}$ for $q \geq p$ and odd $p \geq 1$.

As a result, we are able to find
\begin{equation*}
    \begin{aligned}
        \tilde{p}(y, t; \diffusionOrdinary)
        = \frac{S_{d - 2}} {2 (2 \pi)^{d}} \left [ \integralSinExpP_{d - 2, d - 1} + \integralSinExpM_{d - 2, d - 1} \right ]
        = \frac{S_{d - 2}} {(2 \pi)^{d}} \integralSinExpP_{d - 2, d - 1}
    \end{aligned}
\end{equation*}
through the recurrence relation we have derived for any odd dimension $d$. The last equality holds because both $\integralSinExpP_{p, q}$ and $\integralSinExpM_{p, q}$ have identical initial quantities and follow the same recurrence relation.

\subsection{Deriving the initial quantities (outer integral)}\label{appendix:special_cases_half_alpha_derivation_part4}
In this subsection, we derive \Cref{eq:special_cases_half_alpha_recurrence_relation_outer_integral_initial_quantities}. If $p = 1$ and $q \geq 1$, then we have
\begin{equation*}
    \begin{aligned}
        \integralSinExpPM_{1, q}
        = \left . \pm \frac{1}{y i} \integralExpPM_{q - 1}(\cos(\theta) y) \right |_{\theta = 0}^{\pi}
        = \frac{1}{y i} \left ( \integralExpM_{q - 1}(y) - \integralExpP_{q - 1}(y) \right ).
    \end{aligned}
\end{equation*}

Similarly, if $p = 3$ and $q \geq 3$, then we have
\begin{equation*}
    \begin{aligned}
        \integralSinExpPM_{3, q}
        &= \left . \pm \frac{\sin^{2}(\theta)}{y i} \integralExpPM_{q - 1}(\cos(\theta) y) \right |_{\theta = 0}^{\pi}
        \mp \frac{2}{y i} \int_{0}^{\pi} \sin(\theta) \cos(\theta) \integralExpPM_{q - 1}(\cos(\theta) y) \, d \theta\\
        &= \mp \frac{2}{y i} \int_{0}^{\pi} \sin(\theta) \cos(\theta) \integralExpPM_{q - 1}(\cos(\theta) y) \, d \theta\\
        &= \mp \frac{2}{y i} \left . \left ( \pm \frac{\cos(\theta)}{y i} \integralExpPM_{q - 2}(\cos(\theta) y) \right |_{\theta = 0}^{\pi}
        \pm \frac{1}{y i} \int_{0}^{\pi} \sin(\theta) \integralExpPM_{q - 2}(\cos(\theta) y) \, d \theta \right )\\
        &= - \frac{2}{y^{2}} \left ( \integralExpM_{q - 2}(y) + \integralExpP_{q - 2}(y) \right )
        + \frac{2}{y^{2}} \integralSinExpPM_{1, q - 2}\\
        &= - \frac{2}{y^{2}} \left ( \integralExpM_{q - 2}(y) + \integralExpP_{q - 2}(y) \right )
        + \frac{2}{y^{3} i} \left ( \integralExpM_{q - 3}(y) - \integralExpP_{q - 3}(y) \right ).
    \end{aligned}
\end{equation*}

\subsection{Examples for various dimensions}\label{appendix:special_cases_half_alpha_various_examples}
We can easily verify that \Cref{eq:special_cases_half_alpha_recurrence_relation_inner_integral} generates
\begin{equation*}
    \left \{
    \begin{aligned}
        \integralExpPM_{2}(z)
        &= - \frac{\diffusionFractional t \pm z i}{4 \left ( \diffusionOrdinary t \right )^{2}}
        + \frac{\sqrt{\pi} \left ( \left ( \diffusionFractional t \pm z i \right )^{2} + 2 \diffusionOrdinary t \right )}{8 \left ( \diffusionOrdinary t \right )^{5 / 2}} \textnormal{erfcx} \left ( \frac{\diffusionFractional t \pm z i}{2 \sqrt{\diffusionOrdinary t}} \right ),\\
        \integralExpPM_{3}(z)
        &= \frac{\left ( \diffusionFractional t \pm z i \right )^{2} + 4 \diffusionOrdinary t}{8 \left ( \diffusionOrdinary t \right )^{3}}\\
        &\quad \, - \frac{\sqrt{\pi} \left ( \left ( \diffusionFractional t \pm z i \right )^{3} + 6 \diffusionOrdinary t \left ( \diffusionFractional t \pm z i \right ) \right )}{16 \left ( \diffusionOrdinary t \right )^{7 / 2}} \textnormal{erfcx} \left ( \frac{\diffusionFractional t \pm z i}{2 \sqrt{\diffusionOrdinary t}} \right ).
    \end{aligned}
    \right .
\end{equation*}

In the three-dimensional case, we have
\begin{equation*}
   \begin{aligned}
        \tilde{p}(y, t; \diffusionOrdinary)
        = \frac{S_{1}} {(2 \pi)^{3}} \integralSinExpP_{1, 2}
        &= \frac{1}{16 \left ( \diffusionOrdinary t \pi \right )^{3 / 2} y i}
        \Bigg [
            - \left ( \diffusionFractional t - y i \right ) \textnormal{erfcx} \left ( \frac{\diffusionFractional t - y i}{2 \sqrt{\diffusionOrdinary t}} \right )\\
            &\qquad \qquad \qquad \qquad + \left ( \diffusionFractional t + y i \right ) \textnormal{erfcx} \left ( \frac{\diffusionFractional t + y i}{2 \sqrt{\diffusionOrdinary t}} \right )
        \Bigg ].
   \end{aligned}
\end{equation*}

In the five-dimensional case, we have
\begin{equation*}
    \begin{aligned}
        &\quad \, \tilde{p}(y, t; \diffusionOrdinary)
        = \frac{S_{3}}{(2 \pi)^{5}} \integralSinExpP_{3, 4}\\
        &= \frac{\diffusionFractional t}{16 \pi \left ( \diffusionOrdinary t \pi \right )^{2} y^{2}}
        - \resizebox{0.6\linewidth}{!}{$\displaystyle
        \frac{\diffusionFractional t}{32 \pi \left ( \diffusionOrdinary t \pi \right )^{3 / 2} y^{3} i}
        \left [
            \textnormal{erfcx} \left ( \frac{\diffusionFractional t - y i}{2 \sqrt{\diffusionOrdinary t}} \right )
            - \textnormal{erfcx} \left ( \frac{\diffusionFractional t + y i}{2 \sqrt{\diffusionOrdinary t}} \right )
        \right ]
        $}\\
        &\quad \, - \resizebox{0.8\linewidth}{!}{$\displaystyle
        \frac{1}{64 \left ( \diffusionOrdinary t \pi \right )^{5 / 2} y^{2}}
        \left [
            \left ( \diffusionFractional t - y i \right )^{2} \textnormal{erfcx} \left ( \frac{\diffusionFractional t - y i}{2 \sqrt{\diffusionOrdinary t}} \right )
            + \left ( \diffusionFractional t + y i \right )^{2} \textnormal{erfcx} \left ( \frac{\diffusionFractional t + y i}{2 \sqrt{\diffusionOrdinary t}} \right )
        \right ]
        $}.
    \end{aligned}
\end{equation*}

%%%%%%%%%%%%%%%%%%%%%%%%%%%Appendix C%%%%%%%%%%%%%%%%%%%%%%%%%%

\section{Explicit formulas of solutions for the pure FFPE with rational \texorpdfstring{$\alpha$}{α}}\label{appendix:special_cases_pure_FFPE_rational_alpha}
Explicit formulas of solutions for the pure FFPE, i.e., $\diffusionOrdinary = 0$,  with rational $\alpha$ can be obtained by using special functions. 
Let $S_{d - 1}$ denote the surface area of the unit $(d - 1)$-sphere embedded in $d$-dimensional Euclidean space and ${}_{p} F_{q}(a_{1}, \ldots, a_{p}; b_{1}, \ldots, b_{q}; z)$ is the generalized hypergeometric function defined as
\begin{equation*}
    {}_{p} F_{q}(a_{1}, \ldots, a_{p}; b_{1}, \ldots, b_{q}; z)
    = \sum_{n = 0}^{\infty} \frac{(a_{1})_{n} \cdots (a_{p})_{n}}{(b_{1})_{n} \cdots (b_{q})_{n}} \frac{z^{n}}{n!},
\end{equation*}
where $(a)_{n}$ is the rising Pochhammer symbol defined as $(a)_{0} := 1$ and $(a)_{n} := a (a + 1)_{n - 1}$ for $n \geq 1$.

Let $\alpha = p / q \in (0, 1)$, where $p, q \in \mathbb{N}_{+}$ and $\gcd(p, q) = 1$, the solution $p(\boldsymbol{x}, t) = \tilde{p}(y, t)$ to \Cref{eq:Dirac_general_FFPE_hd} for the case $\diffusionOrdinary = 0$ is given as
\begin{equation}\label{eq:special_cases_pure_FFPE}
    \begin{aligned}
        &\tilde{p}(y, t; 0)
        = \frac{(\diffusionFractional t)^{\frac{- d}{2 p / q}}}{(2 \pi)^{d}}
        \sum_{i = 0}^{p - 1} \frac{S_{d + 2 i - 1} \Gamma \left ( \frac{d + 2 i}{2 p / q} + 1 \right) z^{i}}{i! (d + 2 i) \pi^{i}}\\
        &\quad \ {}_{q - 1} F_{2 p - 2}
        \resizebox{0.7\linewidth}{!}{$\displaystyle
        \left (
            \left \{ \frac{d + 2 i}{2 p} + \frac{j}{q} \right \}_{j = 1}^{q - 1};
            \left \{ \frac{i + 1 + j}{p} \right \}_{\substack{j = 0\\i + 1 + j \neq p}}^{p - 1}, \left \{ \frac{d + 2 i + 2 j}{2 p} \right \}_{j = 1}^{p - 1};
            \frac{z^{p} q^{q}}{p^{2 p}}
        \right )
        $},
    \end{aligned}
\end{equation}
where $\displaystyle z = \frac{- (y / 2)^{2}}{(\diffusionFractional t)^{q / p}}$.

The formula is simple when $p = 1$. For example, when $\alpha = 1 / 3$, the solution is given by
\begin{equation*}
    \begin{aligned}
        \tilde{p}(y, t; 0)
        &= \frac{S_{d - 1} \Gamma \left ( \frac{3 d}{2} + 1 \right )}{(2 \pi)^{d} (\diffusionFractional t)^{\frac{3 d}{2}} d}
        \ {}_{2} F_{0} \left ( \frac{d}{2} + \frac{1}{3}, \frac{d}{2} + \frac{2}{3}; ; - \frac{3^{3} (y / 2)^{2}}{(\diffusionFractional t)^{3}} \right ).
    \end{aligned}
\end{equation*}

The derivation in this section is based on the solution representation in \Cref{eq:Dirac_general_FFPE_Bessel}.

\subsection{The classical case for \texorpdfstring{$\alpha = 1 / 2$}{α = 1 / 2}}\label{appendix:special_cases_pure_FFPE_rational_alpha_background}
The integral
\begin{equation*}
    \int_{0}^{\infty} \exp(- \tau r) J_{\nu}(y r) r^{\mu - 1} \, d r
\end{equation*}
has been studied for over a hundred years. It is well known \cite{watson1922treatise} that for any $\mu$ and $\nu$ such that $\mu + \nu$ has a positive real part, we have
\begin{equation}\label{eq:rational_order_special_case_integral}
    \begin{aligned}
        &\quad \, \int_{0}^{\infty} \exp(- \tau r) J_{\nu}(y r) r^{\mu - 1} \, d r\\
        &= \int_{0}^{\infty} \exp(- \tau r) r^{\mu - 1} \sum_{m = 0}^{\infty} \frac{(- 1)^{m}}{m! \Gamma(\nu + m + 1)} \left ( \frac{y r}{2} \right )^{\nu + 2 m} \, d r\\
        &= \sum_{m = 0}^{\infty} \frac{(- 1)^{m} (y / 2)^{\nu + 2 m}}{m! \Gamma(\nu + m + 1)} \int_{0}^{\infty} r^{\mu + \nu + 2 m - 1} \exp(- \tau r) \, d r\\
        &= \sum_{m = 0}^{\infty} \frac{(- 1)^{m} (y / 2)^{\nu + 2 m}}{m! \Gamma(\nu + m + 1)} \frac{\Gamma(\mu + \nu + 2 m)}{\tau^{\mu + \nu + 2 m}}\\
        &= \frac{(y / 2)^{\nu} \Gamma(\mu + \nu)}{\tau^{\mu + \nu} \Gamma(\nu + 1)} \sum_{m = 0}^{\infty} \frac{\left ( \mu + \nu \right )_{2 m} / 2^{2 m}}{\left ( \nu + 1 \right )_{m} m!} \left ( - \frac{y^{2}}{\tau^{2}} \right )^{m}\\
        &= \frac{(y / 2)^{\nu} \Gamma(\mu + \nu)}{\tau^{\mu + \nu} \Gamma(\nu + 1)} \sum_{m = 0}^{\infty} \frac{\left ( \frac{\mu + \nu}{2} \right )_{m} \left ( \frac{\mu + \nu + 1}{2} \right )_{m}}{\left ( \nu + 1 \right )_{m} m!} \left ( - \frac{y^{2}}{\tau^{2}} \right )^{m}\\
        &= \frac{(y / 2)^{\nu} \Gamma(\mu + \nu)}{\tau^{\mu + \nu} \Gamma(\nu + 1)} \ {}_{2} F_{1} \left ( \frac{\mu + \nu}{2}, \frac{\mu + \nu + 1}{2}; \nu + 1;  - \frac{y^{2}}{\tau^{2}} \right ).
    \end{aligned}
\end{equation}
By further applying the Euler transformation, we can express it in the following forms:
\begin{equation*}
    \begin{aligned}
        &\quad \, \int_{0}^{\infty} \exp(- \tau r) J_{\nu}(y r) r^{\mu - 1} \, d r\\
        &= \frac{(y / 2)^{\nu} \Gamma(\mu + \nu)}{\tau^{\mu + \nu} \Gamma(\nu + 1)} \left ( 1 + \frac{y^{2}}{\tau^{2}} \right )^{\frac{1}{2} - \mu} \ {}_{2} F_{1} \left ( \frac{\nu - \mu + 1}{2}, \frac{\nu - \mu}{2} + 1; \nu + 1;  - \frac{y^{2}}{\tau^{2}} \right )\\
    \end{aligned}
\end{equation*}
By taking $\mu = \nu + 2$, we obtain
\begin{equation}
    \int_{0}^{\infty} \exp(- \tau r) J_{\nu}(y r) r^{\nu + 1} \, d r
    = \frac{(2 \tau) (2 y)^{\nu} \Gamma \left ( \nu + \frac{3}{2} \right )}{\left [ \tau^{2} + y^{2} \right ]^{\nu + \frac{3}{2}} \sqrt{\pi}},
\end{equation}
if we further set $\nu = (d - 2) / 2$ and $\tau = \diffusionFractional t$, we obtain the multivariate Cauchy distribution for the special case where $\diffusionOrdinary = 0$ and $\alpha = 1 / 2$ (as in \Cref{eq:special_cases_half_alpha_Cauchy_distribution}):
\begin{equation*}
    \begin{aligned}
        &\quad \, \frac{1}{y^{(d - 2) / 2} (2 \pi)^{d / 2}} \int_{0}^{\infty} r^{d / 2} J_{(d - 2) / 2}(y r) \exp \left ( - \diffusionFractional t r \right ) \, d r\\
        &= \frac{1}{y^{(d - 2) / 2} (2 \pi)^{d / 2}} \frac{(2 \diffusionFractional t) (2 y)^{(d - 2) / 2} \Gamma \left ( \frac{d + 1}{2} \right )}{\left [ (\diffusionFractional t)^{2} + y^{2} \right ]^{\frac{d + 1}{2}} \sqrt{\pi}}\\
        &= \frac{\Gamma \left ( \frac{d + 1}{2} \right )}{\pi^{\frac{d + 1}{2}}} \frac{\diffusionFractional t}{\left [ (\diffusionFractional t)^{2} + y^{2} \right ]^{\frac{d + 1}{2}}}.
    \end{aligned}
\end{equation*}

\subsection{General cases for rational \texorpdfstring{$\alpha$}{α}}
Modify \Cref{eq:rational_order_special_case_integral} to treat the general $\alpha$, we have
\begin{equation*}
    \begin{aligned}
        &\quad \, \int_{0}^{\infty} \exp(- \tau r^{2 \alpha}) J_{\nu}(y r) r^{\mu - 1} \, d r\\
        &= \int_{0}^{\infty} \exp(- \tau r^{2 \alpha}) r^{\mu - 1} \sum_{m = 0}^{\infty} \frac{(- 1)^{m}}{m! \Gamma(\nu + m + 1)} \left ( \frac{y r}{2} \right )^{\nu + 2 m} \, d r\\
        &= \sum_{m = 0}^{\infty} \frac{(- 1)^{m} (y / 2)^{\nu + 2 m}}{m! \Gamma(\nu + m + 1)} \int_{0}^{\infty} r^{\mu + \nu + 2 m - 1} \exp(- \tau r^{2 \alpha}) \, d r\\
        &= \sum_{m = 0}^{\infty} \frac{(- 1)^{m} (y / 2)^{\nu + 2 m}}{m! \Gamma(\nu + m + 1)} \frac{\Gamma \left ( \frac{\mu + \nu + 2 m}{2 \alpha} \right)}{\tau^{\frac{\mu + \nu + 2 m}{2 \alpha}} (2 \alpha)}\\
        &= \frac{1}{2 \alpha} \sum_{m = 0}^{\infty} \frac{(- 1)^{m} (y / 2)^{\nu + 2 m} \Gamma \left ( \frac{\mu + \nu + 2 m}{2 \alpha} \right)}{m! \Gamma(\nu + m + 1) \tau^{\frac{\mu + \nu + 2 m}{2 \alpha}}}.
    \end{aligned}
\end{equation*}

Let $\alpha = p / q \in (0, 1)$, where $p, q \in \mathbb{N}_{+}$ and $\gcd(p, q) = 1$, i.e., $\alpha$ is a positive rational number.
We need to separate the summation into $p$ hypergeometric functions:
\begin{equation*}
    \begin{aligned}
        &\quad \, \int_{0}^{\infty} \exp(- \tau r^{2 p / q}) J_{\nu}(y r) r^{\mu - 1} \, d r\\
        &= \frac{1}{2 p / q}
        \sum_{m = 0}^{\infty} \frac{(- 1)^{m} (y / 2)^{\nu + 2 m} \Gamma \left ( \frac{\mu + \nu + 2 m}{2 p / q} \right)}{m! \Gamma(\nu + m + 1) \tau^{\frac{\mu + \nu + 2 m}{2 p / q}}}\\
        &= \frac{q}{2 p} \sum_{i = 0}^{p - 1}
        \sum_{m = 0}^{\infty} \frac{(- 1)^{p m + i} (y / 2)^{\nu + 2 p m + 2 i} \Gamma \left ( \frac{\mu + \nu + 2 i}{2 p / q} + q m \right)}{\Gamma(p m + i + 1) \Gamma(\nu + p m + i + 1) \tau^{\frac{\mu + \nu + 2 p m + 2 i}{2 p / q}}}\\
        &= \frac{q}{2 p} \sum_{i = 0}^{p - 1} \frac{(- 1)^{i} (y / 2)^{\nu + 2 i} \Gamma \left ( \frac{\mu + \nu + 2 i}{2 p / q} \right)}{\Gamma(i + 1) \Gamma(\nu + i + 1) \tau^{\frac{\mu + \nu + 2 i}{2 p / q}}}\\
        &\qquad \qquad \sum_{m = 0}^{\infty} \frac{\left ( \frac{\mu + \nu + 2 i}{2 p / q} \right)_{q m} (1)_{m}}{(i + 1)_{p m} (\nu + i + 1)_{p m} m!} \left ( (- 1)^{p} \frac{(y / 2)^{2 p}}{\tau^{q}} \right )^{m}\\
        &= \frac{q}{2 p} \sum_{i = 0}^{p - 1} \frac{(- 1)^{i} (y / 2)^{\nu + 2 i} \Gamma \left ( \frac{\mu + \nu + 2 i}{2 p / q} \right)}{i! \Gamma(\nu + i + 1) \tau^{\frac{\mu + \nu + 2 i}{2 p / q}}}\\
        &\qquad \qquad \sum_{m = 0}^{\infty} \frac{(1)_{m} \prod_{j = 0}^{q - 1} \left ( \frac{\mu + \nu + 2 i}{2 p} + \frac{j}{q} \right)_{m}}{\prod_{j = 0}^{p - 1} \left ( \frac{i + 1 + j}{p} \right )_{m} \prod_{j = 0}^{p - 1} \left ( \frac{\nu + i + 1 + j}{p} \right )_{m} m!} \left ( (- 1)^{p} \frac{q^{q} y^{2 p}}{p^{p} p^{p} 2^{2 p} \tau^{q}} \right )^{m}\\
        &= \frac{q}{2 p} \sum_{i = 0}^{p - 1} \frac{(- 1)^{i} (y / 2)^{\nu + 2 i} \Gamma \left ( \frac{\mu + \nu + 2 i}{2 p / q} \right)}{i! \Gamma(\nu + i + 1) \tau^{\frac{\mu + \nu + 2 i}{2 p / q}}}\\
        &\qquad \ {}_{q} F_{2 p - 1}
        \resizebox{0.8\linewidth}{!}{$\displaystyle
        \left (
            \left \{ \frac{\mu + \nu + 2 i}{2 p} + \frac{j}{q} \right \}_{j = 0}^{q - 1};
            \left \{ \frac{i + 1 + j}{p} \right \}_{\substack{j = 0\\i + 1 + j \neq p}}^{p - 1}, \left \{ \frac{\nu + i + 1 + j}{p} \right \}_{j = 0}^{p - 1};
            \frac{(- 1)^{p} q^{q} y^{2 p}}{(2 p)^{2 p} \tau^{q}}
        \right )
        $}.
    \end{aligned}
\end{equation*}
Notice that $\displaystyle 1 \in \left \{ \frac{i + 1 + j}{p} \right \}_{j = 0}^{p - 1}$ for any $0 \leq i \leq p - 1$ and $(1)_{m} = m!$.

Moreover, in our scenario, we have $\nu = (d - 2) / 2$ and $\mu = (d + 2) / 2$. Substituting $\mu = \nu + 2$ gives
\begin{equation*}
    \begin{aligned}
        &\quad \, \int_{0}^{\infty} \exp(- \tau r^{2 p / q}) J_{\nu}(y r) r^{\nu + 1} \, d r\\
        &= \frac{q}{2 p} \sum_{i = 0}^{p - 1} \frac{(- 1)^{i} (y / 2)^{\nu + 2 i} \Gamma \left ( \frac{\nu + i + 1}{p / q} \right)}{i! \Gamma(\nu + i + 1) \tau^{\frac{\nu + i + 1}{p / q}}}\\
        &\qquad \ {}_{q - 1} F_{2 p - 2}
        \resizebox{0.8\linewidth}{!}{$\displaystyle
        \left (
            \left \{ \frac{\nu + i + 1}{p} + \frac{j}{q} \right \}_{j = 1}^{q - 1};
            \left \{ \frac{i + 1 + j}{p} \right \}_{\substack{j = 0\\i + 1 + j \neq p}}^{p - 1}, \left \{ \frac{\nu + i + 1 + j}{p} \right \}_{j = 1}^{p - 1};
            \frac{(- 1)^{p} q^{q} y^{2 p}}{(2 p)^{2 p} \tau^{q}}
        \right )
        $}.
    \end{aligned}
\end{equation*}
Further substitute $\nu = (d - 2) / 2$ and $\tau = \diffusionFractional t$, we have
\begin{equation*}
    \begin{aligned}
        &\quad \, \tilde{p}(y, t; 0)
        = \frac{1}{y^{(d - 2) / 2}} \frac{1}{(2 \pi)^{d / 2}} \int_{0}^{\infty} \exp(- \diffusionFractional t r^{2 p / q}) J_{(d - 2) / 2}(y r) r^{d / 2} \, d r\\
        &= \frac{q}{p} \frac{1}{2^{d} \pi^{d / 2} (\diffusionFractional t)^{\frac{d}{2 p / q}}}
        \sum_{i = 0}^{p - 1} \frac{(- 1)^{i} (y / 2)^{2 i} \Gamma \left ( \frac{d + 2 i}{2 p / q} \right)}{i! \Gamma \left ( \frac{d + 2 i}{2} \right ) (\diffusionFractional t)^{\frac{2 i}{2 p / q}}}\\
        &\qquad \ {}_{q - 1} F_{2 p - 2}
        \resizebox{0.7\linewidth}{!}{$\displaystyle
        \left (
            \left \{ \frac{d + 2 i}{2 p} + \frac{j}{q} \right \}_{j = 1}^{q - 1};
            \left \{ \frac{i + 1 + j}{p} \right \}_{\substack{j = 0\\i + 1 + j \neq p}}^{p - 1}, \left \{ \frac{d + 2 i + 2 j}{2 p} \right \}_{j = 1}^{p - 1};
            \frac{(- 1)^{p} q^{q} y^{2 p}}{(2 p)^{2 p} (\diffusionFractional t)^{q}}
        \right )
        $}\\
        &= \frac{(\diffusionFractional t)^{\frac{- d}{2 p / q}}}{(2 \pi)^{d}}
        \sum_{i = 0}^{p - 1} \frac{S_{d + 2 i - 1} \Gamma \left ( \frac{d + 2 i}{2 p / q} + 1 \right) z^{i}}{i! (d + 2 i) \pi^{i}}\\
        &\qquad \ {}_{q - 1} F_{2 p - 2}
        \resizebox{0.7\linewidth}{!}{$\displaystyle
        \left (
            \left \{ \frac{d + 2 i}{2 p} + \frac{j}{q} \right \}_{j = 1}^{q - 1};
            \left \{ \frac{i + 1 + j}{p} \right \}_{\substack{j = 0\\i + 1 + j \neq p}}^{p - 1}, \left \{ \frac{d + 2 i + 2 j}{2 p} \right \}_{j = 1}^{p - 1};
            \frac{z^{p} q^{q}}{p^{2 p}}
        \right )
        $},
    \end{aligned}
\end{equation*}
where $\displaystyle z = \frac{- (y / 2)^{2}}{(\diffusionFractional t)^{q / p}}$.

\subsection{Examples for various parameters}
\Cref{eq:special_cases_pure_FFPE} is simple when $p = 1$. When $\alpha = 1 / 3$, the solution is given by
\begin{equation*}
    \begin{aligned}
        \tilde{p}(y, t; 0)
        &= \frac{S_{d - 1} \Gamma \left ( \frac{3 d}{2} + 1 \right )}{(2 \pi)^{d} (\diffusionFractional t)^{\frac{3 d}{2}} d}
        \ {}_{2} F_{0} \left ( \frac{d}{2} + \frac{1}{3}, \frac{d}{2} + \frac{2}{3}; ; - \frac{3^{3} (y / 2)^{2}}{(\diffusionFractional t)^{3}} \right ).
    \end{aligned}
\end{equation*}
For the case $\alpha = 1 / 4$, we can write the solution as
\begin{equation*}
    \begin{aligned}
        \tilde{p}(y, t; 0)
        &= \frac{S_{d - 1} \Gamma \left ( 2 d + 1 \right )}{(2 \pi)^{d} (\diffusionFractional t)^{2 d} d}
        \ {}_{3} F_{0} \left ( \frac{d}{2} + \frac{1}{4}, \frac{d}{2} + \frac{2}{4}, \frac{d}{2} + \frac{3}{4}; ; - \frac{4^{4} (y / 2)^{2}}{(\diffusionFractional t)^{4}} \right ).
    \end{aligned}
\end{equation*}

For the case $\alpha = 2 / 3$, the solution is given by
\begin{equation}\label{eq:special_cases_alpha_2_divided_by_3}
    \begin{aligned}
        &\quad \, \tilde{p}(y, t; 0)\\
        &= \frac{1}{(2 \pi)^{d} (\diffusionFractional t)^{\frac{3 d}{4}}}
        \Bigg [
            \frac{S_{d - 1} \Gamma \left ( \frac{3 d}{4} + 1 \right )}{d}
            \ {}_{2} F_{2} \left ( \frac{d}{4} + \frac{1}{3}, \frac{d}{4} + \frac{2}{3}; \frac{1}{2}, \frac{d + 2}{4}; \frac{3^{3} z^{2}}{2^{4}} \right )\\
            &\quad + \frac{S_{d + 1} \Gamma \left ( \frac{3 (d + 2)}{4} + 1 \right ) z}{(d + 2) \pi}
            \ {}_{2} F_{2} \left ( \frac{d + 2}{4} + \frac{1}{3}, \frac{d + 2}{4} + \frac{2}{3}; \frac{3}{2}, \frac{d + 4}{4}; \frac{3^{3} z^{2}}{2^{4}} \right )
        \Bigg ],
    \end{aligned}
\end{equation}
where $\displaystyle z = \frac{- (y / 2)^{2}}{(\diffusionFractional t)^{3 / 2}}$.

\begin{remark}
    While theoretically sound, \Cref{eq:special_cases_alpha_2_divided_by_3} may encounter significant computational difficulties in practice, especially when $t > 0$ is small and $y$ is large.
    Specifically, the hypergeometric functions ${}_{2} F_{2}$ (in general, ${}_{q - 1} F_{2 p - 2}$ for even $p$) involved in the expressions can result in terms where very large quantities must be subtracted from each other.
    Such operations are highly susceptible to numerical precision issues, potentially leading to substantial inaccuracies in the computed values.
    As a result, the direct computation of this solution may be infeasible without implementing specialized numerical techniques to mitigate these precision errors. Therefore, we avoid using the formula to compute the solution.
\end{remark}

For the case $\alpha = 3 / 4$ (which corresponds to the Holtsmark distribution), we can write the solution as
\begin{equation*}
    \begin{aligned}
        &\quad \, \tilde{p}(y, t; 0)\\
        &= \frac{1}{(2 \pi)^{d} (\diffusionFractional t)^{\frac{2 d}{3}}}
        \Bigg [
            \frac{S_{d - 1} \Gamma \left ( \frac{2 d}{3} + 1 \right )}{d}
            \ {}_{3} F_{4}
            \resizebox{0.52\linewidth}{!}{$\displaystyle
            \left (
                \frac{d}{6} + \frac{1}{4}, \frac{d}{6} + \frac{2}{4}, \frac{d}{6} + \frac{3}{4}; \frac{1}{3}, \frac{2}{3}, \frac{d + 2}{6}, \frac{d + 4}{6}; \frac{4^{4} z^{3}}{3^{6}}
            \right )
            $}\\
            &\quad + \frac{S_{d + 1} \Gamma \left ( \frac{2 (d + 2)}{3} + 1 \right ) z}{(d + 2) \pi}
            \ {}_{3} F_{4}
            \resizebox{0.6\linewidth}{!}{$\displaystyle
            \left (
                \frac{d + 2}{6} + \frac{1}{4}, \frac{d + 2}{6} + \frac{2}{4}, \frac{d + 2}{6} + \frac{3}{4}; \frac{2}{3}, \frac{4}{3}, \frac{d + 4}{6}, \frac{d + 6}{6}; \frac{4^{4} z^{3}}{3^{6}}
            \right )
            $}\\
            &\quad + \frac{S_{d + 3} \Gamma \left ( \frac{2 (d + 4)}{3} + 1 \right ) z^{2}}{2 (d + 4) \pi^{2}}
            \ {}_{3} F_{4}
            \resizebox{0.6\linewidth}{!}{$\displaystyle
            \left (
                \frac{d + 4}{6} + \frac{1}{4}, \frac{d + 4}{6} + \frac{2}{4}, \frac{d + 4}{6} + \frac{3}{4}; \frac{4}{3}, \frac{5}{3}, \frac{d + 6}{6}, \frac{d + 8}{6}; \frac{4^{4} z^{3}}{3^{6}}
            \right )
            $}
        \Bigg ],
    \end{aligned}
\end{equation*}
where $\displaystyle z = \frac{- (y / 2)^{2}}{(\diffusionFractional t)^{4 / 3}}$.

%%%%%%%%%%%%%%%%%%%%%%%%%%%Appendix D%%%%%%%%%%%%%%%%%%%%%%%%%%

\section{Window Function Comparison}\label{appendix:window_function_comparison}
We extend the results we shown in \Cref{tab:Dirac_windowing_function} in this section to better understand the performance of various window functions.
We consider the following model problem same as the one we used in \Cref{tab:Dirac_windowing_function}:
\begin{equation*}
    \begin{aligned}
        I
        &= \int_{1}^{\infty} \cos(z) \exp(- 0.01 z^{0.8}) \, d z
        \approx - 0.82835952299051113433,\\
        I_{\textnormal{WF}}(M)
        &= \int_{1}^{M} \cos(z) \exp(- 0.01 z^{0.8}) w_{M}^{WF}(z) \, d z,\\
        E_{\textnormal{WF}}(M)
        &= |I_{\textnormal{WF}}(M) - I|,
    \end{aligned}
\end{equation*}
where $\textnormal{WF}$ denotes the choice of window function applied to the integrand.

\begin{figure}
    \centering
    \hfill
    \includegraphics[width=0.485\linewidth]{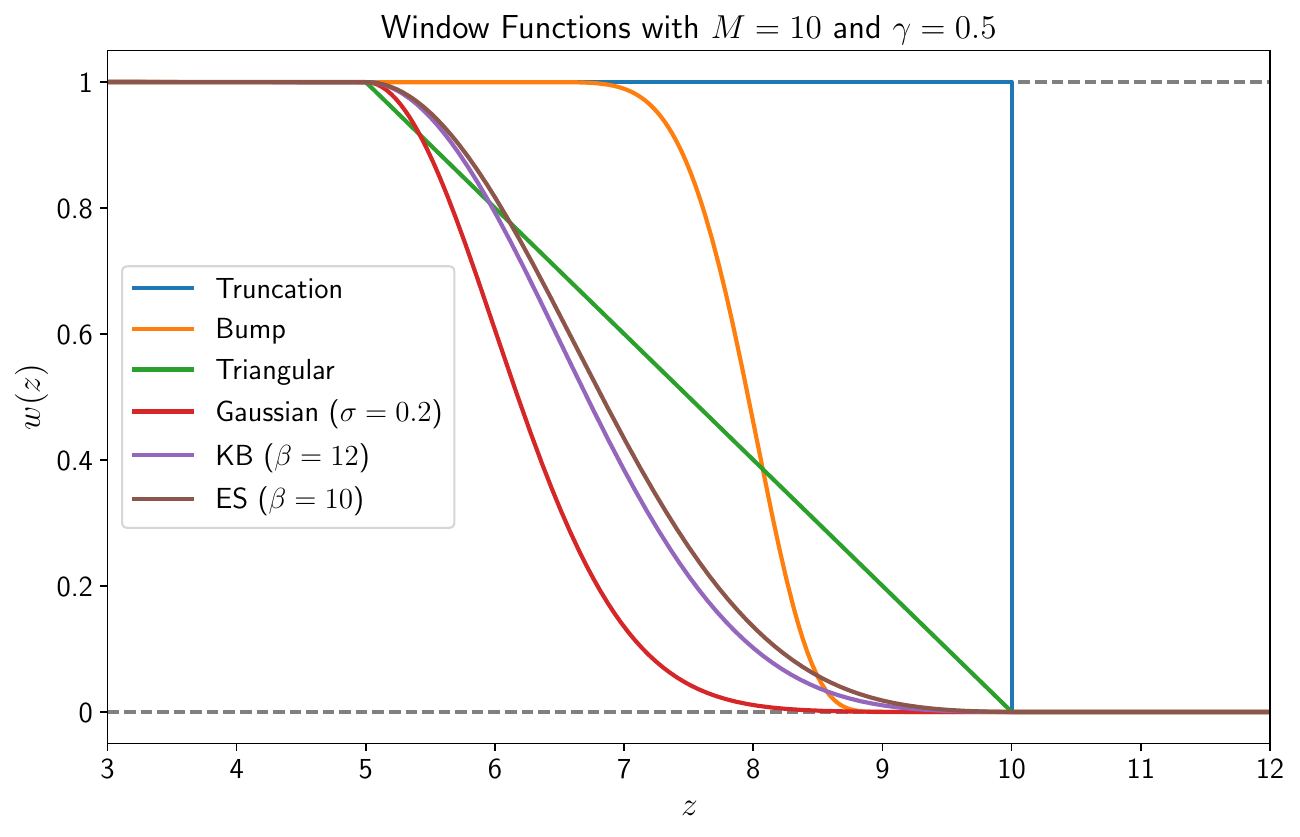}
    \hfill
    \includegraphics[width=0.485\linewidth]{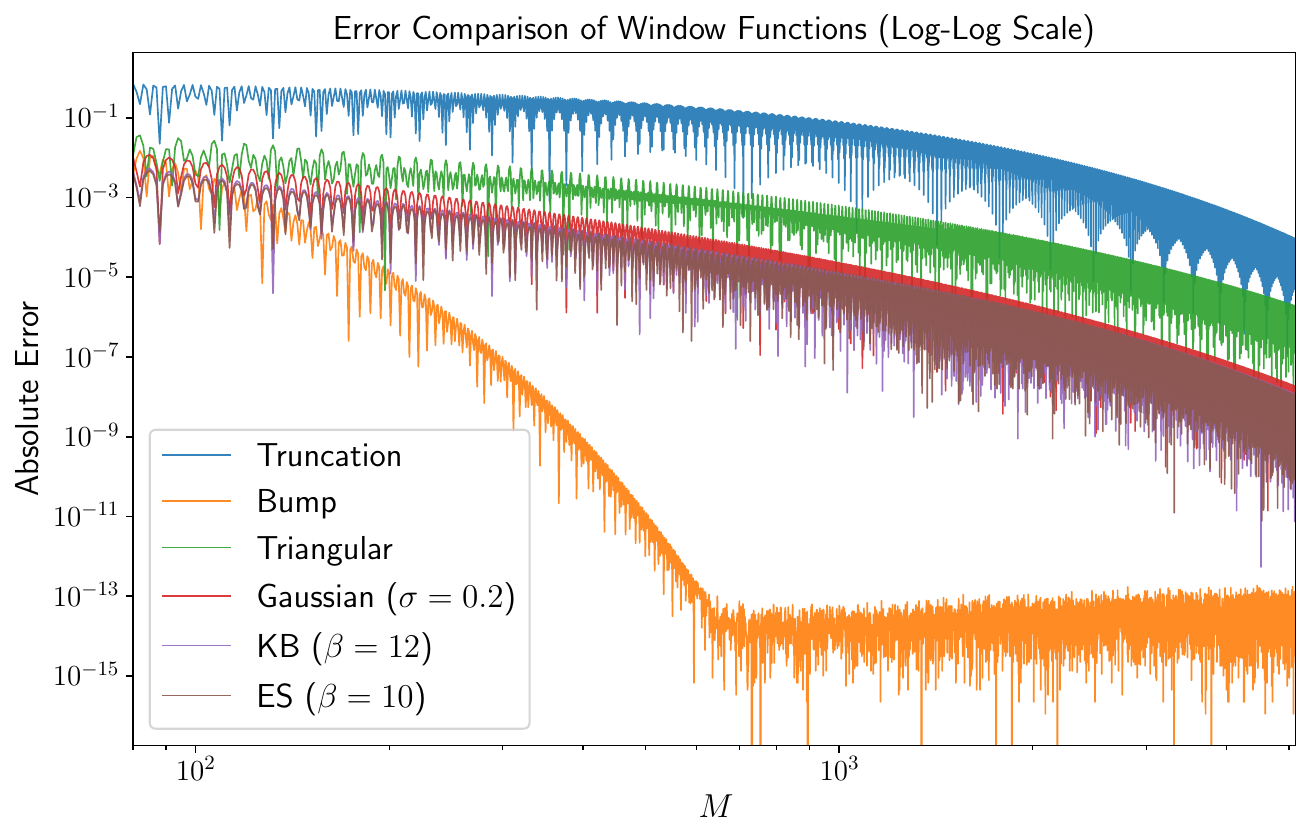}
    \hfill
    \caption{
        Left panel: Transition profiles of various window functions, plotted for $M = 10$ and $\gamma = 0.5$.
        Right panel: Log–log plot of the absolute integration error $E_{\textnormal{WF}}(M)$ for each window function, illustrating convergence behavior as $M$ increases.
    }
    \label{fig:window_function_comparison}
\end{figure}

We evaluate a variety of window functions commonly used in the spectral method literature, including the Gaussian, Kaiser-Bessel (KB), and exponential-of-semicircle (ES) windows \cite{barnett2019parallel}, as well as the classical triangular window and the bump window function introduced in \Cref{sec:far_field_integration}.
For each window, parameters such as $\sigma$ and $\beta$ are chosen to achieve near-optimal numerical performance for the given $M$ and $\gamma$ on the model problem.
All tested window functions satisfy the support conditions $w_{M}(z) = 1$ for $|z| \leq \gamma M$ and $w_{M}(z) = 0$ for $|z| \geq M$, and differ only on the transition region $(\gamma M, M)$.
The corresponding transition profiles are listed below, and representative plots illustrating their behavior are shown in the left panel of \Cref{fig:window_function_comparison}.
\begin{itemize}
    \item Truncation:
    \begin{equation*}
        w_{M, \gamma}^{\textnormal{Truncation}}(z)
        = 1.
    \end{equation*}
    \item Bump window (our choice):
    \begin{equation*}
        w_{M, \gamma}^{\textnormal{Bump}}(z)
        = \exp \left ( - 2 \frac{\exp(- 1 / s^{2})}{(1 - s)^{2}} \right ),
        \qquad
        s(z)
        = \frac{|z| - \gamma M}{M - \gamma M}.
    \end{equation*}
    \item Triangular window:
    \begin{equation*}
        w_{M, \gamma}^{\textnormal{Triangular}}(z)
        = \frac{|z|}{(\gamma - 1) M} + \frac{1}{1 - \gamma}.
    \end{equation*}
    \item Gaussian window (with $\sigma = 0.2$):
    \begin{equation*}
        w_{M, \gamma, \sigma}^{\textnormal{Gaussian}}(z)
        = \exp \left ( - \left ( \frac{|z|}{\gamma M} - 1 \right )^{2} / ( 2 \sigma^{2} ) \right ).
    \end{equation*}
    \item Kaiser-Bessel (KB) window (with $\beta = 12$):
    \begin{equation*}
        w_{M, \gamma, \beta}^{\textnormal{KB}}(z)
        = I_{0} \left ( \beta \sqrt{1 - \left ( \frac{|z|}{\gamma M} - 1 \right )^{2}} \right ) / I_{0}(\beta),
    \end{equation*}
    where $I_{0}$ is a modified Bessel function of the first kind defined as $\displaystyle I_{0}(z) = \sum_{k = 0}^{\infty} \frac{(z^{2} / 4)^{k}}{(k!)^{2}}$.
    \item Exponential of semicircle (ES) window (with $\beta = 10$)
    \begin{equation*}
        w_{M, \gamma, \beta}^{\textnormal{ES}}(z)
        = \exp \left ( \beta \left ( \sqrt{1 - \left ( \frac{|z|}{\gamma M} - 1 \right )^{2}} - 1 \right ) \right ).
    \end{equation*}
\end{itemize}
\begin{table}[htp]
    \centering
    \resizebox{0.9\linewidth}{!}{
    \begin{tabularx}{1.1\linewidth}{c|XXXXX}
        $M$ & $80$ & $160$ & $320$ & $640$ & $1280$\\
        \noalign{\hrule height 1pt}
        $E_{\textnormal{truncation}}(M)$ & \textit{7.12097e-01} & \textit{1.24460e-01} & \textit{1.56845e-01} & \textit{1.33684e-01} & \textit{4.59293e-02}\\
        $E_{\textnormal{bump}}(M)$ & \textit{3.80618e-03} & \textit{1.14729e-04} & \textit{3.26691e-08} & \textit{5.06262e-14} & \textit{4.44089e-16}\\
        $E_{\textnormal{triangular}}(M)$ & \textit{1.15552e-02} & \textit{5.77197e-03} & \textit{5.46327e-03} & \textit{6.87109e-04} & \textit{1.84488e-04}\\
        $E_{\textnormal{Gaussian}}(M)$ & \textit{1.02130e-02} & \textit{2.81350e-03} & \textit{1.25763e-04} & \textit{3.81985e-05} & \textit{8.02043e-06}\\
        $E_{\textnormal{KB}}(M)$ & \textit{4.53994e-03} & \textit{1.31051e-03} & \textit{7.23810e-05} & \textit{9.62665e-06} & \textit{1.41024e-06}\\
        $E_{\textnormal{ES}}(M)$ & \textit{3.88335e-03} & \textit{1.18340e-03} & \textit{8.16596e-05} & \textit{2.76710e-06} & \textit{7.62770e-07}\\
    \end{tabularx}
    }
    \caption{
        Absolute integration errors $E_{\textnormal{WF}}(M)$ for different window functions applied to the model problem.
        The bump window achieves near machine-precision accuracy for sufficiently large $M$, significantly outperforming other window functions, which exhibit slower convergence.
    }\label{tab:windowing_function_error_comparison}
\end{table}

\Cref{tab:windowing_function_error_comparison} and the right panel of \Cref{fig:window_function_comparison} summarize the numerical integration errors for each window function across increasing values of $M$.
As seen from the table and the plot, the bump function outperforms all other tested window functions in terms of accuracy.
In particular, for sufficiently large $M$, the bump window achieves near machine-precision errors, whereas other window functions stay at higher error levels.

In summary, while a variety of window functions from spectral analysis are available, we find that the bump window function we adopted, despite traditionally being used for analytical constructions such as partitions of unity, proves to be particularly well suited for our integration setting.

\end{document}